\documentclass[12pt]{article}
	
	\usepackage[margin=1in]{geometry}  
	\usepackage{graphicx}              
	\usepackage{amsmath}               
	\usepackage{amsfonts}              
	\usepackage{amsthm}                
	\usepackage{mathtools}
	\usepackage{xcolor}
	\usepackage{hyperref}
	\usepackage{tikz-cd}
    \usepackage{amssymb}
	
	\newtheorem{thm}{Theorem}[section]
	\newtheorem{lem}[thm]{Lemma}
	\newtheorem{prop}[thm]{Proposition}
	\newtheorem{cor}[thm]{Corollary}
	
	\newtheorem{defn}[thm]{Definition}

	\newcommand{\RR}{\mathbb{R}}      

    \newcommand{\Ga}{\mathrm{G}}
    \newcommand{\Pol}{\operatorname{Pol}}
    \newcommand{\wLambda}{\widetilde{\Lambda}}

\newcommand{\ddc}{\sqrt{-1}\partial\bar\partial}
\newcommand{\tr}{\operatorname{tr}}
\newcommand{\eps}{\epsilon}
\newcommand{\Ucal}{\Upsilon_{g_2}}
\newcommand{\Dcal}{\mathcal D_0}
\numberwithin{equation}{section}

\begin{document}

\title{Geodesics for generalised complex Hessian equations}

\author{Kartick Ghosh, XiShen Jin, Ziyue Wang and Yulun Xu}
\maketitle

\begin{abstract}
 We study geodesics associated with a class of generalized complex Hessian equations on compact K\"ahler manifolds. The defining polynomial is allowed to have degree \(1\leq d\leq n\), and is assumed to be right-Noetherian. We establish the existence of unique smooth solutions to an elliptic regularization of the geodesic equation and derive uniform a priori estimates, independent of the regularization parameter. A key ingredient is a partial Legendre transform, which reduces the regularized equation to a form that allows us to obtain uniform spatial and mixed Hessian estimates. These estimates yield the existence of unique weak $C^{1,1}$ geodesics as limits of the regularized solutions, with convergence in \(C^{1,\alpha}\) for every $0<\alpha<1$. We further establish convexity properties of associated functionals: the \(I_2\)-functional is affine and the \(J_1\)-functional is convex along weak geodesics. Finally, we study the induced \( \mathcal{E}^{p}\)-Finsler geometry on the space of admissible potentials and prove that the lengths of the regularized geodesics converge to the corresponding path distances, the \( \mathcal{E}^{p}\)-energies of the regularized geodesics converge uniformly in time to the \(p\)-th power of the corresponding path distance, and that the resulting path distance is a genuine metric for every \(1\leq p<\infty\).
\end{abstract}

\tableofcontents 

\section{Introduction}
Geodesic equations in spaces of Kähler potentials have played an important role in complex differential geometry, particularly in the study of canonical metrics, variational problems, and the geometry of spaces of solutions to nonlinear partial differential equations. A fundamental example is the space of Kähler potentials equipped with the Mabuchi metric \cite{Mabuchi}, whose geodesic equation is governed by a homogeneous complex Monge–Ampère equation on a product manifold (\cite{DonaldsonHMA}, \cite{Semmes}). The weak geodesic theory developed by X. Chen \cite{X-Chen} provides a powerful framework for understanding this infinite-dimensional geometry. X. Chen proved the existence of $C^{1,\bar{1}}$ geodesics connecting two K\"ahler potentials. The later work of Chu, Tosatti, and Weinkove \cite{Chu-Tosatti-Weicove} proves the full $C^{1,1}$ estimate.
The counterexample of Darvas and Lempert \cite{Darvas-Lempert}  proves that the $C^{1,1}$ estimate cannot be
improved to $C^{2}$. More recently, analogous geodesic constructions have been studied for various nonlinear equations arising in complex geometry ( \cite{Collins-Yau}, \cite{Chu-Collins-Lee}, \cite{CGN} ) and applied in the existence-coerciveness correspondence (\cite{Chu-Lee}, \cite{CGW}).

In this paper, we develop a geodesic theory for a class of generalized complex Hessian equations associated with lower-degree symmetric polynomials. Our setting is motivated by the generalized Monge–Ampère equations considered in (\cite{CGN}), but allows the defining polynomial to have arbitrary degree \(d\), with \(1\leq d\leq n\). The lower-degree setting introduces additional difficulties in the analysis of the geodesic equation, particularly in establishing uniform second-order estimates and constructing an appropriate elliptic regularization.

Let \((X,\omega)\) be a compact Kähler manifold of complex dimension \(n\), and let \(\chi\) be a smooth closed real \((1,1)\)-form. For a smooth real-valued function \(\varphi\), set
\[
\chi_\varphi=\chi+\sqrt{-1}\partial\bar\partial\varphi.
\]
Consider a polynomial
\begin{equation}\label{e f2 g2}
f_{2}(x)=\sum_{k=0}^{d}c_{k}^{(2)} \binom{n}{k}x^{k},\qquad g_2(A)=\sum_{k=0}^{d}c_k^{(2)}\sigma_k(A),
\qquad c_d^{(2)}>0,\quad 1\leq d\leq n,
\end{equation}
where \(\sigma_k\) denotes the \(k\)-th elementary symmetric polynomial. We assume that the corresponding one-variable polynomial $f_{2}$ is right-Noetherian (see definition \ref{right-Noetherian} ). Let \(\Upsilon_{g_2}\) denote the distinguished admissible component associated with the polarization of \(g_2\) of $f_{2}$. We consider the space of admissible potentials
\[
\mathcal H=
\left\{
\varphi\in C^\infty(X,\mathbb R):
\lambda(\omega^{-1}\chi_\varphi)\in\Upsilon_{g_2}
\right\}.
\]
On this space, we introduce the Riemannian metric
\[
\langle\psi_1,\psi_2\rangle_\varphi
=
\int_X\psi_1\psi_2\,g_2(\chi_\varphi).
\]
A direct first-variation calculation shows that the geodesic equation is
\begin{equation}
\label{introgeodesicequation}
\ddot\varphi\,g_2(\chi_\varphi)
-
\sqrt{-1}\partial\dot\varphi
\wedge\bar\partial\dot\varphi
\wedge g_2'(\chi_\varphi)=0.
\end{equation}
Here \(g_2'(\chi_\varphi)\) denotes the formal derivative of \(g_2(\chi_\varphi)\) with respect to $\chi_{\varphi}$.

As in the classical Monge–Ampère setting, this equation admits a formulation as a fully nonlinear degenerate elliptic equation in one complex dimension higher. To describe this, let
\[
\mathcal A=\{z_0\in\mathbb C:1<|z_0|<e\},
\qquad
\mathcal X=X\times\mathcal A,
\]
and consider \(S^1\)-invariant functions on \(\mathcal X\), with \(t=\log|z_0|\). In the corresponding product metric, the complex Hessian of \(\varphi\) can be represented by a block matrix
\[
H_\varphi=
\begin{pmatrix}
a&b^*\\
b&A
\end{pmatrix},
\qquad
a=\frac{\varphi_{tt}}4,\quad
b=\frac{\partial_X\varphi_t}{2},
\quad A=\chi_\varphi.
\]
Define
\[
P_0(H_\varphi)
=
a\,g_2(A)-b^*Dg_2(A)b.
\]
Then equation \ref{introgeodesicequation} is equivalent to
\[
P_0(H_\varphi)=0.
\]
This is proved in Lemma \ref{newgeodesiclemma}. This block-matrix formulation is particularly useful because it allows us to exploit the convexity and concavity properties of the admissible cone.

A key point in our analysis is the construction of a suitable elliptic regularization of the degenerate geodesic equation. We introduce the admissible block domain
\[
\mathcal D_0=
\left\{
H_\varphi:
A\in\Upsilon_{g_2},\quad
a>b^*D\log g_2(A)b
\right\}
\]
and consider the approximate equation
\begin{equation}
\label{introepsilongeodesicequation}
a^{d-1}P_0(H_{\varphi_\varepsilon})
=\varepsilon,
\qquad
\varphi_\varepsilon(\cdot,0)=\varphi_0,
\quad
\varphi_\varepsilon(\cdot,1)=\varphi_1.
\end{equation}
Solutions of \ref{introepsilongeodesicequation} are called $\epsilon$-geodesics. The main structural observation is that the operator
\[
\widetilde{G}(H)=\log P_0(H)+(d-1)\log a
\]
is concave and strictly elliptic on \(\mathcal D_0\) ( see Lemma \ref{ex:schur}). This provides the appropriate framework for applying the maximum principle and developing uniform estimates. Our main theorem is the following.
\begin{thm}\label{thm existence geo}
Let $(X,\omega)$ be a compact K\"ahler manifold and let $\chi$ be a real closed smooth $(1,1)$-form. Assume that $f_2$, which is defined by (\ref{e f2 g2}), is a right-Noetherian polynomial. Then, for any $\varphi_0, \varphi_1\in \mathcal{H}$, the following statements hold:
\begin{enumerate}
    \item For any sufficiently small $\epsilon>0$, there exists a smooth $S^1$-invariant $\epsilon$-geodesic $\varphi_{\epsilon}$ satisfying $H_{\varphi_{\epsilon}} \in \mathcal D_0$ and (\ref{introepsilongeodesicequation}) on $\overline{\mathcal{X}}=X\times\overline{\mathcal{A}}$ with $\varphi_{\epsilon}(\cdot,0)=\varphi_{0}$ and $\varphi_{\epsilon}(\cdot,1)=\varphi_{1}$. $\varphi_{\epsilon}$ satisfies 
\[ ||\varphi_{\epsilon}||_{C^1(\overline{\mathcal{X}})}+||\nabla^2\varphi_{\epsilon}||_{L^{\infty}(\overline{\mathcal{X}})}\le C\]
 for some $C$ independent of $\epsilon$ and $\varphi_{\epsilon_1}\ge\varphi_{\epsilon_2}$ for all $0<\epsilon_1\le\epsilon_2$.
    \item  There exists a unique geodesic connecting $\varphi_0$ and $\varphi_1$, satisfying
 \[ ||\varphi||_{C^1(\overline{\mathcal{X}})}+||\nabla^2\varphi||_{L^{\infty}(\overline{\mathcal{X}})}< \infty,
 \]
 $H_{\varphi}\in \overline{\mathcal D_0}$ and (\ref{introgeodesicequation}). For any $\alpha \in (0,1)$,
 \begin{equation*}
     \lim_{\epsilon \rightarrow 0} ||\varphi_{\epsilon}-\varphi||_{C^{1,\alpha}(\overline{\mathcal{X}})}=0.
 \end{equation*}
 
\end{enumerate}
\end{thm}
 The principal difficulty is obtaining a uniform bound for the full complex Hessian. We overcome this difficulty by employing a partial Legendre transform in the time variable. Kiselman's principle ensures that the minimum of a family of plurisubharmonic functions that satisfies some convexity condition is still plurisubharmonic. Such a technique was first used by Guan \cite{guan1999modified} on toric bundles. Darvas-Rubinstein \cite{darvas2016kiselman} represent geodesics in the space of K\"ahler potentials as partial Legendre transforms of envelopes.

 This estimate plays a central role in the continuity method, through which we establish the existence of a unique smooth admissible solution to the regularized equation for every sufficiently small positive \(\varepsilon\). Passing to the limit, we obtain a geodesic connecting any two prescribed admissible endpoints. Moreover, the approximate geodesics converge to the limiting geodesic in \(C^{1,\alpha}\) for every \(0<\alpha<1\) (see Proposition \ref{c1alphaconvergence}).

We next investigate the convexity properties of natural functionals along these geodesics. Given a second polynomial
\[
f_{1}(x)=\sum_{k=0}^{d}c_{k}^{(1)}\binom{n}{k} x^{k},\qquad g_1(A)=\sum_{k=0}^{d}c_k^{(1)}\sigma_k(A),
\]
satisfying the appropriate right-Noetherian and proper position assumptions relative to \(g_2\), we consider the associated functionals \(I_1\) and \(I_2\). We prove the following theorem.
\begin{thm}\label{thm convex}
    Let $f_1$ and $f_2$ be two right-Noetherian polynomials satisfying $f_2\prec f_1$ and $c_{d}^{(1)}=c_{d}^{(2)}=1$. If $I_2$ is the functional associated with $f_2$, and $J_1$ is the functional associated with $-f_1$, then $J_1(t)$ is Lipschitz and convex, and $I_2(t)$ is Lipschitz and affine along the geodesic $\varphi_t$ connecting any two potentials in $\mathcal{H}$.
\end{thm}
Finally, we study the induced \( \mathcal{E}^{p}\)-geometry. For \(1\leq p<\infty\), define the Finsler norm
\[
\|\xi\|_{p,\varphi}
=
\left(
\int_X|\xi|^p g_2(\chi_\varphi)
\right)^{1/p},
\]
and  
\[E_{p}(\varphi,t)=\int_{X}\big|\frac{d\varphi}{dt}\big|^{p}g_{2}(\chi_{\varphi}),\qquad
 \ell_p(\varphi)=\int_0^1E_{p}(\varphi,t)^{\frac{1}{p}}dt. \]
Let \(d_p\) be the corresponding path-length distance on \(\mathcal H\). With this notation, we prove the following theorem.

\begin{thm}\label{thm:metric}
For every $1\leq p<\infty$ and $\varphi_0,\varphi_1\in\mathcal H$,
the solutions of \ref{introepsilongeodesicequation} satisfy
\begin{equation}\label{eq:length-limit}
(1)\ \lim_{\epsilon\to0}\ell_p(\varphi_\epsilon)=d_p(\varphi_0,\varphi_1),
\end{equation}
\begin{equation}\label{eq:energy-limit}
 (2)\ \lim_{\epsilon\to0}\sup_{t\in[0,1]}
 \left|E_p(\varphi_\epsilon,t)-d_p(\varphi_0,\varphi_1)^p\right|=0,
\end{equation}
and
\begin{equation}\label{eq:endpoint-bound}
 (3)\ d_p(\varphi_0,\varphi_1)^p\geq\max\left\{\begin{aligned}&
  \int_{\{\varphi_0>\varphi_1\}}
       |\varphi_0-\varphi_1|^pg_{2}(\chi_{\varphi_0}),\\&
  \int_{\{\varphi_1>\varphi_0\}}
       |\varphi_1-\varphi_0|^pg_{2}(\chi_{\varphi_1})
 \end{aligned}\right\}.
\end{equation}
\end{thm}
As an immediate corollary of this theorem, we prove the following.
\begin{cor}\label{cor:metric}
For every $1\leq p<\infty$, $(\mathcal H,d_p)$ is a metric space.
\end{cor}
Our $ \mathcal{E}^{p}$-geometry is essential for the existence-coerciveness correspondence for the generalised complex Hessian equation \cite{CLN}.
The paper is organized as follows. In Section 2, we collect the algebraic and analytic properties of the admissible cone associated with the defining polynomial. In Section 3, we derive the geodesic equation and introduce its elliptic regularization. Section 4 establishes the basic \(C^0\) estimates, and develops the partial Legendre transform and proves the essential spatial Hessian estimates. We also obtain the boundary, time, and full real Hessian estimates. Section 5 establishes the existence of smooth $\epsilon$-geodesics by the method of continuity. In Section 6, we prove that the limit of $ \ epsilon$-geodesics is a geodesic. In Section 7, we prove the uniqueness of the geodesic. In Section 8, we prove Theorem \ref{thm existence geo}. In Section 9, we prove the affinity or convexity of some functionals along geodesics, which corresponds to Theorem \ref{thm convex}. In Section 10, we develop the \( \mathcal{E}^{p}\)-theory, proving the asymptotic length-minimizing property of the regularized geodesics and the metric properties of the induced distance, corresponding to Theorem \ref{thm:metric} and Corollary \ref{cor:metric}.

\noindent{\bf Tool and computational resource disclosure:}
ChatGPT 6.1 Sol was used as an adversarial reader and improve writing of the arguments. ChatGPT suggested that the new $\epsilon$-geodesic (\ref{introepsilongeodesicequation}) is useful in our case.

\section{Preliminaries}

This section recalls the required materials from
\cite{FangMaGarding}, \cite{Chao-Ming-Lin}, and \cite{LinInverseSigmaConvexityCorrigendum} in notation adapted to the present paper and derives some consequences which will be used later.

A polynomial is \emph{multi-affine} if its degree in each individual
variable is at most one.

\begin{defn}[Definition 3.1(7),(9) of \cite{FangMaGarding}]
Let \[\sigma_k(x_1,\ldots,x_n)=\sum_{i_1<\cdots<i_k}x_{i_1}\cdots x_{i_k}\] be the \(k\)-th elementary symmetric
polynomial, with \(\sigma_0=1\).

If \(h(x)=\sum_{k=0}^{d}c_k \binom nk x^k\) and \(n\ge d\), its ordinary polarization is
\[
  \Pol_n(h)(x_1,\ldots,x_n)
  =
  \sum_{k=0}^{d}c_k\sigma_k(x_1,\ldots,x_n).
\]
It is the unique symmetric multi-affine polynomial satisfying
\[
  \Pol_n(h)(x,\ldots,x)=h(x).
\]

More generally, let
\(\boldsymbol{\kappa}=(\kappa_1,\ldots,\kappa_m)\in\mathbb N^m\), and
suppose \(h(x_1,\ldots,x_m)=\sum_{\alpha\le\boldsymbol{\kappa}}
c_\alpha x^\alpha\), where
\(\alpha=(\alpha_1,\ldots,\alpha_m)\),
\(x^\alpha=\prod_i x_i^{\alpha_i}\), and
\(\alpha\le\boldsymbol{\kappa}\) means
\(\alpha_i\le\kappa_i\) for every \(i\).  The
\(\boldsymbol{\kappa}\)-polarization is
\[
  \Pi^\uparrow_{\boldsymbol{\kappa}}(h)
  =
  \sum_{\alpha\le\boldsymbol{\kappa}}c_\alpha
  \prod_{i=1}^m
  \frac{\sigma_{\alpha_i}(x_{i1},\ldots,x_{i\kappa_i})}
       {\binom{\kappa_i}{\alpha_i}}.
\]
Its output is multi-affine in all the variables
\(\{x_{ij}:1\le i\le m,\ 1\le j\le\kappa_i\}\), and setting each
block \(x_{i1}=\cdots=x_{i\kappa_i}=x_i\) recovers \(h\).
\end{defn}

\begin{defn}[{Definitions 2.2 and 2.3} of \cite{FangMaGarding}]
For \(n\ge1\), write
\begin{align*}
  \Gamma_n^+
  &=\{x=(x_1,\ldots,x_n)\in\RR^n:x_i>0\text{ for every }i\},\\
  \overline{\Gamma_n^+}
  &=\{x=(x_1,\ldots,x_n)\in\RR^n:x_i\ge0\text{ for every }i\}.
\end{align*}
A set \(\Upsilon\subseteq\RR^n\) passes the \emph{positive ray test}
\emph{(PRT)} if
\[
  \Upsilon+\overline{\Gamma_n^+}\subseteq \Upsilon.
\]
It is \emph{negative ray terminating} \emph{(NRT)} if, for
every \(x\in \Upsilon\), the set
\[
  (\{x\}-\overline{\Gamma_n^+})\cap \Upsilon
\]
is bounded.
\end{defn}

We now introduce the definition of a G{\aa}rding polynomial from \cite{FangMaGarding}. To distinguish this notion from the classical one, we call such a polynomial a Fang-Ma-G{\aa}rding polynomial. Since the original definition in Definitions~4.8 and 4.9 of \cite{FangMaGarding} is more complicated, we use another characterization here, which is equivalent to the original definition by Theorem~1.1(3) of \cite{FangMaGarding}.

\begin{defn}[Theorem 1.1(3) of
\cite{FangMaGarding}]
\label{FMGdefinition}
Let \(p\in\RR[x_1,\ldots,x_n]\) be nonzero.  For a
multi-index
\(\alpha=(\alpha_1,\ldots,\alpha_n)\in\mathbb Z_{\ge0}^n\), write
\[
  \partial^\alpha p
  =
  \frac{\partial^{|\alpha|}p}
       {\partial x_1^{\alpha_1}\cdots\partial x_n^{\alpha_n}},
  \qquad |\alpha|=\alpha_1+\cdots+\alpha_n.
\]
Then \(p\) is called a \emph{Fang-Ma-G{\aa}rding polynomial} if and only if,
for every \(\alpha\) such that \(\partial^\alpha p\not\equiv0\),
\begin{enumerate}
  \item the set \(\{\partial^\alpha p>0\}\) has a unique connected
  component \(\Upsilon_{\partial^\alpha p}\) passing PRT; and
  \item this component is contained in the corresponding components
  of all its first partial derivatives:
  \[
    \Upsilon_{\partial^\alpha p}
    \subseteq \Upsilon_{\partial_i\partial^\alpha p}
    \qquad(1\le i\le n).
  \]
\end{enumerate}
Here \(\Upsilon_h=\RR^n\) when \(h\equiv c>0\), and by convention \(\Upsilon_0=\RR^n\).
The component \(\Upsilon_p\) is called the \emph{Fang-Ma-G{\aa}rding component} of
\(p\), which is a generalization of the $\Upsilon$-stable component of \cite{Chao-Ming-Lin}.  We denote the class of such polynomials in \(n\) variables by
\(\Ga_n\), and put \(\Ga=\bigcup_{n\ge1}\Ga_n\).

\end{defn}

\begin{thm}[{Theorem 8.13} of \cite{FangMaGarding}]
\label{kappa-polarization}
Let \(p\in\Ga_n\) have multidegree at most
\(\boldsymbol{\kappa}\), meaning
\(\deg_{x_i}p\le\kappa_i\) for every \(i\).  If its Fang-Ma-G{\aa}rding
component \(\Upsilon_p\) is NRT, then
\[
  \Pi^\uparrow_{\boldsymbol{\kappa}}(p)\in\Ga.
\]
\end{thm}

Fang and Ma have provided the following characterization of NRT polynomials:

\begin{lem}[Lemma 9.5 of \cite{FangMaGarding}]
\label{NRT-redundant}
Let $f\in \Ga_n$. Then the Fang--Ma--G{\aa}rding component $\Upsilon_f$ is NRT if and only if $\deg_{x_i} f\ge 1$ for every $i=1,2,\ldots,n$.
\end{lem}

We next recall the definition introduced by Lin.

\begin{defn}[Definition 2.1 of \cite{Chao-Ming-Lin}]
\label{right-Noetherian}
For a nonconstant univariate polynomial \(h\), let \(r(h)\) denote its
largest real root, provided such a root exists.  If
\(\deg h=d\), its \emph{root sequence} is
\[
  R(h)=\bigl(r(h),r(h'),\ldots,r(h^{(d-1)})\bigr),
\]
provided all entries exist.  The polynomial \(h\) is
\emph{right-Noetherian} if all these entries exist
and
\[
  r(h)\ge r(h')\ge\cdots\ge r(h^{(d-1)}).
\]

The polynomial \(h\) with degree $d\ge 2$ is
\emph{strictly right-Noetherian} if all these entries exist
and
\[
  r(h)> r(h')\ge\cdots\ge r(h^{(d-1)}).
\]

\end{defn}

The following is immediate from Section~9.2 of \cite{FangMaGarding}.
\begin{lem}[Section 9.2 of \cite{FangMaGarding}]
\label{FangMaLinCorrespondence}
A monic univariate polynomial $p(x)$ belongs to $\Ga_1$ if and only if it is right-Noetherian.
\end{lem}

Next, let us recall the definition of proper position from \cite{FangMaGarding}.

\begin{defn}[{Definition 10.9} of \cite{FangMaGarding}]
 Suppose $f_2(x), f_1(x) \in \Ga[x] \setminus \{0\}$. The pair $(f_2,f_1)$ is said to be \emph{in proper position} if
 \[f_2(x)y +f_1(x) \in \Ga[x,y].\]
In this case, we also write $f_2 \prec f_1$.
\end{defn}

Using this definition and Theorem \ref{kappa-polarization}, we can prove the following corollary. 
        \begin{cor}
        \label{polarizationproperposition}
        Suppose that $f_2, f_1\in \Ga_1$ are univariate monic polynomials of degree $1\le d\le n$. If $f_2\prec f_1$, then $\Pol_n f_2 \prec \Pol_n f_1$.
        \end{cor}
        \begin{proof}
        By definition, $f_2(x)y+f_1(x)\in \Ga[x,y]$. By Lemma \ref{NRT-redundant}, its Fang-Ma-G{\aa}rding component is NRT. By Theorem \ref{kappa-polarization}, if we choose $\kappa=(n,1)$, then the $\kappa$-polarization 
        \[
        g_2(\mathbf{x})y+g_1(\mathbf{x})\in \Ga[\mathbf{x},y],
        \]
        where $g_2=\Pol_n f_2$ and $g_1=\Pol_n f_1$. Thus, $g_2\prec g_1$ by definition.
        \end{proof}
    The main reason for us to use the proper position condition is the following:
    \begin{thm}[Theorem 6.14 of \cite{FangMaGarding}]
\label{PositiveDerivative}
Let $g_2$ and $g_1$ be two nontrivial multi-affine Fang--Ma--G{\aa}rding polynomials. Then the following are equivalent:
\begin{enumerate}
    \item $g_2\prec g_1$;
    \item $\Upsilon_{g_2} \cap\{g_1>0\}=\Upsilon_{g_1}$;
    \item For every $i=1,2,\ldots,n$, $g_1\partial_i g_2-g_2\partial_i g_1 \le 0$ on $\Upsilon_{g_2}$.
\end{enumerate}
\end{thm}

\begin{lem}\label{ex:ideal}
The Fang-Ma-G\aa rding component $\Ucal$ is convex and stable under addition of positive
semidefinite matrices, i.e. if \(A \in \Upsilon_{g_2}\) and \(B \succeq 0\), then \(A + B \in \Upsilon_{g_2}\). Moreover, we have $g_2(A+B)\geq g_2(A)$. On $\Ucal$, $Dg_2(A)>0$, and $g_2^{1/d}$ and $\log g_2$ are
concave. Its finite boundary lies in $\{g_2=0\}$.
\end{lem}
\begin{proof}
Theorem~8.2 of \cite{FangMaIdealGarding} explicitly includes $d\leq n$ and makes the
polarization $g_2$ ideal G\aa rding. Its Fang-Ma-G\aa rding component $\Ucal$ and the
components of its nontrivial derivatives w.r.t. eigenvalues $\lambda_i, 1\leq i\leq n$, $\Upsilon_{\partial^{\alpha} g_2}$ are convex, pass the positive ray
test, and satisfy derivative-component inclusion $\Ucal \subset \cap^{n}_{i=1} \Upsilon_{\partial^{i} g_2} $ . Each first partial $\partial^{i} g_2>0$ because $c_d^{(2)}>0$ and $d\geq1$. Hence, $Dg_2(A)$ i.e. differentiating w.r.t. elements of $A$ is positive definite on $\Ucal$.
Theorem~12.1 of \cite{FangMaIdealGarding} gives the concavity of $g_2^{1/d}$ and $\log g_2$ w.r.t. eigenvalues $\lambda_i$. To lift these assertions to matrices, use the Ky Fan--Lidskii inequality\cite{Horn_Johnson_2012}:
the eigenvalues of a convex combination are majorized by the corresponding
convex combination of ordered eigenvalues, i.e. $\lambda(\alpha A + (1 - \alpha)B) \prec \alpha \lambda(A) + (1 - \alpha)\lambda(B)$ for $0\leq \alpha \leq 1$. By \cite{BONDAR}, since $g_2^{1/d}$ and $\log g_2$ are symmetric and concave, they are schur-concave. Hence, 
\[\log g_2(\lambda(\alpha A + (1 - \alpha)B))\geq \log g_2(\alpha \lambda(A) + (1 - \alpha)\lambda(B))\geq \alpha\log g_2( \lambda(A))+(1 - \alpha)\log g_2( \lambda(B))
\]
This proves the concavity w.r.t. matrix.

Then we want to prove that $\Ucal$ is stable under addition of positive
semidefinite matrices, which means that if \(A \in \Upsilon_{g_2}\) and \(B \succeq 0\), i.e. positive semidefinite, then \(A + B \in \Upsilon_{g_2}\).
Firstly, because \(B \succeq 0\), by Weyl's monotonicity theorem, the eigenvalues satisfy:
    \[
        \lambda_i(A+B) \ge \lambda_i(A), \quad \forall i
    \]
    This shows that \(\lambda(A)\) is component-wise dominated by \(\lambda(A+B)\), i.e., \(\lambda(A+B) - \lambda(A) \in \mathbb{R}^n_{\ge 0}\).
    Because \(\Upsilon_{g_2}\) satisfies PRT in the eigenvalue space, i.e., \(\Upsilon_{g_2} + \mathbb{R}^n_{\ge 0} \subset \Upsilon_{g_2}\).
    Therefore, \(\lambda(A) + (\lambda(A+B) - \lambda(A)) = \lambda(A+B) \in \Upsilon_{g_2}\).
    Thus, \(A + B \in \Upsilon_{g_2}\). Moreover, $g_2(A+B)\geq g_2(A)$ follows from $Dg_2>0$ on $\Ucal$.

A finite boundary point with $g_2>0$
would have a connected positive neighborhood meeting the same component,
which is impossible.
\end{proof}
\begin{lem}\label{ex:center}
There exists a fixed translation of $\chi$ and $g_2$, leaving $P_0(H_\varphi)$ unchanged and
reduces to $c_{d-1}^{(2)}=0$. After that, on $\Ucal$, 
\begin{equation}\label{ex:trace}
 \tr A>0,\qquad
 -(n-1)\lambda_{\max}(A)<\lambda_i(A)\leq\lambda_{\max}(A).
\end{equation}
If $A_\ell\in\Ucal$, $R_\ell\to +\infty$, and $A_\ell/R_\ell\to B$,
then $B\in\overline\Gamma_d$.
\end{lem}
\begin{proof}
Set
\[
 h=\frac{c_{d-1}^{(2)}}{(n-d+1)c_d^{(2)}},\qquad
 \tilde\chi=\chi+h\omega,\qquad \tilde g_2(A)=g_2(A-hI).
\]
The $\sigma_{d-1}$ coefficient of $\tilde g_2$ vanishes. Its scalar
polynomial is $f_2(r-h)$, so right-Noetherianity is preserved. For the original block matrix \(H_\phi = \begin{pmatrix} a & b^* \\ b & A \end{pmatrix}\), where \(A = \chi_\phi = \chi + \sqrt{-1}\partial\bar{\partial}\phi\).
Under the new \(\tilde{\chi}\), the corresponding block matrix becomes:
\[
\tilde{H}_\phi = \begin{pmatrix} a & b^* \\ b & \tilde{A} \end{pmatrix}, \quad \text{where } \tilde{A} = \tilde{\chi}_\phi = \chi_\phi + h\omega = A + hI.
\]
while $a$ is unchanged. Both $g_2$ and $Dg_2$ are invariant under the translation,
and hence $P_0(H_\varphi)$ is unchanged.

Since \(c_{d-1}^{(2)} = 0\), the polynomial is mainly dominated by the \(d\)-th order terms.
For the monomial \(\sigma_d(A)\), taking partial derivatives with respect to \(d-1\) distinct variables leaves a 1st-order elementary symmetric polynomial, i.e. the trace:
\[
\frac{\partial^{d-1} \sigma_d(A)}{\partial x_{i_1} \cdots \partial x_{i_{d-1}}} = \sum_{j \notin \{i_1, \dots, i_{d-1}\}} x_j = \operatorname{tr} A - \sum_{k \in \{i_1, \dots, i_{d-1}\}} x_k
\]

On $\Ucal$, the sum of these derivatives equals \(c_d^{(2)}\) multiplied by the sum of all complementary eigenvalues. Since there are a total of \(\binom{n-1}{d-1}\) choices, the total sum is:
\[
c_d^{(2)}\sum_{|I|=d-1} \partial_I \sigma_d(\lambda) =c_d^{(2)} \sum_{j=1}^{n} \binom{n-1}{d-1} \lambda_j =  
c_d^{(2)} \binom{n-1}{d-1} \operatorname{tr} A
\]
Because \(g_2\) is within the Gårding cone \(\Upsilon_{g_2}\), and \(c_{d-1}^{(2)} = 0\), the derivative is always positive, so we must have \(\operatorname{tr} A > 0\).
Since \(\operatorname{tr} A > 0\), the maximum eigenvalue \(\lambda_{\max}(A)\) must be \(> 0\).
If some eigenvalue \(\lambda_i(A) \le -(n-1)\lambda_{\max}(A)\), then the sum \(\operatorname{tr} A = \lambda_i + \sum_{j \neq i} \lambda_j \le -(n-1)\lambda_{\max} + (n-1)\lambda_{\max} = 0\), which contradicts \(\operatorname{tr} A > 0\). Therefore:
\[
-(n-1)\lambda_{\max}(A) < \lambda_i(A) \le \lambda_{\max}(A)
\]

The derivative $D_I^r g_2$, $0\leq r<d$,
has leading term
$c_d^{(2)}(n-d+r)!\,\sigma_{d-r}/(n-d)!$, i.e.
\[
D_I^r g_2(\lambda) = c_d^{(2)} \frac{(n - d + r)!}{(n - d)!} \sigma_{d-r}(\lambda) + \text{lower order terms (from } \sigma_k \text{ with } k < d \text{)}
\] 
Since any order of derivative of \(g_2\) on the positive components $\Ucal$ is strictly greater than \(0\), we have \(D_I^r g_2(A_\ell) \ge 0\).
Divide its value at $A_\ell$ by $R_\ell^{d-r}$ and pass to the limit. We have
\[
    \lim_{\ell \to \infty} \frac{D_I^r g_2(A_\ell)}{R_\ell^{d-r}} = \lim_{\ell \to \infty} c_d^{(2)} \frac{(n - d + r)!}{(n - d)!} \sigma_{d-r}\left(\frac{A_\ell}{R_\ell}\right)\ge 0
    \] 
Since \(A_\ell / R_\ell \to B\), the limit equals:
    \[
    c_d^{(2)} \frac{(n - d + r)!}{(n - d)!} \sigma_{d-r}(B) \ge 0
    \]
Thus $\sigma_j(B)\geq0$ for $1\leq j\leq d$, i.e. $B\in\overline\Gamma_d$. 
\end{proof}

We work with the centered data henceforth. The potentials and their
equation are unchanged by this fixed normalization.
We retain the symbols $f_2,g_2,c_k^{(2)},\widetilde c_k,\chi$ for the
translated data, with $\widetilde c_k=c_{k-1}^{(2)}$ as above.

\begin{lem}\label{ex:schur}
Suppose $A\in \Ucal$, one has
\begin{equation}
     g_2(A-\tau bb^*/a)
 =(1-\tau)g_2(A)+\tau P_0(H)/a>0,
 \qquad0\leq\tau\leq1.
\end{equation}
on $\Dcal$, which is defined by
\begin{equation}\label{ex:schur-id}
 \Dcal=\{H:a>0,\ S=A-bb^*/a\in\Ucal\},\qquad P_0(H)=a g_2(S).
\end{equation}
The domain $\Dcal$ is convex, and
$G=\log P_0+(d-1)\log a$ is concave and strictly elliptic there.
\end{lem}
\begin{proof}
For rank-one matrix, we have $det(I-tB)=1-t \rm{tr}(B)$. Hence, let $M:=bb^{*}$, by Cauchy’s formula for the determinant of a rank-one perturbation \cite{Horn_Johnson_2012}:
\begin{align*}
\det(I + t(A - M)) 
&= \det(I + tA) - t b^{*} \operatorname{adj}(I + tA)b.
\end{align*}
Now, expand both sides as a power series in \(t\). The left-hand side is:
\[
\sum_{k=0}^{n} \sigma_k(A - M)t^k
\]
The right-hand side is:
\[
\sum_{k=0}^{n} \sigma_k(A)t^k - t b^{*} \operatorname{adj}(I + tA)b
\]
$g_2$ is defined by
\[
g_2(A - M) = \sum_{k=0}^{n} c^{(2)}_k \sigma_k(A - M)
\]
For $H\in\Dcal$, since $g_2(A-\tau bb^*/a)$ is linear in $\tau$ by comparing equation:
\[
 g_2(A-\tau bb^*/a)
 =g_2(A)-\tau b^{*} Dg_2(A)b/a 
 =(1-\tau)g_2(A)+\tau P_0(H)/a>0,
 \qquad0\leq\tau\leq1.
\]
Since $A\in\Ucal$, $P_0(H)=a g_2(S)>0$, the path $\tau \mapsto A-\tau bb^*/a\in \Ucal, 0\leq\tau\leq1$. The Schur map $A-bb^*/a$ is Loewner
concave\cite{Horn_Johnson_2012}, i.e.let $H_t=(1-t)H_0+tH_1,A_t=(1-t)A_0+tA_1,a_t=(1-t)a_0+ta_1,b_t=(1-t)b_0+tb_1,S_t=A_t-b_tb_t^*/a_t,S_t=(1-t)S_0+tS_1$, we have $S_t-(1-t)S_0-tS_1$ is semipositive. By \ref{ex:ideal}, monotonicity and convexity of $\Ucal$ prove convexity of
$\Dcal$; monotonicity and concavity of $\log g_2$ prove concavity of
$\log P_0$. Hence, $G=\log P_0+(d-1)\log a$ is concave. With $Q=D\log g_2(S)$ and $\xi=b/a$,
\[
 D\log P_0(H)=
 \begin{pmatrix}a^{-1}+\xi^*Q\xi&-\xi^*Q\\-Q\xi&Q\end{pmatrix}.
\]
Since this is a triangular congruence of a  $\operatorname{diag}(a^{-1},Q)$ using $\begin{pmatrix}1&0\\\xi&I\end{pmatrix}$.
By \ref{ex:ideal}, $Dg_2$ is positive definite. Combined with $g_2>0$ on $\Ucal$, we know $Q$ and $D\log P_0(H)$ are positive definite. 
Adding $(d-1)\log a$ preserves strict ellipticity.
\end{proof}

\begin{lem}
\label{bound_mat}
    Suppose \[|B| \le M\] and  the unscaled matrix \(R^2 B\in \Ucal\), we have \[
\sigma_k(B) \ge -O(R^{-2})= -o(1), 1\leq k\leq d.
\]
\end{lem}
\begin{proof}
Let \(B\) be the matrix of \(\alpha\) in a \(\beta\)-unitary frame. Since
 \(R_j^2 B_j \in \Ucal\), derivative-component inclusion gives  $R_j^2 B_j \in \Upsilon_{\partial_\alpha g_2} $ for any $\alpha$ such that $\partial_\alpha g_2 \not\equiv 0$. Therefore, we have
\[
D_I^{d-k} g_2(R^2 B) > 0,
\]
for $k \le d-1$, where $D_I$ is the differentiation in the identity direction .
Using \(D_I \sigma_\ell = (n-\ell+1)\sigma_{\ell-1}\), we have
\[
D_I^{d-k} g_2(A) = c_d^{(2)} \frac{(n-k)!}{(n-d)!} \sigma_k(A) + \text{terms of degree at most } k-1.
\]

Each remaining term is a fixed multiple of some \(\sigma_q(A), q < k\). Its rescaled value satisfies
\[
R^{-2k} \sigma_q(R^2 B_j) = R^{-2(k-q)} \sigma_q(B) = O(R^{-2}),
\]
because \(|B| \le M\) bounds every \(\sigma_q(B)\) uniformly.

Consequently,
\[
0 < R^{-2k} D_I^{d-k} g_2(R^2 B) = c_d^{(2)} \frac{(n-k)!}{(n-d)!} \sigma_k(\alpha; \beta) + O(R^{-2}).
\]

Since the leading coefficient is positive,
\[
\sigma_k(B) \ge -O(R^{-2})= -o(1).
\]
\end{proof}

\section{The Geodesic equation}
For $i=1,2$, set $c_d^{(i)}=1$. We define
	   \[
	   f_i(x) = \sum_{k=0}^{d} c_k^{(i)} \binom{n}{k} x^k, \quad
	   g_i(\lambda) = \sum_{k=0}^{d} c_k^{(i)} \sigma_k(\lambda),
	   \]
	   \[
		   g_i(\chi_\varphi) = \sum_{k=0}^{d} c_k^{(i)}\binom{n}{k} \chi_\varphi^k \wedge \omega^{n-k},\quad i=1,2.
		   \]
		   Assume that $f_1, f_2$ are right-Noetherian polynomials such that $f_2\prec f_1$.
  
        Let $\mathcal{H} = \{ \varphi \in C^{\infty}(X,\mathbb{R})\mid \lambda(\omega^{-1}\chi_\varphi)\in \Upsilon_{g_2} \}$. The tangent space at $\varphi\in\mathcal{H}$ is $T_{\varphi}\mathcal{H}=C^{\infty}(X,\mathbb{R})$. Consider the following Riemannian metric:
        \[ \langle \Psi_1, \Psi_2 \rangle_{\varphi} = \int_{X} \Psi_1 \Psi_2 g_2(\chi_\varphi)
        =\int_{X} \Psi_1 \Psi_2 \sum_{k=0}^{d} \binom{n}{k} c_k^{(2)} \chi_{\varphi}^{k}\wedge  \omega^{n-k}.\]
        Let $\varphi(t)$ be a smooth curve in $\mathcal{H}$ with $\varphi(0)=\varphi_{0}$ and $\varphi(1)=\varphi_{1}$. Let us also consider the family of curves $\varphi(t,s)$ in $\mathcal{H}$ such that $\varphi(0,s)=\varphi_{0},\varphi(1,s)=\varphi_{1}$ and $\varphi(t,0)=\varphi(t)$. 
        Consider the energy
        \[E(s) = \int_{0}^{1} \int_{X} \left( \frac{d\varphi}{dt} \right)^{2}\sum_{k=0}^{d} \binom{n}{k} c_k^{(2)} \chi_{\varphi}^{k}\wedge  \omega^{n-k}dt.\]
        We now differentiate with respect to $s$.
        
        \begin{align*}
            \frac{dE}{ds}\Big|_{s=0} & = \int_0^1 \int_{X} 2\frac{d\varphi}{dt}\frac{d}{dt}\left(\frac{d}{ds}\Big|_{s=0}\varphi\right)\sum_{k=0}^{d} \binom{n}{k} c_k^{(2)} \chi_{\varphi}^{k}\wedge  \omega^{n-k}dt\\
            & +\int_{0}^{1}\int_{X}\big(\frac{d\varphi}{dt}\big)^{2}\sum_{k=1}^{d}\binom{n}{k}c_k^{(2)}k\chi_{\varphi}^{k-1}\wedge\omega^{n-k}\wedge\sqrt{-1}\partial\bar{\partial}\left(\frac{d}{ds}\Big|_{s=0}\varphi\right)dt
        \end{align*}
        
        Integrating the second term by parts, we obtain
        \begin{align*}
            & = \int_0^1 \int_{X} 2\frac{d\varphi}{dt}\frac{d}{dt}\left(\frac{d}{ds}\Big|_{s=0}\varphi\right)\sum_{k=0}^{d} \binom{n}{k} c_k^{(2)} \chi_{\varphi}^{k}\wedge  \omega^{n-k}dt\\
            &+2\int_{0}^{1}\int_{X}(\frac{d}{ds}\Big|_{s=0}\varphi)\sqrt{-1}\partial(\frac{d\varphi}{dt})\wedge \bar{\partial}(\frac{d\varphi}{dt})\wedge\sum_{k=1}^{d}\binom{n}{k}c_k^{(2)}k\chi_{\varphi}^{k-1}\wedge\omega^{n-k}dt\\
            &+2\int_{0}^{1}\int_{X}(\frac{d}{ds}\Big|_{s=0}\varphi)(\frac{d\varphi}{dt}) \sqrt{-1}\partial\bar{\partial}(\frac{d\varphi}{dt})\wedge\sum_{k=1}^{d}\binom{n}{k}c_k^{(2)}k\chi_{\varphi}^{k-1}\wedge\omega^{n-k}dt
        \end{align*}
        Integrating the first term by parts in time, we obtain
        \begin{align*}
            &\int_0^1 \int_{X} 2\frac{d\varphi}{dt}\frac{d}{dt}\left(\frac{d}{ds}\Big|_{s=0}\varphi\right)\sum_{k=0}^{d} \binom{n}{k} c_k^{(2)} \chi_{\varphi}^{k}\wedge  \omega^{n-k}dt\\
            &=-\int_{0}^{1}\int_{X}2\frac{d^{2}\varphi}{dt^{2}}\left(\frac{d}{ds}\Bigg|_{s=0}\varphi\right)\sum_{k=0}^{d} \binom{n}{k} c_k^{(2)} \chi_{\varphi}^{k}\wedge  \omega^{n-k}dt\\
            &-\int_{0}^{1}\int_{X}2\frac{d\varphi}{dt}\left(\frac{d}{ds}\Bigg|_{s=0}\varphi\right)\sum_{k=1}^{d}\binom{n}{k}c_k^{(2)}k\chi_{\varphi}^{k-1}\wedge\omega^{n-k}\wedge\sqrt{-1}\partial\bar{\partial}\left(\frac{d\varphi}{dt}\right)dt.
        \end{align*}
        Therefore,
        \begin{align*}
            \frac{dE}{ds}\Bigg|_{s=0} &=2\int_{0}^{1}\int_{X}(\frac{d}{ds}\Big|_{s=0}\varphi)\sqrt{-1}\partial(\frac{d\varphi}{dt})\wedge \bar{\partial}(\frac{d\varphi}{dt})\wedge\sum_{k=1}^{d}\binom{n}{k}c_k^{(2)}k\chi_{\varphi}^{k-1}\wedge\omega^{n-k}dt\\
            &-\int_{0}^{1}\int_{X}2\frac{d^{2}\varphi}{dt^{2}}\left(\frac{d}{ds}\Bigg|_{s=0}\varphi\right)\sum_{k=0}^{d} \binom{n}{k} c_k^{(2)} \chi_{\varphi}^{k}\wedge  \omega^{n-k}dt
        \end{align*}

        Thus, we obtain the geodesic equation
        \begin{equation}
                \frac{d^2\varphi}{dt^2} \sum_{k=0}^{d}c_k^{(2)}\binom{n}{k}\chi_{\varphi}^k\wedge\omega^{n-k} 
                -\sqrt{-1}\partial\frac{d\varphi}{dt}\wedge  \bar{\partial}(\frac{d\varphi}{dt})\wedge 
                \sum_{k=1}^{d}kc_k^{(2)}\binom{n}{k}\chi_{\varphi}^{k-1}\wedge\omega^{n-k} =0.
            \label{GeodesicEq2}        
        \end{equation}
    
        \subsection{Geodesic equation in one dimension higher and regularization}
Let us write 
\[g_2(A)=\sum_{j=0}^d c_j^{(2)}\sigma_j(A)\]
for a Hermitian matrix $A$.
We see that $g_2(\chi_\varphi)$ is an $(n,n)$-form, whereas $g_2(A)$ is
the scalar polynomial evaluated on a Hermitian matrix. In an
$\omega$-unitary frame, $A$ is the coefficient matrix of
$\chi_\varphi$ and $g_2(\chi_\varphi)=g_2(A)\omega^n$.
Throughout the matrix calculations, spatial forms are identified with
these matrices. Matrix derivatives such as $Dg_2$ refer to the formal derivative of the 
polynomial $g_{2}(\chi_{\varphi})$ with respect to $\chi_{\varphi}$. In particular, the fixed endpoints
$\varphi_0,\varphi_1\in\mathcal H$ have matrices
$A_j=\chi_{\varphi_j}\in\Ucal$, $j=0,1$.

As in \cite{CGN}, let
\[
 \mathcal A=\{z_0\in\mathbb C:1<|z_0|<e\},\qquad
 \mathcal X=X\times\mathcal A,\qquad t=\log|z_0|\in(0,1).
\]
Potentials are $S^1$-invariant. We use $\partial,\bar\partial$ on $X$
and $D,\bar D$ on $\mathcal X$. Time derivatives are  as
$\partial_t\varphi$ and $\partial_t^2\varphi$. Fix the product metric
\[
 \Omega=\pi_X^*\omega+
          \frac{\sqrt{-1}\,dz_0\wedge d\bar z_0}{|z_0|^2}.
\]
In the logarithmic time frame, we write
\begin{equation}\label{ex:blocks}
\begin{gathered}
 H_\varphi=\pi_X^*\chi+\sqrt{-1}D\bar D \varphi=
 \begin{pmatrix}a&b^*\\b&A\end{pmatrix},\qquad
 a=\frac{\partial_t^2\varphi}4,\quad b=\frac{\partial_X \partial_t\varphi}2,\quad A=\chi_\varphi,\\
 P_0(H_{\varphi})=a g_2(A)-b^{*}Dg_2(A)b,
\end{gathered}
\end{equation}
where $b^{*}Dg_2(A)b\ \omega^{n}$ corresponds to the form $\frac{1}{4}\sqrt{-1}\partial\frac{d\varphi}{dt}\wedge  \bar{\partial}(\frac{d\varphi}{dt})\wedge 
               g_{2}'(\chi_{\varphi})$.
\begin{lem}
\label{newgeodesiclemma}
    With the above notations, we have that the geodesic equation \ref{GeodesicEq2} can be written as $P_{0}(H_{\varphi})=0$. 
\end{lem}
\begin{proof}
    We calculate 
    \begin{align*}
        &P_{0}(H_{\varphi})\omega^{n}\\
        &=ag_{2}(A)\omega^{n}-b^{*}Dg_{2}(A)b\omega^{n}\\
        &=\frac{1}{4}\frac{d^{2}\varphi}{dt^{2}}g_{2}(\chi_{\varphi})-\frac{1}{4}\sqrt{-1}\partial\frac{d\varphi}{dt}\wedge  \bar{\partial}(\frac{d\varphi}{dt})\wedge 
               g_{2}'(\chi_{\varphi}) \\
                &=\frac{1}{4}\times (\text{ L.H.S of geodesic equation}). 
    \end{align*}
\end{proof}
We have seen in Lemma \ref{ex:schur} that $\Dcal=\{H_{u}: a>0, S=A-\frac{bb^{*}}{a}\in\Upsilon_{g_{2}}\}$ and $\mathcal{F}=\log P_{0}+(d-1)\log a$ is concave and strictly elliptic on $\Dcal$. Hence we have the following definition which is an elliptic  regularization of the geodesic equation. Observe that it is different from the corresponding regularization in \cite{CGN}.  
\begin{defn}\label{ex:admissibility}
The admissible block domain and the $\epsilon$-geodesic are
\[
 \Dcal=\{H:A\in\Ucal,\quad a>b^{*}D\log g_2(A)b\},
\]
\begin{equation}\label{ex:family}
 (*_\eps)\qquad a^{d-1}P_0(H_{\varphi_\eps})=\eps,\qquad H_{\varphi_\eps}\in\Dcal,
 \qquad \varphi_\eps(\cdot,j)=\varphi_j\quad(j=0,1).
\end{equation}
\end{defn}

When estimating one fixed approximation, we suppress its subscript
and write $\varphi$ for $\varphi_\eps$.

\section{Hessian estimate for $\epsilon$-geodesic}
Denote $F(A)=g_2(A)^{1/d}$ and $G(u) \triangleq F(\omega^{-1} \chi_u)+ 4 \epsilon^{1/d} u_{pp}$. In this section, we write $\varphi$ for $\epsilon$-geodesic.
\subsection{$C^0$ estimate}

\begin{lem}\label{exs:barriers}
There are fixed $\eps_0,K_0,C>0$ and a smooth function
\[
 \underline{\varphi}=(1-t)\varphi_0+t\varphi_1-K_0t(1-t)
\]
such that every solution $\varphi$ of \eqref{exs:family}, $0<\eps\leq\eps_0$,
satisfies
\begin{equation}\label{exs:comparison}
 \underline{\varphi}\leq \varphi\leq(1-t)\varphi_0+t\varphi_1,
 \qquad \|\varphi\|_\infty+\|\partial_t\varphi\|_\infty\leq C.
\end{equation}

\end{lem}
\begin{proof}
The spatial matrices $(1-t)\omega^{-1}\chi_{\varphi_0}+t \omega^{-1} \chi_{\varphi_1}$ form a compact subset of $\Ucal$.
 The mixed space-time derivatives of $H_{\underline{\varphi}}$ are fixed, whereas its
time entry is $K_0/2$. Thus one fixes a large $K_0$ so that $H_{\underline{\varphi}}\in \Dcal$. Lemma~\ref{ex:schur} then gives a fixed positive lower bound $s_0$
for $(\frac{\partial_{t}^2 \underline{\varphi}}{4})^{d-1}P_0 (H_{\underline{\varphi}})$. Take $\eps_0=s_0$. Comparison principle of the elliptic operator
$\log P_0+(d-1)\log (\frac{\partial_t^2 \varphi}{4})$ gives $\underline{\varphi}\leq \varphi$.
Admissibility implies $\partial_t^2\varphi>0$, so time convexity gives the upper
bound in \eqref{exs:comparison}. 
The endpoint derivative inequalities and time convexity bound $\partial_t\varphi$
throughout the cylinder.
\end{proof}

\subsection{Partial Legendre transform}

On the domain
\[
 \mathcal O=\{(x,p):\partial_t\varphi(x,0)<p<\partial_t\varphi(x,1)\},
\]
use the partial Legendre transform and its time variables
\begin{equation}\label{exs:dual-variables}
 \begin{gathered}
 p=\partial_t\varphi(x,t),\qquad \psi(x,p)=\inf_{s\in [0,1]} (\varphi(x,s)-ps)=  \varphi(x,t)-pt,\\
 -\psi_p=t,\qquad -\psi_{pp}=\frac1{\partial_t^2\varphi},
 \end{gathered}
\end{equation}

We extend $\varphi_0$ on $\mathcal O$ by $\varphi_0(x,p)= \varphi_0(x)$.

The regularity of $\psi$ is provided by the following Lemma:
\begin{lem}[Darvas--Rubinstein, \cite{darvas2019minimum} Lemma~3.8]
\label{lem:darvas-rubinstein-minimum}
Let $I\subset\mathbb R$ be an open interval,
$D\subset\mathbb R^n$ a domain, and
$f\in C^2(I\times D)$.
Assume that, for every $x\in D$, the function
$t\mapsto f(t,x)$ is strongly convex and attains
its minimum at $t(x)\in I$.
Define
\[
 g(x)=\inf_{t\in I}f(t,x)=f(t(x),x).
\]
Then $g\in C^2(D)$ and
\begin{equation}
\label{eq:darvas-rubinstein-minimum-hessian}
 \nabla^2 g(x)
 =\left(
     \nabla_x^2 f
     -\frac{1}{f_{tt}}
       (\nabla_{tx}^2 f)^{\mathsf T}\nabla_{tx}^2 f
   \right)(t(x),x).
\end{equation}
\end{lem}

\begin{lem}\label{lem psi hess}
  Fixing $p$, we have that:  $\psi_{i\bar{j}}=\varphi_{i\bar{j}} - \frac{\varphi_{t i} \varphi_{t \bar j}}{\varphi_{tt}}$
\end{lem}
\begin{proof}
    Denote $t=t(x,p)$. Then, we can write 
    \begin{equation*}
        \psi(x,p)= \varphi(x,t(x,p))-p t(x,p).
    \end{equation*}
    Take the derivative of the above formula in the space direction. We can get:
    \begin{equation*}
        \psi_i = \varphi_{i} + \varphi_t t_i -p t_i= \varphi_i.
    \end{equation*}
    Take another derivative of the above formula in the space direction. We can get:
    \begin{equation}\label{e psi ij}
        \psi_{i\bar{j}} = \varphi_{i \bar j}+ \varphi_{i t }t_{\bar j}.
    \end{equation}
    Take space derivative on  the formula $p= \partial_t \varphi(x,t(x,p))$. We get:
    \begin{equation}\label{e varphi ti tt}
         0 = \varphi_{t \bar{j}} + \varphi_{tt} t_{\bar j}.
    \end{equation}
    Plugging (\ref{e varphi ti tt}) into (\ref{e psi ij}), we can conclude the proof of this Lemma.
\end{proof}

\begin{lem}\label{e equ psi1}
  Fixing $p$, we have that:  $P_0(H_{\varphi})=\frac{\partial^2_t \varphi}{4} g_2(\omega^{-1}\chi_{\psi})$.
\end{lem}
\begin{proof}
Take a normal coordinate of $\omega$. Denote $A= (\chi_{\varphi})_{i \bar{j}}$ and $b=\frac{1}{2}( \varphi_{1 t}, \varphi_{2 t}, ... , \varphi_{n t})$.
    Using the Lemma \ref{ex:schur}, we can get that:
    \begin{equation*}
        g_2(A- \frac{4 b b^*}{\partial_t^2 \varphi})=g_2(A)-\frac{4}{\partial_t^2 \varphi}b^* D g_2(A) b.
    \end{equation*}
    Combining this formula with Lemma \ref{lem psi hess}, we can get:
    \begin{equation*}
        \frac{\partial_t^2 \varphi}{4} g_2(\omega^{-1}\chi_{\psi})= P_0(H_{\varphi}).
    \end{equation*}
\end{proof}

\begin{lem}\label{e equ psi}
  $\psi$ satisfies  
  \begin{equation}\label{equ psi}
      G(\psi)=0.
  \end{equation}
\end{lem}
\begin{proof}
By Lemma \ref{e equ psi1}, we can get: $P_0(H_{\varphi})=a g_2(\omega^{-1}\chi_{\psi})$. Plugging this into (\ref{exs:family}), we can get:
\begin{equation*}
    (\frac{\partial_t^2 \varphi}{4})^d g_2(\omega^{-1}\chi_{\psi})= \epsilon.
\end{equation*}
Since $ \frac{\varphi_{tt}}{4}=-\frac{1}{4\psi_{pp}}$, we can get that:
\begin{equation*}
g_2(\omega^{-1}\chi_{\psi})= (-4 \psi_{pp})^d \epsilon.
\end{equation*}
Take the power of $1/d$ on both sides of the above formula. We can conclude the proof of this Lemma.
\end{proof}

Define the linearized operator of (\ref{equ psi}) as $\mathcal L u= F^{i \bar j}(\omega^{-1}\chi_{\psi}) u_{i \bar j}+ 4 \epsilon^{1/d} u_{pp}$.  Denote $\mathcal F = F^{i\bar j}\omega_{i \bar j}$. Then, we have:
\begin{lem}\label{lem lpsi phi0}
  There exist constants $\eta, m_0 >0$ such that  $$ \mathcal L \psi_p=0, \,\,\, -\mathcal L (\psi- \varphi_0) \ge \eta \mathcal F +m_0$$
\end{lem}
\begin{proof}
    Taking derivative with respect to $p$ in (\ref{equ psi}), we can get $\mathcal L \psi_p =0$. 
    \begin{equation}\label{e lpsi varphi0}
        -\mathcal L (\psi- \varphi_0)= -F^{i \bar j}(\omega^{-1}\chi_{\psi})(\psi- \varphi_0)_{i \bar j} -4 \epsilon^{1/d}\psi_{pp}=-F^{i \bar j}(\omega^{-1}\chi_{\psi})(\psi- \varphi_0)_{i \bar j} + F(\omega^{-1}\chi_{\psi}).
    \end{equation}
    Using the concavity of  $F=g_2^{1/d}$ on $\Ucal$ by Lemma \ref{ex:ideal}, we can get that:
    \begin{equation}\label{e conc f}
        F(\omega^{-1}\chi_{\psi}) + F^{i \bar{j}}(\omega^{-1}\chi_{\psi}) \big( (\varphi_0)_{i \bar{j}} -\eta \omega - \psi_{i\bar j} \big) \ge F\big(\omega^{-1}(\chi_{\varphi_0}-\eta \omega)\big) \ge m_0 >0.
    \end{equation}
    Here $\eta$ is a small positive constant such that $\chi_{\varphi_0}-\eta \omega \in \Upsilon_{g_2}$. $m_0= \inf F\big(\omega^{-1}(\chi_{\varphi_0}-\eta \omega) \big)$. Combining (\ref{e lpsi varphi0}) and (\ref{e conc f}), we can prove  $ -\mathcal L (\psi- \varphi_0) \ge \eta \mathcal F +m_0$.
\end{proof}

\begin{lem}\label{lem math f low}
    For any $A\in \Upsilon_{g_2}$, we have that:
    \begin{equation*}
        \mathcal F(A) \ge (c_d^{(2)} \binom nd )^{1/d}.
    \end{equation*}
\end{lem}
\begin{proof}
    For any $t \ge 0$, $A+ t I \in \Upsilon_{g_2}$. Then, we can use the concavity of $F$ to get:
    \begin{equation*}
        F(A+ t I) \le F(A) + F^{i \bar j} (A+ t I -A)_{i \bar{j}} .
    \end{equation*}
    Divide both sides of the above formula by $t$ and let $t$ go to $\infty$, we can use the fact that
    \begin{equation*}
        \lim_{t \rightarrow \infty}\frac{F(A+ t I)- F(A)}{t}= (c_d^{(2)} \binom nd )^{1/d}
    \end{equation*}
    to conclude the proof of this Lemma.
\end{proof}

Let $\lambda_1 \ge \lambda_2 \ge \ldots \ge \lambda_n$ be eigenvalues of $\chi_{\psi}$ with respect to $\omega$. In a coordinate where $\chi_{i \bar j}= \lambda_i \delta_{ij}$, we write $F_i = F^{i \bar i}$. We write $g_{i \bar j}= (\chi_{\psi})_{i \bar j}.$

\begin{lem}
    Take $K=1+ \sup_{\mathcal O} |\partial_X (\psi-\varphi_0)|_{\omega}^2$. Then, we have:
    \begin{equation}\label{e lpartial z}
        \mathcal L |\partial (\psi-\varphi_0)|_{\omega}^2 \ge \frac{1}{2} \sum_i F_i \lambda_i^2 + \frac{1}{2} \sum_{i,j} F^{i \bar i} |(\psi- \varphi_0)_{i j}|^2 -CK \mathcal F.
    \end{equation}
    
    At a point where $\lambda_1$ is smooth, we have that:
    \begin{equation}\label{e llam1}
       \mathcal L \lambda_1 \ge -F^{i \bar{j}, k \bar l} g_{i \bar j1} g_{k \bar l \bar 1} -C \lambda_1 \mathcal F
    \end{equation}
    and
    \begin{equation}\label{e concave est}
        \frac{-F^{i\bar j, k \bar l} g_{i \bar {j} 1} g_{k \bar l \bar 1}}{\lambda_1} \ge \sum_{F^{i \bar i} > 2 F^{1 \bar 1} } \frac{\lambda_1 (F^{i \bar i}- F^{1 \bar 1})}{\lambda_1-\lambda_i} |(\log \lambda_1)_i|^2.
    \end{equation}
\end{lem}
\begin{proof}
First, we estimate (\ref{e lpartial z}).
Comparing (98) of \cite{szekelyhidi2018fully} with our case, the additional term we have is $4 \epsilon^{1/d}(|\partial (\psi-\varphi_0)|_{\omega}^2)_{pp}$. 
\begin{equation*}
    4 \epsilon^{1/d}(|\partial (\psi-\varphi_0)|_{\omega}^2)_{pp} = 8 \epsilon^{1/d} |\partial (\psi-\varphi_0)_p|^2 + 8 \epsilon^{1/d} Re (\psi- \varphi_0)_i (\psi- \varphi_0)_{\bar i pp}.
\end{equation*}
We can combine the term $ 8 \epsilon^{1/d} Re (\psi- \varphi_0)_i (\psi- \varphi_0)_{\bar i p p}$  with the second order term $2 Re F^{k \bar l} (\psi- \varphi_0)_i (\psi- \varphi_0)_{k \bar l  \bar i}$ from $F^{i \bar j}(|\partial (\psi- \varphi_0)|_{\omega}^2)_{i \bar j}$, and use the first order derivative of (\ref{equ psi}) to deal with them.
Then, we can use (98) of \cite{szekelyhidi2018fully} to get (\ref{e lpartial z}).

(\ref{e llam1}) follows from (81) of \cite{szekelyhidi2018fully}. Since we assume that $\omega$ is K\"ahler, the torsion term there disappears in our case.

(\ref{e concave est}) is \cite[(67)]{szekelyhidi2018fully}.
\end{proof}

In the calculations below,
$g_{i\bar j}$ denotes the components of $\chi_\psi$.

\begin{prop}\label{exs:quadratic-spatial}
There are $\kappa_0,C>0$, depending only on the $M, \varphi_0, \varphi_1, \omega, \chi$, such that every $\epsilon$-geodesic above with
$0<\epsilon \leq \epsilon_0$ satisfies
\begin{equation}\label{exs:quadratic-bound}
 \|\chi_\psi\|_{L^{\infty}(\mathcal O)} \leq 
 C\left(1+\sup_{\mathcal O}|\partial_X(\psi-\varphi_0)|_\omega^2\right).
\end{equation}
\end{prop}
\begin{proof}
The proof of this Proposition is similar to \cite[Proposition 13]{szekelyhidi2018fully}. The difference is that in our case there are terms which are derivatives in $p$ direction. We will emphasize how to deal with those terms.
 Take $K=1+ \sup_{\mathcal O} |\partial_X (\psi-\varphi_0)|_{\omega}^2$. consider
\[
 \begin{gathered}
 \Phi=\log\lambda_1+H(|\partial(\psi-\varphi_0)|_\omega^2)+\xi(\psi-\varphi_0)+\gamma(-\psi_p),\\
 H(r)=-\frac1{12}\log(2K-r),\qquad
 \xi(r)=-Ar+\frac1{r+||\psi-\varphi_0||_{\infty}+1},\qquad
 \gamma(r)=-\frac1{12}\log(2-r).
 \end{gathered}
\]
Then $(24K)^{-1}\leq H'\leq(12K)^{-1}$,
$H''=12(H')^2$, $-A-1\leq\xi'\leq-A$,
$\xi''\geq c_{\xi}=2/(2||\psi-\varphi_0||_{\infty}+1)^3$, and
$\gamma''=12(\gamma')^2$. 

Let $(x_0,p_0)\in \mathcal O$ be a maximum point of $\Phi$. We do the calculation in a normal coordinate of $\omega$ around $(x_0,p_0).$ Define $B=diag(0,\epsilon_1,\ldots, \epsilon_1)$. Then the biggest eigenvalue $\widetilde{\lambda}_1$ of $(g_{i \bar j})-B$ is locally smooth near $(x_0,p_0)$ and it is equal to $\lambda_1$ at $(x_0,p_0)$. As a result, without loss of generality, we can assume that $\lambda_1$ is smooth near $(x_0,p_0)$. First, we use \ref{e llam1} to get
\begin{align}
\mathcal L\log\lambda_1
 &=\frac{F^{i\bar j}(\lambda_1)_{i\bar j}}{\lambda_1}
   -\frac{F^{i\bar j}(\lambda_1)_i(\lambda_1)_{\bar j}}
          {\lambda_1^2}
   +4\epsilon^{1/d}
    \left(\frac{(\lambda_1)_{pp}}{\lambda_1}
          -\frac{(\lambda_1)_p^2}{\lambda_1^2}\right)
   \notag\\
 &=\frac{L\lambda_1}{\lambda_1}
   -\frac{F^{i\bar j}(\lambda_1)_i(\lambda_1)_{\bar j}}
          {\lambda_1^2}
   -\frac{4\epsilon^{1/d}(\lambda_1)_p^2}{\lambda_1^2}
   \notag\\
 &\geq
   -\frac{F^{p\bar q,r\bar s}
           g_{p\bar q 1}g_{r\bar s\bar 1}}{\lambda_1}
   -C_0\mathcal F
   -\frac{F^{k\bar k}|g_{1\bar 1 k}|^2}{\lambda_1^2}
   -\frac{4\epsilon^{1/d}(\lambda_1)_p^2}{\lambda_1^2}.
 \label{draft:log-eigenvalue}
\end{align}
Differentiating
\[
 F(\chi_\psi)+4\epsilon^{1/d}\psi_{pp}=0
\]
gives
\begin{equation}
 F^{k\bar k}g_{k\bar k\ell}
       +4\epsilon^{1/d}\psi_{pp\ell}=0.
 \label{draft:differentiated-equation}
\end{equation}
For the gradient term, the chain rule gives
\begin{align}
\mathcal L H\bigl(|\partial(\psi-\varphi_0)|_\omega^2\bigr)
 &=H' \mathcal L\bigl(|\partial(\psi-\varphi_0)|_\omega^2\bigr)
   \notag\\
 &\quad+H''F^{i\bar j}
   \bigl(|\partial(\psi-\varphi_0)|_\omega^2\bigr)_i
   \bigl(|\partial(\psi-\varphi_0)|_\omega^2\bigr)_{\bar j}
   \notag\\
 &\quad+4\epsilon^{1/d}H''
   \bigl(|\partial(\psi-\varphi_0)|_\omega^2\bigr)_p^2.
 \label{draft:gradient-chain}
\end{align}
Using (\ref{e lpartial z}), we obtain
\begin{align}
\mathcal L H\bigl(|\partial(\psi-\varphi_0)|_\omega^2\bigr)
 &\geq \frac12 H'\sum_i F^{i\bar i}\lambda_i^2 
       -CH'K\mathcal F + \frac{1}{2} H' \sum_{i,j} F^{i \bar i} |(\psi- \varphi_0)_{i j}|^2
   \notag\\
 &\quad+H''F^{i\bar j}
   \bigl(|\partial(\psi-\varphi_0)|_\omega^2\bigr)_i
   \bigl(|\partial(\psi-\varphi_0)|_\omega^2\bigr)_{\bar j}
   \notag\\
 &\quad+4\epsilon^{1/d}H''
   \bigl(|\partial(\psi-\varphi_0)|_\omega^2\bigr)_p^2.
 \label{draft:gradient-weight-bound}
\end{align}
For the next term,
\begin{align}
\mathcal L\xi(\psi-\varphi_0)
 &=\xi' L(\psi-\varphi_0)
   +\xi''F^{i\bar j}(\psi-\varphi_0)_i
                         (\psi-\varphi_0)_{\bar j}
   \notag\\
 &\quad+4\epsilon^{1/d}\xi''(\psi-\varphi_0)_p^2.
 \label{draft:xi-chain}
\end{align}
Using Lemma \ref{lem lpsi phi0},
\begin{align}
\mathcal L\xi(\psi-\varphi_0)
 &\geq -\xi' (\eta \mathcal F + m_0)+\xi''F^{i\bar j}(\psi-\varphi_0)_i
                             (\psi-\varphi_0)_{\bar j}
   \notag\\
 &\quad
       +4\epsilon^{1/d}\xi''(\psi-\varphi_0)_p^2.
 \label{draft:xi-bound}
\end{align}

Finally, using Lemma \ref{lem lpsi phi0}, $\gamma'>0$, and $\gamma''>0$,
\begin{align}
\mathcal L\bigl(\gamma(-\psi_p)\bigr)
 &=\gamma' L(-\psi_p)
   +\gamma''F^{i\bar j}\psi_{pi}\psi_{p\bar j}
   +4\epsilon^{1/d}\gamma''\psi_{pp}^2
   \notag\\
 &=\gamma''F^{i\bar j}\psi_{pi}\psi_{p\bar j}
   +4\epsilon^{1/d}\gamma''\psi_{pp}^2.
 \label{draft:gamma-chain}
\end{align}

Combining the above calculations gives
\begin{align}
\mathcal L\Phi
 &\geq
   -\frac{F^{p\bar q,r\bar s}
           g_{p\bar q 1}g_{r\bar s\bar 1}}{\lambda_1}
   -C_0\mathcal F
   -\frac{F^{k\bar k}|g_{1\bar 1 k}|^2}{\lambda_1^2}
   -\frac{4\epsilon^{1/d}(\lambda_1)_p^2}{\lambda_1^2}
   \notag\\
 &\quad+\frac12 H'\sum_i F^{i\bar i}\lambda_i^2
   +\frac12 H'\sum_{i,j}F^{i\bar i}
                       |(\psi-\varphi_0)_{ij}|^2
   \notag\\
 &\quad+H''F^{i\bar j}
   \bigl(|\partial(\psi-\varphi_0)|_\omega^2\bigr)_i
   \bigl(|\partial(\psi-\varphi_0)|_\omega^2\bigr)_{\bar j}
   \notag\\
 &\quad+4\epsilon^{1/d}H''
   \bigl(|\partial(\psi-\varphi_0)|_\omega^2\bigr)_p^2
   -\xi'(\eta \mathcal F +m_0)
   \notag\\
 &\quad+\xi''F^{i\bar j}(\psi-\varphi_0)_i
                             (\psi-\varphi_0)_{\bar j}
   +4\epsilon^{1/d}\xi''(\psi-\varphi_0)_p^2
   \notag\\
 &\quad+\gamma''F^{i\bar j}\psi_{pi}\psi_{p\bar j}
   +4\epsilon^{1/d}\gamma''\psi_{pp}^2.
 \label{draft:combined}
\end{align}

According to the calculation in \cite[Proposition 13]{szekelyhidi2018fully}, the key point is to estimate the third-order derivatives $-\frac{F^{k\bar k}|g_{1\bar 1 k}|^2}{\lambda_1^2}
   -\frac{4\epsilon^{1/d}(\lambda_1)_p^2}{\lambda_1^2}$.

\medskip
\noindent\textbf{Case 1: $-\lambda_n>\delta\lambda_1$.}
At a maximum point, $\nabla\Phi=0$ implies
\begin{align}
 0
 &=\frac{g_{1\bar 1 k}}{\lambda_1}
   +H'\sum_i\Bigl(
       (\psi-\varphi_0)_i(\psi-\varphi_0)_{\bar i k}
       +(\psi-\varphi_0)_{ik}(\psi-\varphi_0)_{\bar i}
     \Bigr)
   \notag\\
 &\quad+\xi'(\psi-\varphi_0)_k-\gamma'\psi_{pk}.
 \label{draft:spatial-stationarity}
\end{align}
Hence
\begin{align}
 -\frac{F^{k\bar k}|g_{1\bar 1 k}|^2}{\lambda_1^2}
 &=-F^{k\bar k}\Biggl|
       \gamma'\psi_{pk}-\xi'(\psi-\varphi_0)_k
   \notag\\[-0.3em]
 &\hspace{5em}
       -H'\sum_i\Bigl(
         (\psi-\varphi_0)_i(\psi-\varphi_0)_{\bar i k}
         +(\psi-\varphi_0)_{ik}(\psi-\varphi_0)_{\bar i}
       \Bigr)\Biggr|^2
   \notag\\
 &\geq-3F^{k\bar k}(\gamma')^2|\psi_{pk}|^2
       -3F^{k\bar k}|\xi'|^2|(\psi-\varphi_0)_k|^2
   \notag\\
 &\quad-3F^{k\bar k}(H')^2
       \left|\sum_i\Bigl(
         (\psi-\varphi_0)_i(\psi-\varphi_0)_{\bar i k}
         +(\psi-\varphi_0)_{ik}(\psi-\varphi_0)_{\bar i}
       \Bigr)\right|^2.
 \label{draft:three-term-square}
\end{align}
The first term is absorbed by $\gamma''F^{i\bar j}\psi_{pi}\psi_{p\bar j}$ in (\ref{draft:combined}), and
the third term is absorbed by
 $$H''F^{i\bar j}
   \bigl(|\partial(\psi-\varphi_0)|_\omega^2\bigr)_i
   \bigl(|\partial(\psi-\varphi_0)|_\omega^2\bigr)_{\bar j}=12 (H')^2 F^{i\bar j}
   \bigl(|\partial(\psi-\varphi_0)|_\omega^2\bigr)_i
   \bigl(|\partial(\psi-\varphi_0)|_\omega^2\bigr)_{\bar j}$$ in (\ref{draft:combined})
 The middle term is bounded by $C A^2 K\mathcal F$.

Differentiating $\Phi$ with respect to $p$ similarly gives
\begin{align}
 0
 &=\frac{g_{1\bar 1 p}}{\lambda_1}
   +H'\sum_i\Bigl(
       (\psi-\varphi_0)_i(\psi-\varphi_0)_{\bar i p}
       +(\psi-\varphi_0)_{ip}(\psi-\varphi_0)_{\bar i}
     \Bigr)
   \notag\\
 &\quad+\xi'(\psi-\varphi_0)_p-\gamma'\psi_{pp}.
 \label{draft:p-stationarity}
\end{align}

Using (\ref{draft:p-stationarity}), we obtain
\begin{align}
 -\frac{4\epsilon^{1/d}(\lambda_1)_p^2}{\lambda_1^2}
 &=-4\epsilon^{1/d}\Biggl|
       \gamma'\psi_{pp}-\xi'(\psi-\varphi_0)_p
   \notag\\[-0.3em]
 &\hspace{5em}
       -H'\sum_i\Bigl(
         (\psi-\varphi_0)_i(\psi-\varphi_0)_{\bar i p}
         +(\psi-\varphi_0)_{ip}(\psi-\varphi_0)_{\bar i}
       \Bigr)\Biggr|^2
   \notag\\
 &\geq-12\epsilon^{1/d}(\gamma')^2\psi_{pp}^2
       -12\epsilon^{1/d}(\xi')^2(\psi-\varphi_0)_p^2
   \notag\\
 &\quad-12\epsilon^{1/d}(H')^2
       \left|\sum_i\Bigl(
         (\psi-\varphi_0)_i(\psi-\varphi_0)_{\bar i p}
         +(\psi-\varphi_0)_{ip}(\psi-\varphi_0)_{\bar i}
       \Bigr)\right|^2.
 \label{draft:p-three-term-square}
\end{align}
The first term is absorbed by $4\epsilon^{1/d}\gamma''\psi_{pp}^2$ in (\ref{draft:combined}), and the second term is bounded because $|(\psi-\varphi_0)_p|= |\psi_p|\le 1$,
and the last term is absorbed by
 $$4\epsilon^{1/d}H''
   \bigl(|\partial(\psi-\varphi_0)|_\omega^2\bigr)_p^2=48\epsilon^{1/d}(H')^2
   \bigl(|\partial(\psi-\varphi_0)|_\omega^2\bigr)_p^2$$ in (\ref{draft:combined}). The rest of the proof in the Case 1 is the same as in \cite{szekelyhidi2018fully}

\medskip
\noindent\textbf{Case 2: $-\lambda_n\leq\delta\lambda_1$.} Using \cite{szekelyhidi2018fully}'s argument, $-\sum_k\frac{F^{k\bar k}|g_{1\bar1 k}|^{2}}{\lambda_1^2}$ can be decomposed into two parts. One part is estimated using the argument in the Case 1. Another part is estimated using (\ref{e concave est}).

The estimate of
\[
 -\frac{4\epsilon^{1/d}(\lambda_1)_p^2}{\lambda_1^2}
\]
is the same as in Case~1. The key point is that the coefficient of the $\partial^2_p$ term of $L$ is uniformly bounded. The rest of the proof is the same as in \cite{szekelyhidi2018fully}.

In conclusion, if $(x_0,p_0)$ is an interior point of $\mathcal O$, then the above computation implies that:
\begin{equation*}
    \Phi(x_0,p_0)\le CK.
\end{equation*}
Recall that Lemma \ref{exs:barriers} provides the bound on $||\psi||_{\infty}$. This implies (\ref{exs:quadratic-bound}). 

If $(x_0,p_0)\in \partial \mathcal O$, we can use Lemma \ref{lem psi hess} to get $\chi_{\psi}\le \chi_{\varphi}$. This implies that:
\begin{equation*}
   \lambda_1 (x_0,p_0) \le C K
\end{equation*}
or 
\begin{equation*}
   \chi_{\psi} (x_0,p_0) \le \chi_{\varphi_1}(x_0).
\end{equation*}
This shows that $\Phi(x_0,p_0)$ is uniformly bounded, which also proves (\ref{exs:quadratic-bound}). 
\end{proof}

Extend the partial Legendre transform of a smooth admissible solution by
\[
 \psi(x,p)=\min_{0\leq t\leq1}\{\varphi(x,t)-pt\},\qquad p\in\mathbb R.
\]
Extend $t$ as well by:
\[
t(x,p)=
\begin{cases}
-\psi_p(x,p), & (x,p)\in\mathcal O,\\[2pt]
0, & p\leq \partial_t\varphi(x,0),\\[2pt]
1, & p\geq \partial_t\varphi(x,1).
\end{cases}
\]
\begin{lem}\label{lem t lip}
    $t(x,p)$ is a Lipschitz continuous function on $X\times \mathbb R$.
\end{lem}
\begin{proof}
    Since $\partial_t^2\varphi>0$ on $X\times [0,1]$, we can denote
\[
 \mu:=\min_{X\times[0,1]}\varphi_{tt}>0.
\]

Take two arbitrary points $(x,p)$ and $(y,q)$.
Suppose first that $t(x,p)>t(y,q)$. Then $t(x,p)>0$ and $t(y,q)<1$,
this implies
\[
 p\geq\varphi_t(x,t(x,p)),
 \qquad
 q\leq\varphi_t(y,t(y,q)).
\]
Equality holds when the corresponding minimizer is interior;
the inequalities also cover $t(x,p)=1$ and $t(y,q)=0$.
Consequently,
\[
 \begin{aligned}
 \mu(t(x,p)-t(y,q))
 &\leq\varphi_t(x,t(x,p))-\varphi_t(x,t(y,q))\\
 &\leq p-q+\varphi_t(y,t(y,q))-\varphi_t(x,t(y,q))\\
 &\leq |p-q|+M|x-y|.
 \end{aligned}
\]
Interchanging the two points handles $t(y,q)>t(x,p)$, while
$t(x,p)=t(y,q)$ is immediate. Therefore,
\begin{equation}\label{eq:clipped-inverse-lipschitz}
 |t(x,p)-t(y,q)|
 \leq\frac{|p-q|+M|x-y|}{\mu}.
\end{equation}
Thus the clipped inverse is jointly locally Lipschitz in
$(x,p)$, including across the contact graphs.
The constants $\mu^{-1}$ and $M$ may depend on the fixed
solution; no uniformity in $\epsilon$ is asserted.

The clipped inverse $t(x,p)$ is therefore locally
Lipschitz. 
\end{proof}

\begin{lem}\label{exg:contact-extension}
$\psi \in C_{loc}^{1,1}(X \times \mathbb{R})$. 
Then the full extension satisfies
\begin{equation}\label{exg:extended-spatial-bound}
 |\chi_\psi|_\omega\leq C(1+\sup|\partial_X \psi|_\omega^2)
 \quad \text{almost everywhere}.
\end{equation}
\end{lem}
\begin{proof}
Write $p_j(x)=\partial_t\varphi(x,j)$. The minimizing time is $0$ for $p\leq p_0(x)$,
$1$ for $p\geq p_1(x)$, and is determined by $p=\partial_t\varphi(x,t)$ between these
graphs. By Lemma \ref{lem t lip}, $t$ is Lipschitz continuous. Then, we can apply the identities
\begin{equation}\label{e psi varphi spa}
 \psi_p=-t(x,p),\qquad
 \partial_X\psi=\partial_X \varphi(x,t(x,p))
\end{equation}
to get that $\psi \in C^{1,1}(\overline{\mathcal O})$ regularity. On the contacts, $\psi=\varphi_0$ or $\varphi_1-p$;
its first derivatives match on $\partial \mathcal O$. Hence the spatial
distributional Hessian equals its bounded, piecewise Hessian and has no
interface measure. Proposition~\ref{exs:quadratic-spatial} implies
(\ref{exg:extended-spatial-bound}).
\end{proof}

\begin{lem}\label{res:5-5}
There exist constants $C, \epsilon_0>0$ such that for any $0 < \epsilon < \epsilon_0$, the corresponding $\epsilon$-geodesic $\psi$, after being extended to $X\times \mathbb R$ satisfies:
\begin{equation}\label{eq:148}
\|\partial_X \psi\|_{C^0(X\times \mathbb R)}\le C.
\end{equation}
\end{lem}
\begin{proof}
\textbf{Step 1:} Blow up of $\psi$. 
Suppose that there exists a sequence of $\epsilon$-geodesics $\{\varphi_j\}$ whose partial Legendre transform $\psi_j$ 
satisfy
\[
K_j=\|\partial_X \psi_j\|_{C^0(X\times \mathbb R)}\longrightarrow\infty,
\qquad |\partial_X \psi_j(x_j, p_j)|_\omega=K_j.
\]

Note that $\mathcal O \subset \{(x,p): \min_X \partial_t \varphi (x,0 < \max_X \partial_t \varphi(x,1))\}$. This implies that
\begin{equation}
    \begin{split}
        \psi(x,p) &= \varphi(x,1)-p \text{ for } p > \max_X \partial_t \varphi(x,1) \\
        \psi(x,p) &= \varphi(x,0) \text{ for } p< \min \partial_t \varphi(x,0).
    \end{split}
\end{equation}
As a result, we may assume that $ \inf_{\epsilon< \epsilon_0} \min_{x\in X} \partial_t \varphi(x,0) \le p_j \le \sup_{\epsilon< \epsilon_0} \max_{x\in X} \partial_t \varphi(x,1)$. Using Lemma \ref{exs:barriers}, we can see that $ \inf_{\epsilon< \epsilon_0} \min_{x\in X} \partial_t \varphi(x,0)$ and $\sup_{\epsilon< \epsilon_0} \max_{x\in X} \partial_t \varphi(x,1)$ are uniformly bounded. Thus, up to a subsequence, we can assume that $$\qquad (x_j,p_j)\longrightarrow (x_0,p_0)\in X\times \mathbb R.$$
Choose holomorphic charts $\Psi_j:B_{2r}\to X$ with
$\Psi_j(0)=x_j$ and $(\Psi_j^*\omega)(0)=I$, whose coefficients
and inverses have uniform bounds. On $B_{rK_j}\subset \mathbb C^n$, put

\[
\begin{gathered}
\Phi_j(\zeta)=\Psi_j(\zeta/K_j),\qquad
v_j(z,p)=\psi_j (\Phi_j(z),p),\\
\omega_j=K_j^2\Phi_j^*\omega,\qquad \chi_j=\Phi_j^*\chi.
\end{gathered}
\]
Then $\omega_j\to I$ and $\chi_j\to0$ locally uniformly. Recall that we have uniform $C^0$ estimate for $\psi_j$ using Lemma \ref{exs:barriers}. This implies that $\psi_j$ has uniform $C^0$ bound on $$\mathcal O_1 \triangleq \{(x,p)\in X \times \mathbb R: \inf_{\epsilon< \epsilon_0} \min_{x\in X} \partial_t \varphi(x,0) -1 \le p \le \sup_{\epsilon< \epsilon_0} \max_{x\in X} \partial_t \varphi(x,1) +1\}$$
Define $D_r= B_r \times [\inf_{\epsilon< \epsilon_0} \min_{x\in X} \partial_t \varphi(x,0) -1, \sup_{\epsilon< \epsilon_0} \max_{x\in X} \partial_t \varphi(x,1) +1]$.
Combining this with Lemma~\ref{exg:contact-extension}, we get
\begin{equation}\label{e est vj}
||v_j||_{L^{\infty}(D_r)}+||\partial_X v_j||_{L^{\infty}(D_r)}+||\partial \bar \partial_X v_j||_{L^{\infty}(D_r)}\le C_1,
\qquad |\partial_X v_j|(0,p_j)=1.
\end{equation}
for a constant $C_1$ independent of $j,r$. In particular, $|v_j|_{C^{1,\alpha}(D_r)}(p)\le C_2$ for a uniform constant $C_2$ independent of $r$ and $j$. Since $|\partial_p v_j|\le 1$, we have that 
\begin{equation*}
    |v_j(\cdot, p)- v_j(\cdot, p')|\le |p-p'|.
\end{equation*}
Using interpolation, we have that:
\begin{equation}\label{e gradient 1}
    |\nabla_z v_j(z,p)- \nabla_z v_j(z,p')| \le C_3 |p-p'|^{\alpha/(1+\alpha)},
\end{equation}
for any $(z,p), (z,p')\in D_r$. Here $C_2$  is independent of $j$ and $r$. After diagonal extraction, we can find a function $v$ defined on $\mathbb C^n$ such that 
\begin{equation}\label{e con vj}
    v_j \rightarrow v, \nabla_z v_j \rightarrow \nabla_z v
\end{equation}
locally uniformly on $\mathbb C^n \times \mathbb R$. Denote 
\begin{equation*}
    \mathcal O_2 \triangleq \mathbb C^n \times [\inf_{\epsilon< \epsilon_0} \min_{x\in X} \partial_t \varphi(x,0) -1 , \sup_{\epsilon< \epsilon_0} \max_{x\in X} \partial_t \varphi(x,1) +1]
\end{equation*}

By (\ref{e est vj}), we have that: 
\begin{equation}\label{e bdd v}
    ||v||_{L^{\infty}(\mathcal O_2)}+ ||\partial_X v||_{L^{\infty}(\mathcal O_2)} + ||\partial \bar \partial_X v||_{L^{\infty}(\mathcal O_2)}\le C ,\,\,\, |\partial_X v|(0,p_0)=1
\end{equation}

\textbf{Step 2:} $v$ is $\Gamma_d$-solution. Fix a bounded interval $J\subset\mathbb R$ and a compact spatial
ball contained in the rescaled charts for all sufficiently large $j$.
Set
\[
F_j(z,p)=K_j^{-2}
F\!\left((\omega^{-1}\chi_{\psi_j})(\Phi_j(z),p)\right).
\]
All constants below are independent of $j$.

On the region $\{(x,p): \partial_t \varphi(x,0)< p < \partial_t \varphi(x,1)\}$, the root equation gives
\[
F_j=\frac{4\epsilon_j^{1/d}}{K_j^2}(-v_{j,pp}).
\]
On the region $\{(x,p): p < \partial_t \varphi(x,0) \text{ or } p>  \partial_t \varphi(x,1)\}$,
we have
\[
v_j(z,p)=\varphi_0(\Phi_j(z)), \text{ or } v_j(z,p)=\varphi_1(\Phi_j(z)) - p,
\qquad v_{j,pp}=0
\quad\text{almost everywhere}.
\]
Since the endpoints are fixed, smooth and admissible,
\[
0\le F_j
=K_j^{-2}F(\omega^{-1}\chi_{\varphi_i})(\Phi_j(z))
\le \frac{C}{K_j^2}
\]
there. Consequently, on the entire extended domain,
\begin{equation*}
0\le
F_j+\frac{4\epsilon_j^{1/d}}{K_j^2}v_{j,pp}
\le \frac{C}{K_j^2}
\qquad\text{almost everywhere}.
\end{equation*}

The extended transform is concave in $p$, with
$-1\le v_{j,p}\le0$. Therefore, for almost every fixed $z$,
\[
\int_J(-v_{j,pp})\,dp\le1.
\]
Integrating the preceding inequality yields
\[
\int_JF_j\,dp
\le
\frac{4\epsilon_j^{1/d}+C|J|}{K_j^2}.
\]
The quadratic spatial estimate
$|\chi_{\psi_j}|_\omega\le C(1+K_j^2)$
and the degree-$d$ polynomial growth of $g_2$ imply
$0\le F_j\le C$. Hence
\begin{equation}\label{e Fj zero}
\int_JF_j^d\,dp
\le C^{d-1}\int_JF_j\,dp
\le \frac{C_J}{K_j^2}
\longrightarrow0,
\end{equation}
where we used $0<\epsilon_j\le\epsilon_0$.
The estimate is essentially uniform in $z$ on the fixed spatial ball.

$v$ is a $\Gamma_d$ solution. 
Since $\psi_j$ is a $\epsilon$ geodesic, we can compute that:
\begin{equation}\label{e Fj expansion}
\begin{split}
   F_j^d & = \sum_{r=0}^d c_r^{(2)}K_j^{-2(d-r)}
       \sigma_r(\omega_j^{-1}(\chi_j + dd^c v_j)).
\end{split}
\end{equation}
Note that $\chi_j=O(K_j^{-2})$, $\omega_j \to \beta$.
Using (\ref{e est vj}), (\ref{e Fj zero}), (\ref{e Fj expansion}), we have that:
\begin{equation}\label{e converge vj}
   \lim_{j \rightarrow \infty} \int_J |\sigma_d(v_j)|dp =0. 
\end{equation}

Using (\ref{e est vj}) and Lemma \ref{bound_mat}, we have that $\sigma_k(dd^c v_j)\ge o(1)$ for $k \le d$. As a result, we have that $ v$ is d-subharmonic. 

For any smooth compactly supported function $\eta(z,p)$, we have that:
\begin{equation*}
    \int_{\mathbb C^n \times \mathbb R} \eta \big( (dd^c v_j)^d \wedge \beta^{n-d} - (dd^c v)^d \wedge \beta^{n-d}\big) dp= \int_{\mathbb C^n \times \mathbb R} (v_j- v) dd^c \eta \wedge \sum_{i=0}^{d-1} (dd^c v_j)^i \wedge (dd^c v)^{d-1-i}\wedge \beta^{n-d}.
\end{equation*}
Here $\beta$ is the K\"ahler form for the Euclidean metric on $\mathbb{C}^n$. The second integral converges uniformly to zero because of (\ref{e est vj}), (\ref{e con vj}) and (\ref{e bdd v}). Combining this with (\ref{e converge vj}), we have that:
\begin{equation*}
    \int_{\mathbb C^n \times \mathbb R} \eta (dd^c v)^d \wedge \beta^{n-d} dp=0
\end{equation*}
for any test function $\eta$. Using Fubini's theorem, we have that for almost every $p$:
$(dd^c v)^d \wedge \beta^{n-d} (\cdot, p)=0$.  Combining \cite[Theorem 3.6]{blocki2005weak} and \cite[Theorem 3.2]{dinew2017liouville}, we can get that $v$ 
 is a constant function on those slices. We can also
using the terminology in  \cite[Section 6]{szekelyhidi2018fully} to say that $v$ is a $\Gamma_d$-solution on those slices. Using \cite[Theorem 20]{szekelyhidi2018fully}, $v$ is a constant function on those slices. Since $v$ is continuous, we have that $v$ is a constant function on $\mathbb C^n \times \mathbb R$. This contradicts (\ref{e gradient 1}).
\end{proof}

\begin{prop}\label{lem space hess}
\begin{equation}
 |\chi_\psi|_\omega\leq C
 \quad \text{almost everywhere}.
\end{equation}
\end{prop}
\begin{proof}
    This Proposition follows from the combination of Lemma \ref{exg:contact-extension} and Lemma \ref{res:5-5}. 
\end{proof}

\subsection{Boundary and time Hessian estimate}

In this section, we work with $\varphi$ instead of $\psi$. 

We write $\widetilde{G}(H_{\varphi})= \log P_0(H_{\varphi})+ (d-1)\log (\frac{\varphi_{tt}}{4})$. 
Define its linearized operator as $\widetilde{L}= \widetilde{G}^{\alpha \bar\beta} \frac{1}{2}(\nabla_{\alpha} \nabla_{\bar \beta} + \nabla_{\bar \beta} \nabla_{\alpha}).$ For $j=0,1$, define $\zeta_j= \partial_X (\varphi- \varphi_j).$ 

\begin{lem}\label{lem lzeta}
  Let $\varphi$ be a $\epsilon$-geodesic.  For $j=0,1$, we have
    \begin{equation*}
        |\widetilde{L} \zeta_j|\le C (1+||\partial_X \varphi||_{\infty}) tr_{\Omega} D\widetilde{G}
    \end{equation*}
    for some constant $C$ independent of $\epsilon$.
\end{lem}
\begin{proof}
   \begin{equation}\label{e wide l}
       \begin{split}
          ( \widetilde{L} \zeta_j )_i  = \frac{1}{2}\widetilde{G}^{\alpha \bar \beta} \big( (\varphi- \varphi_j)_{i \alpha \bar \beta } + (\varphi- \varphi_j)_{i \bar \beta \alpha} \big)  = \widetilde{G}^{\alpha \bar \beta} (\varphi- \varphi_j)_{\alpha \bar \beta i} + \frac{1}{2} \widetilde{G}^{\alpha \bar \beta} R_{\alpha \bar \beta i}{}^k (\varphi- \varphi_j)_k
       \end{split}
   \end{equation} 
   Take the derivative of the $\epsilon$-geodesic $\widetilde{G}(H_{\varphi})= \log \epsilon$, we get:
   \begin{equation*}
       \widetilde{G}^{\alpha \bar \beta} (\chi_{\alpha \bar \beta i} + \varphi_{\alpha \bar \beta i})=0.
   \end{equation*}
   Plug this into (\ref{e wide l}). We can conclude the proof of this Lemma.
\end{proof}

\begin{prop}\label{prop varphitt}
There exists a constant $\epsilon_0>0$.     Let $\varphi$ be a $\epsilon$-geodesic with $\epsilon \in (0,\epsilon_0)$. We have that:
    \begin{equation}\label{exb:uniform-bound}
 \sup_{\partial(X\times[0,1])}|D\bar D \varphi|_\Omega\leq C,
 \qquad 0<\partial_t^2\varphi\leq C\quad\text{on }X\times[0,1].
\end{equation}
for a constant $C$ independent of $\epsilon$.
\end{prop}
\begin{proof}
    First, we derive boundary Hessian estimate. We can just estimate $|D\bar D \varphi|_\Omega$ on $X\times \{0\}$. The estimate for $|D\bar D \varphi|_\Omega$ on $X\times \{1\}$ is the same. By concavity of $\widetilde{G}$, according to Lemma \ref{ex:schur}, we have:
    \begin{equation}\label{e lvarphi}
        \widetilde{L}(\varphi- \underline{\varphi})\le \widetilde{G}(H_{\varphi})- \widetilde{G}(H_{\underline{\varphi}} -\epsilon_1 I)- \widetilde{G}^{\alpha \bar \beta} \epsilon_1 \delta_{\alpha \beta} \le \log (\epsilon) - \epsilon_1 tr_{\Omega} D\widetilde{G} + C \le  - \epsilon_1 tr_{\Omega} D\widetilde{G}+ C.
    \end{equation}
    Here $\epsilon_1>0$ is a small constant such that $G(H_{\underline{\varphi}} -\epsilon_1 I)>0$. For the last inequality of the above formula, we assume that $\epsilon \le e^{-C}$. By Lemma \ref{exs:barriers}, we have that $0\le  \varphi - \underline{\varphi}\le K_0 t(1-t).$
    Consider the barrier function:
 \[
 h_1= (1+||\partial_X \varphi||_{\infty})\big(A_0(\varphi- \underline{\varphi})+B_0 t \big)-\bigl((|\zeta_j|^2+\nu^2)^{1/2}-\nu\bigr).
\]
Here $A_0, B_0$ are big constants and $\nu$ is a small constant to be chosen later.  Using (\ref{e lvarphi}), Lemma \ref{lem lzeta}, and $\widetilde{L}(t)=0$, we can compute that:
\begin{equation*}
\begin{split}
\widetilde{L} h_1 &\le A_0(1+||\partial_X \varphi||_{\infty}) (-\epsilon_1 tr_{\Omega} D \widetilde{G})- \frac{ Re \langle \zeta_0, \widetilde{L} \zeta_0 \rangle}{(|\zeta_j|^2 + \nu^2)^{1/2}} \\
&\le -\epsilon_1  A_0(1+||\partial_X \varphi||_{\infty}) tr_{\Omega} D \widetilde{G} + C (1+ ||\partial_X \varphi||_{\infty}) tr_{\Omega} D\widetilde{G}.
\end{split}
\end{equation*}
Let $A_0$ be big such that $A_0 \epsilon_1 \ge C$. We can get that $\widetilde{L} h_1 \le 0$. Then, we can let $B_0$ be big such that at $t=1$, $h_1 \ge 0$. Note that $h_1$ vanishes on $t=0$. We can use the maximum principle on $h_1$ over $X\times [0,1]$ to get that $h_1 \ge 0$ on $X\times [0,1]$. Let $\nu$ in the definition of $h_1$ go to zero, and divide $h_1$ by $t$ and let $t$ go to zero. This implies that $\partial_t h_1 |_{t=0} \ge 0$. This implies that:
\begin{equation}\label{e normal tan}
    |\partial_X \partial_t \varphi|(\cdot,0)\le (1+ ||\partial_X \varphi||_{\infty})(A_0 K_0 + B_0).
\end{equation}
By Lemma \ref{res:5-5} and (\ref{e psi varphi spa}), we have that $||\partial_X \varphi||_{\infty}\le C$. Applying this to (\ref{e normal tan}), we can get:
\begin{equation}\label{e bdry mix hes}
    |\partial_X \partial_t \varphi|(\cdot,0) \le C.
\end{equation}

At an endpoint, we denote $b_i = \frac{\partial_X \partial_t \varphi_i}{2}$ which is bounded by (\ref{e bdry mix hes}). Denote $m= \min_{x\in X, i=0,1}g_2(\chi_{\varphi_i})$. Recall that the $\epsilon$-geodesic equation is:
\begin{equation*}
    (\frac{\partial_t^2 \varphi}{4})^{d-1} \big(\frac{\partial_t^2 \varphi}{4} g_2(\chi_{\varphi_i})- b_i^* D g_2(\chi_{\varphi_i}) b_i \big) = \epsilon.
\end{equation*}
This implies that:
\begin{equation*}
    \frac{\varphi_{tt}}{4}(x,i) \le \max\{\frac{2 b_i^* D g_2(\chi_{\varphi_i})b_i}{m}, (\frac{2\epsilon}{m}^{1/d})\},
\end{equation*}
for $i=0,1$. Take second order derivatives of the $\epsilon$-geodesic equation with respect to $t$. We can get:
\begin{equation*}
    \widetilde{L}(\partial^2_t \varphi) = - \widetilde{G}^{i \bar j, k \bar l} (H_{\varphi,t})_{i \bar j} (H_{\varphi,t})_{k \bar l} \ge 0.
\end{equation*}
Using the maximum principle for $\partial^2_t \varphi$, we can get that:
\begin{equation*}
    ||\partial^2_t \varphi||_{L^{\infty}(X\times [0,1])} \le ||\partial^2_t \varphi||_{L^{\infty}(\partial(X\times [0,1]))}\le C.
\end{equation*}
\end{proof}

\subsection{Complex Hessian estimate}

\begin{lem}\label{lem:dual-trace-estimate}
Assume
$-(n-1)\Lambda<\lambda_i\le\Lambda$ for every $i$, with $\Lambda>0$.
Then
\begin{equation}\label{eq:dual-trace-estimate}
 \begin{aligned}
 \mathcal L(\tr_\omega\chi_\psi)
 \ge{}&-\sum_kF^{i\bar j,r\bar s}
                  g_{i\bar j k}g_{r\bar s\bar k}\\
      &-Cn\Lambda\,\mathcal F,
 \end{aligned}
\end{equation}
where $C$ depends only on $n$ and the curvature bound of $\omega$.
\end{lem}

\begin{proof}
The reference metric is parallel. Covariantly differentiating
(\ref{equ psi}) first in $k$ and then in $\bar k$
therefore gives
\[
 \begin{aligned}
 0={}&\sum_iF_i\nabla_{\bar k}\nabla_k g_{i\bar i}
   +F^{i\bar j,r\bar s}g_{i\bar j\bar k}g_{r\bar s k}
   +4\epsilon^{1/d}\psi_{ppk\bar k}.
 \end{aligned}
\]
By symmetry of $D^2F$, this is
\begin{equation}\label{eq:dual-trace-differentiated}
 \sum_iF_i\nabla_{\bar k}\nabla_k g_{i\bar i}
 +4\epsilon^{1/d}\psi_{ppk\bar k}
 =-F^{i\bar j,r\bar s}g_{i\bar j k}g_{r\bar s\bar k}.
\end{equation}

The commutator acts on both indices of $\chi_\psi$ and gives
\begin{equation}\label{e commute chi psi}
 \begin{aligned}
 ([\nabla_i,\nabla_{\bar k}]\chi_\psi)_{k\bar i}
 &=-\sum_\ell R_{i\bar k k}{}^\ell g_{\ell\bar i}
   +\sum_\ell R_{i\bar k\ell\bar i}g_{k\bar\ell}\\
 &=(\lambda_k-\lambda_i)R_{i\bar i k\bar k}.
 \end{aligned}
\end{equation}

Since $\chi$ and $\omega$ are independent of $p$, $(\tr_\omega\chi_\psi)_{pp}=\sum_k\psi_{ppk\bar k}.$
Summing (\ref{eq:dual-trace-differentiated}) and using
(\ref{e commute chi psi}) proves the exact identity
\begin{equation}\label{eq:dual-trace-exact}
 \begin{aligned}
 \mathcal L(\tr_\omega\chi_\psi)
 ={}-\sum_kF^{i\bar j,r\bar s}g_{i\bar j k}g_{r\bar s\bar k}+\sum_{i,k}F_iR_{i\bar i k\bar k}(\lambda_k-\lambda_i).
 \end{aligned}
\end{equation}
The assumed eigenvalue interval gives $|\lambda_k-\lambda_i|\le n\Lambda$.
Consequently,
\[
 \left|\sum_{i,k}F_iR_{i\bar i k\bar k}(\lambda_k-\lambda_i)\right|
 \le n\Lambda\sum_iF_i\sum_k|R_{i\bar i k\bar k}|
 \le Cn\Lambda\,\mathcal F.
\]
Together with (\ref{eq:dual-trace-exact}), this proves
(\ref{eq:dual-trace-estimate}).
\end{proof}

\begin{lem}[The second derivative in $p$]\label{lem:dual-pp-identity}
The dual solution satisfies
\begin{equation}\label{eq:dual-pp-identity}
 \mathcal L\psi_{pp}
 =-F^{i\bar j,r\bar s}g_{i\bar j p}g_{r\bar s p}\ge0.
\end{equation}
\end{lem}

\begin{proof}
Twice differentiating (\ref{equ psi}) in $p$ gives
\[
 F^{i\bar j}(g_{i\bar j})_{pp}
 +F^{i\bar j,r\bar s}g_{i\bar j p}g_{r\bar s p}
 +4\epsilon^{1/d}\psi_{pppp}=0.
\]
The fixed background is independent of $p$, so
$(g_{i\bar j})_{pp}=(\psi_{pp})_{i\bar j}$.
This proves the equality in (\ref{eq:dual-pp-identity}), and its
nonnegativity follows from concavity of $F$ by Lemma \ref{ex:ideal}.
\end{proof}

\begin{lem}[Gradient of the first derivative in $p$]
\label{lem:dual-p-gradient}
There is a constant $C$, depending only on $n$ and the curvature bound
of $\omega$, such that
\begin{equation}\label{eq:dual-p-gradient-estimate}
 \begin{aligned}
 \mathcal L|\partial\psi_p|_\omega^2
 \ge{}&\sum_{i,k}F_i
       \bigl(|\psi_{pki}|^2+|\psi_{pk\bar i}|^2\bigr) +8\epsilon^{1/d}\sum_k|\psi_{ppk}|^2 -2\operatorname{Re}\sum_k\overline{\psi_{pk}}\,
       F^{i\bar j,r\bar s}g_{i\bar j k}g_{r\bar s p}-C|\partial\psi_p|_\omega^2\mathcal F.
 \end{aligned}
\end{equation}
\end{lem}

\begin{proof}
Differentiating (\ref{equ psi}) once in $p$ gives
\[
 F^{i\bar j}(\psi_p)_{i\bar j}
       +4\epsilon^{1/d}\psi_{ppp}=0.
\]
Take its covariant derivative in $k$. We get:
\begin{equation}\label{eq:dual-mixed-differentiated}
 \sum_iF_i\nabla_k(\psi_p)_{i\bar i}
       +4\epsilon^{1/d}\psi_{pppk}
 =-F^{i\bar j,r\bar s}g_{i\bar j k}g_{r\bar s p}.
\end{equation}

To expand the squared norm, write
$|\partial\psi_p|_\omega^2
=\omega^{k\bar\ell}\psi_{pk}\overline{\psi_{p\ell}}$.
At the normal coordinate center, the spatial product rule gives
\begin{equation}\label{eq:dual-spatial-bochner}
 \begin{aligned}
 \nabla_i\nabla_{\bar i}|\partial\psi_p|_\omega^2
 ={}&\sum_k\bigl(|\psi_{pki}|^2+|\psi_{pk\bar i}|^2\bigr)+2\operatorname{Re}\sum_k
       \overline{\psi_{pk}}\,\nabla_k(\psi_p)_{i\bar i}+\operatorname{Re}\sum_{k,\ell}
       R_{i\bar i k}{}^\ell\psi_{p\ell}\overline{\psi_{pk}}.
 \end{aligned}
\end{equation}

The metric is independent of the real variable $p$, so
\begin{equation}\label{eq:dual-real-p-product}
 \partial_{pp}|\partial\psi_p|_\omega^2
 =2\sum_k|\psi_{ppk}|^2
  +2\operatorname{Re}\sum_k\overline{\psi_{pk}}\psi_{pppk}.
\end{equation}
Multiply (\ref{eq:dual-spatial-bochner}) by $F_i$, sum in $i$, add
$4\epsilon^{1/d}$ times (\ref{eq:dual-real-p-product}), and use
(\ref{eq:dual-mixed-differentiated}). This gives the exact formula
\begin{equation}\label{eq:dual-p-gradient-exact}
 \begin{aligned}
 \mathcal L|\partial\psi_p|_\omega^2
 ={}&\sum_{i,k}F_i
       \bigl(|\psi_{pki}|^2+|\psi_{pk\bar i}|^2\bigr)+8\epsilon^{1/d}\sum_k|\psi_{ppk}|^2-2\operatorname{Re}\sum_k\overline{\psi_{pk}}\,
       F^{i\bar j,r\bar s}g_{i\bar j k}g_{r\bar s p}\\
       &+\operatorname{Re}\sum_{i,k,\ell}
       F_iR_{i\bar i k}{}^\ell\psi_{p\ell}\overline{\psi_{pk}}.
 \end{aligned}
\end{equation}
For each $i$, the curvature endomorphism has operator norm bounded by
a fixed multiple of $|\operatorname{Rm}(\omega)|$. As $F_i>0$,
\[
 \left|\operatorname{Re}\sum_{i,k,\ell}
 F_iR_{i\bar i k}{}^\ell\psi_{p\ell}\overline{\psi_{pk}}\right|
 \le C|\partial\psi_p|_\omega^2\sum_iF_i.
\]
Substituting this bound into (\ref{eq:dual-p-gradient-exact}) proves
(\ref{eq:dual-p-gradient-estimate}).
\end{proof}

\begin{lem}\label{lem boundary mix hess}
    There exists a constant $C$ depending on $\varphi_j$, $j=0,1$, $\omega$, $\chi$ such that
    \begin{equation*}
        \frac{|\partial \psi_p|_{\omega}^2}{-\psi_{pp}}(x, \varphi_t(x,j))\le C
    \end{equation*}
    for any $x\in X$ and $j=0,1$.
\end{lem}
\begin{proof}
By partial Legendre transform, we have:
\begin{equation*}
    \psi_{pi}=  \frac{\varphi_{ti}}{\varphi_{tt}}, -\psi_{pp}= \frac{1}{\varphi_{tt}}, \psi_{i \bar k}= \varphi_{i \bar k}- \frac{\varphi_{ti}\varphi_{t \bar k}}{\varphi_{tt}}.
\end{equation*}
    Using Lemma \ref{lem psi hess}, we have:
    \begin{equation*}
        tr_{\omega} (\chi_{\varphi_j}- \chi_{\psi})= \frac{\varphi_{ti}\varphi_{t \bar k}}{\varphi_{tt}} = \frac{|\partial \psi_p|^2_{\omega}}{-\psi_{pp}}.
    \end{equation*}
    Since $tr_{\omega} \chi_{\varphi_j}$ is bounded and $tr_{\omega} \chi_{\psi}$ is nonnegative, we conclude the proof of this Lemma.
\end{proof}

\begin{lem}\label{in:log-ratio}
Suppose $\chi_\psi\leq\Lambda\omega$. At every point where
$|\partial\psi_p|_\omega^2>0$, one has
\begin{equation}\label{in:log-ratio-inequality}
 \begin{aligned}
 \mathcal L\log\frac{|\partial\psi_p|^2}{-\psi_{pp}}\geq{}&
 \frac{\sum_{i,j}F_i|\psi_{pj\bar i}|^2}{|\partial\psi_p|^2} -2\Re\sum_iF_i\left(\log\frac{|\partial\psi_p|^2}{-\psi_{pp}}\right)_{i}
                  \frac{\sum_j\psi_{pj\bar i}\overline{\psi_{pj}}}{|\partial\psi_p|^2}\\
 &+\frac{F^{r\bar s,\ell\bar m}g_{r\bar s p}g_{\ell\bar m p}}{\psi_{pp}}-\frac2{|\partial\psi_p|^2}\Re\sum_j\overline{\psi_{pj}}\,F^{r\bar s,\ell\bar m}g_{r\bar s j}g_{\ell\bar m p}\\
 &-C\mathcal F-4\eps^{1/d}\left(\log\frac{|\partial\psi_p|^2}{-\psi_{pp}}\right)_{p}^2.
 \end{aligned}
\end{equation}
\end{lem}
\begin{proof}
\eqref{lem:dual-pp-identity} and \eqref{lem:dual-p-gradient} give
\begin{equation}\label{in:raw-log-ratio}
 \begin{aligned}
 \mathcal L\log\frac{|\partial\psi_p|^2}{-\psi_{pp}}\geq{}&
 \sum_iF_i\left\{
 \frac{\sum_j(|\psi_{pji}|^2+|\psi_{pj\bar i}|^2)}{|\partial\psi_p|^2}
 -\frac{\left|\left(|\partial\psi_p|^2\right)_i\right|^2}{|\partial\psi_p|^4}
 +\frac{|\psi_{ppi}|^2}{\psi_{pp}^2}
 \right\}\\
 &+4\eps^{1/d}\left\{
 \frac{2\sum_j|\psi_{ppj}|^2}{|\partial\psi_p|^2}
 -\frac{\left(\partial_p|\partial\psi_p|^2\right)^2}{|\partial\psi_p|^4}
 +\frac{\psi_{ppp}^2}{\psi_{pp}^2}
 \right\}\\
 &+\frac{F^{r\bar s,\ell\bar m}g_{r\bar s p}g_{\ell\bar m p}}{\psi_{pp}}-\frac2{|\partial\psi_p|^2}\Re\sum_j\overline{\psi_{pj}}\,F^{r\bar s,\ell\bar m}g_{r\bar s j}g_{\ell\bar m p}
 -C\mathcal F.
 \end{aligned}
\end{equation}
For the spatial terms, differentiate the norm and use Hermitian symmetry:
\begin{equation}\label{in:spatial-log-identity}
 \begin{aligned}
 \frac{\left(|\partial\psi_p|^2\right)_i}{|\partial\psi_p|^2}
 &=\frac{\sum_j\psi_{pji}\overline{\psi_{pj}}}{|\partial\psi_p|^2}
   +\frac{\sum_j\psi_{pi\bar j}\psi_{pj}}{|\partial\psi_p|^2}=\frac{\psi_{ppi}}{\psi_{pp}}+\left(\log\frac{|\partial\psi_p|^2}{-\psi_{pp}}\right)_{i}.
 \end{aligned}
\end{equation}
Substitute \eqref{in:spatial-log-identity} into this bound and complete
the square. We get:
\begin{equation}\label{in:spatial-completed-square}
 \begin{aligned}
  &\frac{\sum_j(|\psi_{pji}|^2+|\psi_{pj\bar i}|^2)}{|\partial\psi_p|^2}
 -\frac{\left|\left(|\partial\psi_p|^2\right)_i\right|^2}{|\partial\psi_p|^4}
 +\frac{|\psi_{ppi}|^2}{\psi_{pp}^2} \\
 & \ge  \left|\frac{\sum_j\psi_{pji}\overline{\psi_{pj}}}{|\partial\psi_p|^2}\right|^2 +\frac{\sum_j |\psi_{pj\bar i}|^2}{|\partial\psi_p|^2}  -\frac{\left|\left(|\partial\psi_p|^2\right)_i\right|^2}{|\partial\psi_p|^4}
 +\frac{|\psi_{ppi}|^2}{\psi_{pp}^2} \\
 & = |\frac{\psi_{ppi}}{\psi_{pp}}+\left(\log\frac{|\partial\psi_p|^2}{-\psi_{pp}}\right)_{i}- \frac{\sum_j\psi_{pi\bar j}\psi_{pj}}{|\partial\psi_p|^2}|^2 +\frac{\sum_j |\psi_{pj\bar i}|^2}{|\partial\psi_p|^2} - | \frac{\psi_{ppi}}{\psi_{pp}}+\left(\log\frac{|\partial\psi_p|^2}{-\psi_{pp}}\right)_{i} |^2 +\frac{|\psi_{ppi}|^2}{\psi_{pp}^2} \\
 &  = -2 Re \big( \frac{\psi_{ppi}}{\psi_{pp}} +(\log \frac{|\partial \psi_p|}{-\psi_{pp}})_i \big)(\sum_j \frac{\psi_{p\bar i j}\psi_{p \bar j}}{|\partial \psi_p|^2}) + |\sum_j \frac{\psi_{p i \bar{j}}\psi_{pj}}{|\partial \psi_p|^2}|^2 + \sum_j \frac{|\psi_{p j \bar i}|^2}{|\partial \psi_p|^2} + |\frac{\psi_{ppi}}{\psi_{pp}}|^2 \\
 & = -2 Re (\log \frac{|\partial \psi_p|^2}{-\psi_{pp}})_i \sum_j \frac{\psi_{p \bar i j} \psi_{p \bar j}}{|\partial \psi_p|^2} + |\sum_j \frac{\psi_{pi \bar j}\psi_{pj}}{|\partial \psi_p|^2}- \frac{\psi_{ppi}}{\psi_{pp}}|^2 + \sum_j \frac{|\psi_{pj \bar i}|^2}{|\partial \psi_p|^2}.
 \end{aligned}
\end{equation}
For the temporal bracket, differentiation and Cauchy--Schwarz give
\[
 \begin{gathered}
 \partial_p|\partial\psi_p|^2
       =2\Re\sum_j\overline{\psi_{pj}}\psi_{ppj},\\
 \frac{2\sum_j|\psi_{ppj}|^2}{|\partial\psi_p|^2}
       \geq\frac12
          \left(\frac{\partial_p|\partial\psi_p|^2}{|\partial\psi_p|^2}\right)^2.
 \end{gathered}
\]
Since
\[
 \frac{\partial_p|\partial\psi_p|^2}{|\partial\psi_p|^2}
 =\frac{\psi_{ppp}}{\psi_{pp}}+\left(\log\frac{|\partial\psi_p|^2}{-\psi_{pp}}\right)_{p},
\]
we compute that:
\begin{equation}\label{e sec bra}
 \begin{aligned}
 & \frac{2\sum_j|\psi_{ppj}|^2}{|\partial\psi_p|^2}
 -\frac{\left(\partial_p|\partial\psi_p|^2\right)^2}{|\partial\psi_p|^4}
 +\frac{\psi_{ppp}^2}{\psi_{pp}^2} \\
 & \ge \frac{1}{2}\frac{\left(\partial_p|\partial\psi_p|^2\right)^2}{|\partial\psi_p|^4} -\frac{\left(\partial_p|\partial\psi_p|^2\right)^2}{|\partial\psi_p|^4}
 +\frac{\psi_{ppp}^2}{\psi_{pp}^2} \\
 & = - \frac{1}{2}\frac{\left(\partial_p|\partial\psi_p|^2\right)^2}{|\partial\psi_p|^4}
 +\frac{\psi_{ppp}^2}{\psi_{pp}^2} \\
 & \ge - \frac{\psi_{ppp}^2}{\psi_{pp}^2} - |(\log \frac{|\partial \psi_p|^2}{-\psi_{pp}})_p|^2 + \frac{\psi_{ppp}^2}{\psi_{pp}^2} \\
 & =  - |(\log \frac{|\partial \psi_p|^2}{-\psi_{pp}})_p|^2.
 \end{aligned}
\end{equation}
Plug (\ref{in:spatial-log-identity}) and (\ref{e sec bra}) into the first and second brackets of (\ref{in:raw-log-ratio}). We can conclude the proof of this Lemma.
\end{proof}

\begin{prop}\label{in:mixed-bound}
Suppose $\chi_\psi\leq\Lambda\omega$. There are fixed $Q$ and $\eps_2>0$ such that\\
$|\partial\psi_p|_\omega^2/(-\psi_{pp})\leq Q$ for
$0<\eps\leq\eps_2$.
\end{prop}
\begin{proof}
Set
\[
 \begin{gathered}
 H(r)=-\tfrac1{64}\log(n\Lambda+1-r),\qquad
 \gamma(r)=-\tfrac1{64}\log(2-r),\\
 \begin{aligned}
 \Psi={}&\log\frac{|\partial\psi_p|^2}{-\psi_{pp}}+H(\tr_\omega\chi_\psi)-A(\psi-\varphi_0)+\gamma(-\psi_p).
 \end{aligned}
 \end{gathered}
\]
Both $H'$ and $\gamma'$ have fixed positive lower bounds and fixed
upper bounds, and $H''=64(H')^2$, $\gamma''=64(\gamma')^2$.
We compute that 

\begin{equation}\label{e lhlgamma}
\begin{aligned}
 \mathcal L\bigl[H(\operatorname{tr}_\omega\chi_\psi)\bigr]
={}&H'\mathcal L(\operatorname{tr}_\omega\chi_\psi) +
H''\sum_i F_i
\bigl|(\operatorname{tr}_\omega\chi_\psi)_i\bigr|^2
 +4\epsilon^{1/d}H''
\bigl((\operatorname{tr}_\omega\chi_\psi)_p\bigr)^2\\
\mathcal L[\gamma(-\psi_p)]
={}&-\gamma'\mathcal L\psi_p+
\gamma''\sum_i F_i|\psi_{pi}|^2 +4\epsilon^{1/d}\gamma''\psi_{pp}^2=\gamma''\sum_i F_i|\psi_{pi}|^2 +4\epsilon^{1/d}\gamma''\psi_{pp}^2.
\end{aligned}
\end{equation}
In the second line of the above formula, we use $\mathcal L \psi_p=0$ from  Lemma \ref{lem lpsi phi0}
Combining the above formula, Lemma \ref{lem lpsi phi0} and Lemma \ref{in:log-ratio} and Proposition \ref{lem space hess}, we get:
\begin{equation}\label{e lpsi}
    \begin{split}
        \mathcal L(\Psi) \geq{}&
 \frac{\sum_{i,j}F_i|\psi_{pj\bar i}|^2}{|\partial\psi_p|^2} -2\Re\sum_iF_i\left(\log\frac{|\partial\psi_p|^2}{-\psi_{pp}}\right)_{i}
                  \frac{\sum_j\psi_{pj\bar i}\overline{\psi_{pj}}}{|\partial\psi_p|^2}\\
 &+\frac{F^{r\bar s,\ell\bar m}g_{r\bar s p}g_{\ell\bar m p}}{\psi_{pp}}-\frac2{|\partial\psi_p|^2}\Re\sum_j\overline{\psi_{pj}}\,F^{r\bar s,\ell\bar m}g_{r\bar s j}g_{\ell\bar m p}\\
 &-C\mathcal F-4\eps^{1/d}\left(\log\frac{|\partial\psi_p|^2}{-\psi_{pp}}\right)_{p}^2 \\
 & + H' ( -\sum_kF^{i\bar j,r\bar s}
                  g_{i\bar j k}g_{r\bar s\bar k}-Cn\Lambda\,\mathcal F ) +
H''\sum_i F_i
\bigl|(\operatorname{tr}_\omega\chi_\psi)_i\bigr|^2
 +4\epsilon^{1/d}H''
\bigl((\operatorname{tr}_\omega\chi_\psi)_p\bigr)^2\\
&+ \gamma''\sum_i F_i|\psi_{pi}|^2 +4\epsilon^{1/d}\gamma''\psi_{pp}^2 \\
& + A(\eta \mathcal F +m_0).
    \end{split}
\end{equation}

At an interior maximum of $\Psi$ with $|\partial\psi_p|^2>0$,
\begin{equation}\label{in:ratio-stationarity}
 \begin{aligned}
 \left(\log\frac{|\partial\psi_p|^2}{-\psi_{pp}}\right)_{i}
 &=-H'(\tr_\omega\chi_\psi)_i
                 +A(\psi-\varphi_0)_i+\gamma'\psi_{pi},\\
 \left(\log\frac{|\partial\psi_p|^2}{-\psi_{pp}}\right)_{p}
 &=-H'(\tr_\omega\chi_\psi)_p
                 +A\psi_p+\gamma'\psi_{pp}.
 \end{aligned}
\end{equation}

Using the first equation of (\ref{in:ratio-stationarity}), we compute that:
\begin{equation}\label{e re f log}
\begin{aligned}
&-2\operatorname{Re}\sum_i F_i
  \left(\log\frac{|\partial\psi_p|^2}{-\psi_{pp}}\right)_i
  \frac{\sum_j\psi_{pj\bar i}\overline{\psi_{pj}}}
       {|\partial\psi_p|^2}
\\
&={}-2\operatorname{Re}\sum_i F_i
  \Bigl[
    -H'(\operatorname{tr}_{\omega}\chi_{\psi})_i
    +A(\psi-\varphi_0)_i
    +\gamma'\psi_{pi}
  \Bigr] \times
  \frac{\sum_j\psi_{pj\bar i}\overline{\psi_{pj}}}
       {|\partial\psi_p|^2} \\
&\ge 
  -\sum_i F_i
  \frac{\sum_j|\psi_{pj\bar i}|^2}
       {8|\partial\psi_p|^2}
  -8(H')^2\sum_i F_i
    |(\operatorname{tr}_{\omega}\chi_{\psi})_i|^2 -
  \sum_i F_i
  \frac{\sum_j|\psi_{pj\bar i}|^2}
       {8|\partial\psi_p|^2}
  -8(\gamma')^2\sum_i F_i|\psi_{pi}|^2
\\
&\quad-
  2\operatorname{Re}\sum_i F_i A(\psi-\varphi_0)_i
  \frac{\sum_j\psi_{pj\bar i}\overline{\psi_{pj}}}
       {|\partial\psi_p|^2} \\
& = -\sum_i F_i
  \frac{\sum_j|\psi_{pj\bar i}|^2}
       {4|\partial\psi_p|^2} - \frac{H''}{8}\sum_i F_i
    |(\operatorname{tr}_{\omega}\chi_{\psi})_i|^2  - \frac{\gamma''}{8} \sum_i F_i|\psi_{pi}|^2 \\
    & \quad-
  2\operatorname{Re}\sum_i F_i A(\psi-\varphi_0)_i
  \frac{\sum_j\psi_{pj\bar i}\overline{\psi_{pj}}}
       {|\partial\psi_p|^2} .
\end{aligned}
\end{equation}

Apply Cauchy--Schwarz and Young's inequality:
\begin{equation}\label{in:source-allocation}
 \begin{aligned}
 &\frac2{|\partial\psi_p|^2}
   \left|\sum_j\overline{\psi_{pj}}\,F^{r\bar s,\ell\bar m}g_{r\bar s j}g_{\ell\bar m p}\right| \leq\frac{F^{r\bar s,\ell\bar m}g_{r\bar s p}g_{\ell\bar m p}}{4\psi_{pp}}
       -\frac{4(-\psi_{pp})}{|\partial\psi_p|^2}\,\sum_hF^{r\bar s,\ell\bar m}g_{r\bar s h}g_{\ell\bar m\bar h}.
 \end{aligned}
\end{equation}
Both terms on the right are nonnegative by concavity and $\psi_{pp}<0$.

Use the second equation of (\ref{in:ratio-stationarity}), we get
\begin{equation}\label{e log partial psi}
 \begin{aligned}
 4\eps^{1/d}\left(\log\frac{|\partial\psi_p|^2}{-\psi_{pp}}\right)_{p}^2
 \leq{}&12\eps^{1/d}(H')^2(\tr_\omega\chi_\psi)_p^2 +12\eps^{1/d}A^2+12\eps^{1/d}(\gamma')^2\psi_{pp}^2,
 \end{aligned}
\end{equation}

Plug (\ref{e re f log}), (\ref{in:source-allocation}) and (\ref{e log partial psi}) into (\ref{e lpsi}) and use $H''= 64(H')^2$, $\gamma''= 64(\gamma')^2$. We can get:
\begin{equation}\label{in:ratio-energy}
 \begin{aligned}
 0\geq\mathcal L\Psi\geq{}&
 \frac{3\sum_{i,j}F_i|\psi_{pj\bar i}|^2}{4|\partial\psi_p|^2}
 +\frac{3F^{r\bar s,\ell\bar m}g_{r\bar s p}g_{\ell\bar m p}}{4\psi_{pp}}\\
 &-\left(H'-\frac{4(-\psi_{pp})}{|\partial\psi_p|^2}\right)
                        \sum_hF^{r\bar s,\ell\bar m}g_{r\bar s h}g_{\ell\bar m\bar h}\\
 &+56(\gamma')^2\sum_iF_i|\psi_{pi}|^2\\
 &+\{A\eta-C_1-C_2n\Lambda H'\}\mathcal F
       +Am_0-12\eps^{1/d}A^2-\mathcal R,
 \end{aligned}
\end{equation}
where the only remaining coupling is
\[
 \mathcal R=\frac{2A}{|\partial\psi_p|^2}\Re
        \sum_{i,j}F_i(\psi-\varphi_0)_i
                     \psi_{pj\bar i}\overline{\psi_{pj}}.
\]
Split these ordered pairs into
$I=\{(i,j):F_i\leq2F_j\}$ and $J=\{(i,j):F_i>2F_j\}$.
On $I$, Young's inequality gives
\begin{equation}\label{in:comparable-pairs}
 \begin{aligned}
 |\frac{2A}{|\partial\psi_p|^2}\Re
        \sum_{(i,j) \in I}F_i(\psi-\varphi_0)_i
                     \psi_{pj\bar i}\overline{\psi_{pj}}|
 &\leq\frac{\sum_{i,j}F_i|\psi_{pj\bar i}|^2}{4|\partial\psi_p|^2}
   +\frac{4A^2}{|\partial\psi_p|^2}
             \sum_I F_i|(\psi-\varphi_0)_i|^2|\psi_{pj}|^2\\
 &\leq\frac{\sum_{i,j}F_i|\psi_{pj\bar i}|^2}{4|\partial\psi_p|^2}
   +\frac{8A^2K}{|\partial\psi_p|^2}\sum_jF_j|\psi_{pj}|^2.
 \end{aligned}
\end{equation}
For $(i,j)\in J$, one has $\lambda_i<\lambda_j$, and the two-sided
spatial bound gives
\[
 \frac{F_i-F_j}{\lambda_j-\lambda_i}\geq\frac{F_i}{2n\Lambda}.
\]
There are no repeated-eigenvalue pairs in $J$. Since
$(\chi_\psi)_p=\ddc_X\psi_p$, \cite[(67)]{szekelyhidi2018fully} gives
\[
 -F^{r\bar s,\ell\bar m}g_{r\bar s p}g_{\ell\bar m p}
 \geq - \sum_{i \neq j} \frac{f_i - f_j}{\lambda_i - \lambda_j} |\psi_{p j \bar i}|^2 \ge - \sum_{(i,j)\in J} \frac{f_i - f_j}{\lambda_i - \lambda_j} |\psi_{p j \bar i}|^2 \ge \frac{1}{2n\Lambda} \sum_J F_i |\psi_{pj \bar {i}}|^2.
\]
This implies that:
\begin{equation}\label{in:separated-pairs}
 \begin{aligned}
 |\frac{2A}{|\partial\psi_p|^2}\Re
        \sum_{(i,j) \in J}F_i(\psi-\varphi_0)_i
                     \psi_{pj\bar i}\overline{\psi_{pj}}|
 &\leq \frac{\sum_J F_i|\psi_{pj\bar i}|^2}
                 {8n\Lambda(-\psi_{pp})}
 +\frac{8A^2n\Lambda(-\psi_{pp})}{|\partial\psi_p|^4}
       \sum_JF_i|(\psi-\varphi_0)_i|^2|\psi_{pj}|^2\\
 &\leq \frac{F^{r\bar s,\ell\bar m}g_{r\bar s p}g_{\ell\bar m p}}
                 {4\psi_{pp}}
       +\frac{8A^2n\Lambda K(-\psi_{pp})}{|\partial\psi_p|^2}\mathcal F.
 \end{aligned}
\end{equation}
Here
\[
 \sum_J F_i|(\psi-\varphi_0)_i|^2|\psi_{pj}|^2
       \leq K|\partial\psi_p|^2\mathcal F.
\]
Thus \eqref{in:comparable-pairs} spends the mixed Hessian energy,
whereas \eqref{in:separated-pairs} spends the nonnegative spectral
contraction divided by $-\psi_{pp}$.

Choose $A\eta\geq4(C_1+C_2n\Lambda\max H'+1)$. Put
\[
 Q_*=\max\left\{1,\frac8{\min H'},
       \frac{32An\Lambda K}{\eta},
       \frac{8A^2K}{c(\min\gamma')^2}\right\}.
\]
If
\[
 \frac{|\partial\psi_p|^2}{-\psi_{pp}}\geq Q_*,
\]
then
\[
 H'-\frac{4(-\psi_{pp})}{|\partial\psi_p|^2}\geq H'/2.
\]
The last term of \eqref{in:separated-pairs} costs at most
$A\eta\mathcal F/4$. Since $-\psi_{pp}\geq c$, at such a point
$|\partial\psi_p|^2\geq cQ_*$. Therefore the last term of
\eqref{in:comparable-pairs} costs at most
$(\min\gamma')^2\sum_iF_i|\psi_{pi}|^2$.
After these allocations, all remaining concavity and Hessian-square
terms in \eqref{in:ratio-energy} are nonnegative. In particular,
\[
 0\geq\mathcal L\Psi
 \geq\tfrac{A\eta}{2}\mathcal F+Am_0-12\eps^{1/d}A^2>0
\]
provided $0<\eps_2\leq\eps_1$ is small enough that
$4\eps_2^{1/d}\leq m_0/(6A)$.
Hence the quotient at an interior maximum is less than $Q_*$.
If $\partial\psi_p\equiv0$, the conclusion is immediate. Otherwise
\[
 \log\frac{|\partial\psi_p|^2}{-\psi_{pp}}\longrightarrow-\infty
 \quad\text{as }|\partial\psi_p|^2\longrightarrow0,
\]
so the maximum test is valid on $|\partial\psi_p|^2>0$.
Boundary values are controlled by Lemma \ref{lem boundary mix hess}.
The fixed oscillation of
$H(\tr_\omega\chi_\psi)-A(\psi-\varphi_0)+\gamma(-\psi_p)$
converts the maximum estimate into a uniform global bound for
$|\partial\psi_p|_\omega^2/(-\psi_{pp})$.
\end{proof}

\begin{prop}\label{in:full-hessian}
There exist $\eps_0>0$ and $C$, depending only on the fixed geometric,
polynomial and endpoint data, such that every smooth admissible
solution of the weighted family with $0<\eps\leq\eps_0$ satisfies
\[
 \|D\bar D\varphi\|_{L^\infty}\leq C
\]
in the fixed nondegenerate product metric.
\end{prop}
\begin{proof}
By partial Legendre transform, we have:
\begin{equation*}
    \psi_{pi}= \frac{\varphi_{ti}}{\varphi_{tt}}, -\psi_{pp}= \frac{1}{\varphi_{tt}}.
\end{equation*}
This implies that
\begin{equation*}
    \frac{\varphi_{ti}\varphi_{t \bar{j}}}{\varphi_{tt}}= \frac{\psi_{pi}\psi_{p\bar j}}{-\psi_{pp}}.
\end{equation*}
By Proposition \ref{in:mixed-bound}, $\frac{\psi_{pi}\psi_{p\bar j}}{-\psi_{pp}}$ is uniformly bounded, so is $ \frac{\varphi_{ti}\varphi_{t \bar{j}}}{\varphi_{tt}}$. By Proposition \ref{prop varphitt}, \\ $||\varphi_{tt}||_{L^{\infty}(X\times [0,1])}$ is uniformly bounded. As a result, $||\varphi_{ti}||_{L^{\infty}(X\times [0,1])}$ is uniformly bounded. Combining Proposition \ref{in:mixed-bound}, Lemma \ref{lem space hess} and Lemma \ref{lem psi hess}, we can get that $||\partial \bar \partial_X \varphi||_{L^{\infty}(X\times [0,1])}$ is bounded.
\end{proof}

\subsection{Real Hessian Estimate}

In this subsection, we give the uniform real Hessian estimate of $\varphi_\epsilon$.

\begin{prop}\label{in:real-hessian}
Let $\varphi_0,\varphi_1\in\mathcal H$. There are constants
$\epsilon_*>0$ and $C$, depending only on the fixed geometric,
polynomial and endpoint data, such that every smooth admissible
solution $\varphi_\epsilon$ of \eqref{ex:family}, $0<\epsilon\le\epsilon_*$,
satisfies
\begin{equation}\label{in:real-hessian-bound}
 \sup_{X\times[0,1]}|\nabla^2_{X,t}\varphi_\epsilon|\le C.
\end{equation}
Here $\nabla^2_{X,t}$ denotes the real Hessian for the product of
the Riemannian metric of $\omega$ and $dt^2$.
\end{prop}

We write $\varphi=\varphi_\epsilon$. Since
\[\nabla_X^2\varphi(x,t)=\nabla_X^2\psi(x,p)
       +\frac{d_X\psi_p\otimes d_X\psi_p}{-\psi_{pp}},\] 
it is sufficient to prove $\nabla^2_X \psi$ is uniformly bounded, according Proposition \ref{prop varphitt}, Proposition \ref{in:mixed-bound}. Set 
\[w=\psi-\varphi_0-C',\qquad M=\omega^{-1}\chi_\psi\]
 where $C'$ is a constant independently of $\epsilon$ so that $-C_w\le w\le-1$. 
 
Let $N=\nabla_X^2 w$ and $\Lambda=\lambda_1\geq \lambda_2\geq \cdots \geq \lambda_{2n}$ be the eigenvalues of $N$. Since $M$ is uniformly bounded, we have
\[\lambda_{\alpha}=\operatorname{Tr}_gN-\sum_{i\neq \alpha}\lambda_i \geq C-(2n-1)\Lambda\]
where $g$ is the Riemannian metric with respect to $\omega$. Hence $|N|\leq C(1+\max \{\Lambda,0\})$ and it suffice to show $\sup \Lambda$ is uniformly bounded.

Setting
\[
L=1+\sup_{\overline{\mathcal O}}|N|,\qquad \rho=N+Lg
\]  
we consider the following quantity
\[
Q=\log\Lambda+\xi(|\rho|^2)+\zeta(|\partial w|^2 )+e^{-Aw}
 \quad\text{where }\Lambda>0
\]
where 
\[
 \xi(r)=-\tfrac13\log(1-r/(C_nL^2)),\qquad
 \zeta(q)=-\tfrac13\log(1-q/C_P),
\]
and $A\geq 2$ is a constant to be fixed below. Since 
\[|\rho|^2\leq 4L^2 \text{\quad and \quad }|\partial_X w|_\omega+|M|_\omega\le C\quad\text{on }\overline{\mathcal O},\]
we may choose $C_n$ and $C_P$ such that both $\xi(|\rho|^2)$ and $\zeta(|\partial w|^2)$ are in $\left[0,\frac{\log 2}{3}\right]$.
It is easy to check
\begin{equation}\label{in:real-weights}
 \xi'\asymp L^{-2},\qquad \zeta'\asymp1,\qquad
 \xi''=3(\xi')^2,\qquad \zeta''=3(\zeta')^2 .
\end{equation}
The constants here are independent of $A$ and $\epsilon$. 

We first show the uniform upper bound of $Q$ implies that $\Lambda$ is uniformly bounded. Suppose that $\Lambda$ and $Q$ attain there supremum at $y$ and $z_0$. Since $\xi,\zeta$ and $e^{-Aw}$ is uniformly bounded, $Q(y)\leq Q(z_0)$ implies
\[\log \Lambda(y)\leq C_A+\log \Lambda(z_0)\]
and 
\begin{equation}\label{in:real-compare}
    L\leq C_A(1+\max \{\Lambda(y),0\})\leq C_A(1+\Lambda(z_0))
\end{equation}
where $C_A$ is a constant dependent on $A$. In particular, $N\le\nabla_X^2(\varphi_j-\varphi_0)$ on the graph $t=j$. This upper bound suffices if $Q$ attains its maximum at the boundary. Therefore, $z_0$ is an interior maximum and we may asumme $\Lambda\geq 1$ at $z_0$.

To make $\Lambda$ smooth at $z_0$, we follow the argument of Chu-McCleerey. Choose a unit top eigenvector $V_1$ of $N$ at $z_0$. Then we extend it locally by parallel transport on $X$, independently of $p$. Let
\[S=g-V_1^\flat\otimes V_1^\flat \text{ and }\widetilde N=N-S\]
where $V_1^\flat=g(V_1,\cdot)$ is its metric dual. It is obvious that the largest eigenvalue $\Lambda$ of $\widetilde{N}$
is smooth near $z_0$ and $\widetilde{N}V_1=\wLambda V_1$. Let $V_\alpha$, $2\le\alpha\le2n$,
be the remaining orthonormal eigenvectors of $\widetilde{N}$ at $z_0$, and
$\widetilde\Lambda_\alpha$ is its related eigenvalue. So we consider
\[\widetilde{Q}=\log\widetilde\Lambda+\xi(|\rho|^2)+\zeta(|\partial w|^2 )+e^{-Aw}.\]
Since $\wLambda (z_0)=\Lambda(z_0)$ and $\wLambda \leq \Lambda$ near $z_0$, $\widetilde{Q}$ attains its maximum as well at $z_0$. 

At the maximum point $z_0=(x_0,p_0)$, choose an $\omega$-unitary frame
$\{e_1,\ldots,e_n\}$ of $T^{1,0}_{x_0}X$ such that
$Me_i=m_i e_i$, with $m_1\ge\cdots\ge m_n$.
In this frame, $F^{i\bar j}=F_i\delta_{ij}$, where
$F_i=\partial F/\partial m_i$, and $\mathcal F=\sum_i F_i$.
By the ellipticity, symmetry and concavity of $F$ in Lemma \ref{ex:ideal}, we have $0<F_1\le\cdots\le F_n$. According to the maximum principle, we have
\begin{equation}\label{eqn:real-first}
    \frac{\widetilde\Lambda_\nu}{\widetilde\Lambda}
   =-\xi'R_\nu-\zeta'P_\nu+Ae^{-Aw}w_\nu
     \qquad(\nu=i\ \text{or}\ p),
\end{equation}
where $R=|\rho|^2$, $P=|\partial w|^2$ and by the ellipticity of $\mathcal L$
\begin{equation}\label{in:real-master}
\begin{split}
        0\ge\mathcal L \widetilde{Q}={}&\frac{\mathcal L\wLambda }{\wLambda }+\xi'\mathcal L R
 +\zeta'\mathcal L P-Ae^{-Aw}\mathcal L w\\
 &+\sum_iF_i(-(\wLambda _i/\wLambda )^2+\xi''|R_i|^2+\zeta''|P_i|^2+A^2e^{-Aw}|w_i|^2)\\
 &+4 \epsilon^{1/d}[-(\widetilde\Lambda_p/\widetilde\Lambda)^2+\xi''R_p^2+\zeta''P_p^2+
               A^2e^{-Aw}w_p^2]
\end{split}
\end{equation}
at $z_0$. The expansion in \eqref{in:real-master} differs from
\cite[(4.4)]{Chu-McCleerey} only by the terms arising from
$4\epsilon^{1/d}\partial_{pp}$ in $\mathcal L$. All estimates below are
taken at the maximum point $z_0$.

Now, we estimate the terms in \eqref{in:real-master}. For a unitary frame $e_i$ as above, write
$N_{1\alpha;i}=(\nabla_{e_i}N)(V_1,V_\alpha)$ and $M_{V_1}=\nabla_{V_1}M$. 
Set
\begin{align*}
 E&=\sum_{i,j}F_i(|w_{ij}|^2+|w_{i\bar j}|^2),\qquad
 G_1=2\sum_{i,\alpha>1}
       \frac{F_i|N_{1\alpha;i}|^2}{\widetilde\Lambda(\widetilde\Lambda-\widetilde\Lambda_\alpha)},\\
 G_2&=\frac{1}{\widetilde\Lambda}\sum_{i\ne j}
       \frac{F_i-F_j}{m_j-m_i}|(M_{V_1})_{i\bar j}|^2,\qquad
 G_3=\frac{c_A}{\widetilde\Lambda^2}\sum_iF_i|\nabla_{e_i}N|^2.
\end{align*}
The constant
$c_A>0$ in $G_3$ is chosen below.

\begin{lem}\label{in:real-good}
There are $c,C>0$ independent of $A$ such that
\[
 \frac{\mathcal L\widetilde\Lambda}{\widetilde\Lambda}+\xi'\mathcal L R+\zeta'\mathcal L P
 \ \ge\ G_1+G_2+G_3+cE-C\mathcal F .
\]
\end{lem}

\begin{proof}
    Fix $p=p_0$. On the slice $X\times\{p_0\}$, \eqref{equ psi} reads
\[
F(M(x,p_0))=h(x),\qquad
 h(x):=-4\epsilon^{1/d} w_{pp}(x,p_0).
\]
We set $ \mathcal L_X v=F^{i\bar j}v_{i\bar j}$. Applying the calculations in Lemma 4.2 and Lemma 4.6-4.8 in
\cite{Chu-McCleerey} to $w$, with
background form $\chi_{\varphi_0}$. Keeping the derivatives of $h$
in these calculations gives, at $(x_0,p_0)$,
\begin{align*}
 \mathcal L_X P
 &\ge E
   +2\operatorname{Re}\langle\partial_X h,\partial_X w\rangle
   -C\mathcal F,\\
 \frac{\mathcal L_X\wLambda}{\wLambda}
 &\ge G_1+G_2
   +\frac{(\nabla_X^2h)(V_1,V_1)}{\wLambda}
   -C\mathcal F,\\
 \mathcal L_X R
 &\ge 2\sum_iF_i|\nabla_{e_i}N|^2
   +2\langle\rho,\nabla_X^2h\rangle
   -CL^2\mathcal F .
\end{align*}
Here we used (4.1) in
\cite{Chu-McCleerey} for $G_2$.

Now, we add the $p$ derivatives. Since the metric $g$ is independent of $p$, we have
\begin{align*}
 P_{pp}
 &=2\operatorname{Re}
   \langle\partial_X w_{pp},\partial_X w\rangle
   +2|\partial_X w_p|^2\geq 2\operatorname{Re}
   \langle\partial_X w_{pp},\partial_X w\rangle,\\
 R_{pp}
 &=2\langle\rho,N_{pp}\rangle+2|N_p|^2\geq  2\langle\rho,N_{pp}\rangle.
\end{align*}
Then we consider the $p$ derivative of $\wLambda $. Since $\nabla_X S=0$ and $S_p=0$, we have $\nabla \widetilde{N}=\nabla N$ and $\widetilde{N}_p=N_p$. According to variation formula for the simple real eigenvalue $\wLambda $ in Chu-McCleerey,
\[\wLambda _p = \widetilde{N}_p(V_1,V_1)=N_p(V_1,V_1)\]
where $V_1$ is the eigenvector of $\widetilde{N}$ with eigenvalue $\wLambda $. Furthermore, since $S_{pp}$ also vanishes, we have
\[
\begin{split}
\wLambda _{pp}&=\widetilde{N}_{pp}(V_1,V_1)+2\sum_{\alpha=2}^{2n}\frac{|\widetilde{N}_p(V_1,V_\alpha)|^2}{\wLambda -\wLambda _\alpha}\\
&={N}_{pp}(V_1,V_1)+2\sum_{\alpha=2}^{2n}\frac{|{N}_p(V_1,V_\alpha)|^2}{\wLambda -\wLambda _\alpha}\geq {N}_{pp}(V_1,V_1)
\end{split}\]

Since
\[
 \partial_X h=-4\epsilon^{1/d}\partial_X w_{pp},
 \qquad
 \nabla_X^2 h=-4\epsilon^{1/d} N_{pp},
\]
the derivatives of $h$ cancel when we pass from $\mathcal L_X$ to
$\mathcal L=\mathcal L_X+4\epsilon^{1/d}\partial_{pp}$. We get
\[
 \mathcal L P\ge E-C\mathcal F,\qquad
 \frac{\mathcal L\wLambda }{\wLambda }\ge G_1+G_2-C\mathcal F,\qquad
 \mathcal L R\ge 2\sum_iF_i|\nabla_{e_i}N|^2-CL^2\mathcal F .
\]
Finally, \eqref{in:real-weights} gives $\zeta'\ge c$ and
$\xi'L^2\le C$. Since $\wLambda \ge1$ at $z_0$, \eqref{in:real-weights} and \eqref{in:real-compare} also give
$2\xi'\ge c_A/\wLambda ^2$. We multiply the estimates for $P$ and $R$ by $\zeta'$ and $\xi'$,
respectively, and add the estimate for $\mathcal L\wLambda /\wLambda $
to obtain the desired inequality.
\end{proof}

Then we estimate the remaining $p$-terms \[T_p=4\epsilon^{1/d}\left[-(\wLambda_p/\wLambda)^2+\xi''R_p^2+\zeta''P_p^2+A^2 e^{-Aw}w_p^2\right]\] 
in \eqref{in:real-master}

\begin{lem}\label{in:real-p}
$T_p\ge-4C_A\epsilon^{1/d}$.
\end{lem}
\begin{proof}
By \eqref{eqn:real-first} and $(a+b+c)^2\leq 3(a^2+b^2+c^2)$, we have
\[
 \left(\frac{\wLambda_p}{\wLambda}\right)^2
 =\left(-\xi'R_p-\zeta'P_p+Ae^{-Aw}w_p\right)^2
 \leq 3(\xi')^2R_p^2+3(\zeta')^2P_p^2
      +3A^2e^{-2Aw}w_p^2.
\]
Using \eqref{in:real-weights} and $|w_p|=t\leq 1$, we obtain
\begin{align*}
 T_p
 &\geq 4\epsilon^{1/d}\Bigl[
   \bigl(\xi''-3(\xi')^2\bigr)R_p^2
   +\bigl(\zeta''-3(\zeta')^2\bigr)P_p^2
   +\bigl(A^2e^{-Aw}-3A^2e^{-2Aw}\bigr)w_p^2
   \Bigr] \\
 &\geq -12\epsilon^{1/d} A^2e^{-2Aw}w_p^2
 \geq -4C_A\epsilon^{1/d} .\qedhere
\end{align*}
\end{proof}

By Lemma \ref{lem lpsi phi0} and \ref{lem math f low}, we have
\begin{equation}\label{in:real-inputs}
 -\mathcal{L} w\ge\eta\mathcal F+m_0,\qquad \mathcal F\ge c_*>0 .
\end{equation}
Substituting these inequalities and the two lemmas above in
\eqref{in:real-master}, we obtain the analogue of Lemma 4.3 in
\cite{Chu-McCleerey}:
\begin{align}
 0\ge{}&G_1+G_2+G_3-B+cE+\sum_iF_i(\xi''|R_i|^2+\zeta''|P_i|^2+A^2e^{-Aw}|w_i|^2)\notag\\
       &+Ae^{-Aw}(\eta\mathcal F+m_0)-C\mathcal F-4C_A\epsilon^{1/d} 
       \label{in:real-ledger}
\end{align}
where $B=\sum\limits_iF_i\left(\frac{\wLambda _i}{\wLambda }\right)^2=\sum\limits_iF_i \frac{|N_{11;i}|^2}{\wLambda ^2}$.

\begin{lem}\label{in:real-b}
For $0<\gamma<1/2$ and $\wLambda\ge C_A/\gamma$, we have
\begin{equation}\label{in:real-b-absorb}
 B\le G_1+G_2+G_3+U+6\gamma A^2e^{-2Aw}\sum_i F_i|w_i|^2
    +\frac{C}{\gamma}\mathcal F ,
\end{equation}
where $U=\sum\limits_iF_i(\xi''|R_i|^2+\zeta''|P_i|^2)$ and $C$ is independent of $A,\gamma,\epsilon$.
\end{lem}

\begin{proof}
    We follow the argument of \cite[Section~4.2]{Chu-McCleerey}
and include the main steps for the reader's convenience. We write $h_i=\dfrac{\wLambda_i}{\wLambda}
     =\dfrac{N_{11;i}}{\wLambda}$. Then $B=\sum\limits_i F_i |h_i|^2$. Set $\tau=A^{-2}e^{2Aw}<\frac{1}{4}$ since $A\geq 2$ and $-C_w\leq w\leq -1$. Let $i_0$ be the largest index
such that $F_{i_0}<\tau F_{i_0+1}$. If there is no such index,
take $i_0=0$. Set $I=\{i:i>i_0\}$. Since
$F_1\le\cdots\le F_n$ and $i\geq i$ if $i\in I$, we have
\begin{equation}\label{in:real-indices}
 F_i\ge\tau^{n-i}F_n\ge\tau^{n-1}F_n\ge\frac{\tau^{n-1}}n\mathcal F
 \quad(i\in I),\qquad
 F_q\le\tau F_i\quad(i\in I,\ q\notin I).
\end{equation}

As in \cite{Chu-McCleerey}, write $B=B_1+B_2+B_3$,
where
\[
 B_1=\sum_{i\notin I}F_i|h_i|^2,\qquad
 B_2=2\gamma\sum_{i\in I}F_i|h_i|^2,\qquad
 B_3=(1-2\gamma)\sum_{i\in I}F_i|h_i|^2.
\]

\noindent \textbf{Step 1. Estimate of $B_1+B_2$}

According to \eqref{eqn:real-first} and $(a+b+c)^2\leq 3(a^2+b^2+c^2)$,
\begin{align}
 |h_i|^2
 &=|-\xi'R_i-\zeta'P_i+Ae^{-Aw}w_i|^2 \notag\\
 &\le 3(\xi')^2|R_i|^2+3(\zeta')^2|P_i|^2
       +3A^2e^{-2Aw}|w_i|^2 \notag\\
 &=\xi''|R_i|^2+\zeta''|P_i|^2+3A^2e^{-2Aw}|w_i|^2.
 \label{in:real-bi}
\end{align}
Hence, by
\eqref{in:real-bi} and $2\gamma<1$,
\begin{equation}\label{in:real-b12}
\begin{split}
 B_1+B_2
 &\le \left(\sum_{i\notin I}+2\gamma\sum_{i\in I}\right)F_i(\xi''|R_i|^2+\zeta''|P_i|^2)\\
 &\quad
  +3A^2e^{-2Aw}\sum_{i\notin I}F_i|w_i|^2
  +6\gamma A^2e^{-2Aw}\sum_{i\in I}F_i|w_i|^2\\
 &\le U
  +3A^2e^{-2Aw}\tau\sum_{i\notin I}F_n|w_i|^2
  +6\gamma A^2e^{-2Aw}\sum_{i}F_i|w_i|^2 \\
  &\le U
  +6\gamma A^2e^{-2Aw}\sum_{i}F_i|w_i|^2+C\mathcal F.
\end{split}
\end{equation}

\noindent\textbf{Step 2. Control of the top eigenvector.}
     
Using $(a+b+c)^2\leq 3(a^2+b^2+c^2)$ again, \eqref{eqn:real-first} and \eqref{in:real-weights}, 
\begin{equation}\label{in:B-upper bd}
\begin{split}
    B&=\sum_iF_i|-\xi'R_i-\zeta'P_i+Ae^{-Aw}w_i|^2\\
    &\le 3\sum_iF_i\left[(\xi'R_i)^2+(\zeta'P_i)^2+A^2e^{-2Aw}|w_i|^2 \right]\\
& =\sum\limits_iF_i(\xi''|R_i|^2+\zeta''|P_i|^2)+3A^2e^{-2Aw}\sum_iF_i|w_i|^2 = U+3A^2e^{-2Aw}\sum_iF_i|w_i|^2.
\end{split}
\end{equation}

Since $|\partial w|^2<C$, we have
\[
\sum_iF_i|w_i|^2\leq C\mathcal F
\]
Substituting \eqref{in:B-upper bd} into \eqref{in:real-ledger}, we have 
\begin{equation}
    \begin{split}
      0\geq &G_1+G_2+G_3+cE+(A^2e^{-Aw}-3A^2e^{-2Aw})\sum_i F_i|w_i|^2\\
      &+Ae^{-Aw}(\eta\mathcal F+m_0)-C\mathcal F-4C_A\epsilon^{1/d}\\
      \geq&cE-C_A(\mathcal F+1)
    \end{split}
\end{equation}
since $G_i\geq 0$, $4\epsilon^{1/d}<1$ and $|w|<C$. Furthermore, $\mathcal F\geq c_*>0$ implies
\[E\leq C_A\mathcal F\]
where $E=\sum\limits_{i,j}F_i(|w_{ij}|^2)+|w_{i\bar j}|^2$.
Thus, for $i\in I$, the first inequality in \eqref{in:real-indices} gives
\begin{equation}\label{in:real-row}
 \sum_j\bigl(|w_{ij}|^2+|w_{i\bar j}|^2\bigr)
 \le\frac{E}{F_i}
 \le\frac{nC_A}{\tau^{n-1}}
 \le C_A ,
\end{equation}
where $\tau\ge A^{-2}e^{-2AC_w}$.

Write
\[
 W_1=\frac{V_1-\sqrt{-1}JV_1}{\sqrt2}
     =\sum_j\nu_j e_j,\qquad \sum_j|\nu_j|^2=1 .
\]
Since $V_1=(W_1+\overline W_1)/\sqrt2$, Cauchy--Schwarz and
\eqref{in:real-row} give, for $q\in I$,
\begin{align*}
 |N(V_1,\bar e_q)|^2=\frac12\bigl|N(W_1+\overline W_1)\bigr|^2
 &=\frac12\left|
   \sum_j\bigl(\nu_jw_{j\bar q}
              +\overline{\nu_j}\,\overline{w_{jq}}\bigr)
   \right|^2\\
 &\le\sum_j\bigl(|w_{q\bar j}|^2+|w_{qj}|^2\bigr)
 \le C_A .
\end{align*}
On the other hand,
$N(V_1,\bar e_q)=\wLambda g(V_1,\bar e_q)=\wLambda\nu_q/\sqrt2$.
Hence, as Lemma~4.5 \cite{Chu-McCleerey},
\begin{equation}\label{in:real-nu}
 |\nu_q|^2
 =\frac{2|N(V_1,\bar e_q)|^2}{\wLambda^2}
 \le\frac{C_A}{\wLambda^2},
 \qquad q\in I .
\end{equation}

\noindent{\textbf{Step 3. Estimate of $B_3$.}}

We use the decomposition in the proof of Lemma~4.6
\cite{Chu-McCleerey}. Write
$JV_1=\sum_{\alpha>1}\mu_\alpha V_\alpha$, where
$\sum_{\alpha>1}\mu_\alpha^2=1$. Since
$V_1=\sqrt2\,\overline W_1-\sqrt{-1}JV_1$, we have
\[
 N_{11;i}
 =-\sqrt{-1}\sum_{\alpha>1}\mu_\alpha N_{1\alpha;i}
   +\sqrt2\sum_q\overline{\nu_q}
(\nabla_{e_i}N)(V_1,\bar e_q).
\]
For $q\notin I$, $(M_{V_1})_{i\bar q}
=(\nabla_{V_1}\chi_{\varphi_0})_{i\bar q}
 +(\nabla_{V_1}N)(e_i,\bar e_q)$ implies
\[
 \left|(\nabla_{e_i}N)(V_1,\bar e_q)
       -(M_{V_1})_{i\bar q}\right|
 \le C\bigl(|\mathrm{Rm}|\,|d_Xw|
+|\nabla\chi_{\varphi_0}|\bigr)\le C 
\]
by commuting the derivatives in the last sum. We may therefore write
\[
 N_{11;i}=T_{1i}+T_{2i}+T_{3i}+r_i,\qquad |r_i|\le C,
\]
where
\begin{align*}
 T_{1i}&=-\sqrt{-1}\sum_{\alpha>1}\mu_\alpha N_{1\alpha;i},\\
 T_{2i}&=\sqrt2\sum_{q\notin I}\overline{\nu_q}
(M_{V_1})_{i\bar q},\\
 T_{3i}&=\sqrt2\sum_{q\in I}\overline{\nu_q}
(\nabla_{e_i}N)(V_1,\bar e_q).
\end{align*}
By Cauchy--Schwarz and \eqref{in:real-nu},
\begin{align*}
 |T_{2i}|^2
 &\le 2\left(\sum_{q\notin I}|\nu_q|^2\right)
        \sum_{q\notin I}|(M_{V_1})_{i\bar q}|^2
 \le 2\sum_{q\notin I}|(M_{V_1})_{i\bar q}|^2,\\
 |T_{3i}|^2
 &\le 2\left(\sum_{q\in I}|\nu_q|^2\right)
        \sum_{q\in I}|(\nabla_{e_i}N)(V_1,\bar e_q)|^2
 \le \frac{C_A}{\wLambda^2}|\nabla_{e_i}N|^2 .
\end{align*}
Using Young's inequality in the form
\[
 (1-2\gamma)|a+b|^2
 \le (1-\gamma)|a|^2+\gamma^{-1}|b|^2,
\]
with $a=T_{1i}$ and $b=T_{2i}+T_{3i}+r_i$, we obtain
\begin{align*}
 B_3
 &\le \frac{1-\gamma}{\wLambda^2}
       \sum_{i\in I}F_i|T_{1i}|^2
    +\frac{3}{\gamma\wLambda^2}
       \sum_{i\in I}F_i
          \bigl(|T_{2i}|^2+|T_{3i}|^2+|r_i|^2\bigr)\\
 &\le
 \underbrace{\frac{1-\gamma}{\wLambda^2}\sum_{i\in I}F_i
\left|\sum_{\alpha>1}\mu_\alpha N_{1\alpha;i}\right|^2}_{B_{31}}
 +\underbrace{\frac{C}{\gamma\wLambda^2}
       \sum_{\substack{i\in I\\q\notin I}}
  F_i|(M_{V_1})_{i\bar q}|^2}_{B_{32}}\\
 &\qquad
 +\underbrace{\frac{C_A}{\gamma\wLambda^4}
       \sum_{i\in I}F_i|\nabla_{e_i}N|^2}_{B_{33}}
 +\frac{C}{\gamma\wLambda^2}\mathcal F .
\end{align*}

For $B_{31}$, weighted Cauchy--Schwarz gives
\[
 \left|\sum_{\alpha>1}\mu_\alpha N_{1\alpha;i}\right|^2
 \le
 \left(\sum_{\alpha>1}(\wLambda-\wLambda_\alpha)\mu_\alpha^2\right)
 \sum_{\alpha>1}
\frac{|N_{1\alpha;i}|^2}{\wLambda-\wLambda_\alpha}.
\]
Since 
\[
2N(W_1,\overline W_1)
=N(V_1,V_1)+N(JV_1,JV_1)
=\wLambda+N(JV_1,JV_1)\]
and $\widetilde N=N-g$ on $V_1^\perp$,
\[
\begin{split}
    \sum_{\alpha>1}(\wLambda-\wLambda_\alpha)\mu_\alpha^2=&\wLambda \sum_{\alpha>1}\mu_\alpha^2-\widetilde N\left(\sum_{\alpha>1}\mu_\alpha V_\alpha,\sum_{\alpha>1}\mu_\alpha V_\alpha\right)=\wLambda -\widetilde N(JV_1,JV_1)
    \\
 =&\wLambda-N(JV_1,JV_1)+1
 =2\wLambda-2N(W_1,\overline W_1)+1
 \le 2\wLambda+C,
\end{split}
\]
where the last inequality uses the uniform complex Hessian bound.
Consequently,
\[
 B_{31}
 \le \frac{(1-\gamma)(2\wLambda+C)}{2\wLambda}\,G_1
 \le G_1
 \qquad\text{if }2\gamma\wLambda\ge C.
\]

For $B_{32}$, if $i\in I$ and $q\notin I$, then
\[
 F_i-F_q\ge(1-\tau)F_i\ge\tfrac12F_i,
 \qquad 0<m_q-m_i\le 2|M|\le C.
\]
Indeed, $m_q=m_i$ would imply $F_q=F_i$ by symmetry.
Thus
\[
 G_2
 \ge\frac1\wLambda
    \sum_{\substack{i\in I\\q\notin I}}
    \frac{F_i-F_q}{m_q-m_i}|(M_{V_1})_{i\bar q}|^2
 \ge\frac{c}{\wLambda}
    \sum_{\substack{i\in I\\q\notin I}}
       F_i|(M_{V_1})_{i\bar q}|^2 .
\]
It follows that
\[
 B_{32}\le\frac{C}{\gamma\wLambda}G_2
          \le G_2
 \qquad\text{if }\wLambda\ge C/\gamma .
\]

Finally, by the definition of $G_3$,
\[
 B_{33}
 \le\frac{C_A}{\gamma\wLambda^4}
    \sum_iF_i|\nabla_{e_i}N|^2
 =\frac{C_A}{c_A\gamma\wLambda^2}\,G_3
 \le G_3
 \qquad\text{if }\wLambda^2\ge\frac{C_A}{c_A\gamma}.
\]
Since $0<\gamma<1/2$, we may choose $C'_A$ large enough that
$\Lambda\ge C'_A/\gamma$ implies $\Lambda\ge1$ and there hold
$B_{31}\le G_1$, $B_{32}\le G_2$, and $B_{33}\le G_3$. Hence
\begin{equation}\label{in:real-b3}
 B_3\le G_1+G_2+G_3+\frac{C}{\gamma}\mathcal F .
\end{equation}
Adding \eqref{in:real-b12} and \eqref{in:real-b3}, we conclude that
\[
 B\le G_1+G_2+G_3+U
      +6\gamma A^2e^{-2Aw}|\partial w|_F^2
      +\frac{C}{\gamma}\mathcal F ,
\]
where $C$ is independent of $A,\gamma,\epsilon$.
\end{proof}

Now, we can complete the proof of Proposition \eqref{in:real-hessian}. 

\begin{proof}
If $\wLambda\ge C'_A/\gamma$, inserting \eqref{in:real-b-absorb} into
\eqref{in:real-ledger} gives, 
\begin{align*}
 0\ge{}&cE+Ae^{-Aw}(\eta\mathcal F+m_0)-\frac{C_0}{\gamma}\mathcal F-4C_A\epsilon^{1/d}
 \\ &+\bigl(A^2e^{-Aw}-6\gamma A^2e^{-2Aw}\bigr)\sum_i F_i |w_i|^2 ,
\end{align*}
where $C_0$ is independent of $A$. Now we fix $A\ge2$ with
$A\eta\ge6C_0+1$. Then taking $\gamma=e^{Aw}/6<\frac12$, we get $A^2e^{-Aw}-6\gamma A^2e^{-2Aw}=0$ and
\[0\ge{}cE+Ae^{-Aw}(\eta\mathcal F+m_0)-6C_0e^{-Aw}\mathcal F-4C_A\epsilon^{1/d}\geq e^{-Aw}\bigl(\mathcal F+Am_0\bigr)-4C_A\epsilon^{1/d} >0\]
for all sufficiently small $\epsilon$, a contradiction. It follows that $\wLambda<C'_A/\gamma\le6C'_Ae^{AC_w}$ at $z_0$. \eqref{in:real-compare} now bounds $\nabla_X^2w$, hence
$\nabla_X^2\psi$, throughout $\overline{\mathcal O}$. 

Finally, since $|M|\leq C$, Proposition \ref{in:mixed-bound} implies $|d_X\psi_p|_\omega^2/(-\psi_{pp})$ is bounded. Hence, $\nabla_X^2\varphi=\nabla_X^2\psi-
\dfrac{d_X\psi_p\otimes d_X\psi_p}{\psi_{pp}}$ is bounded.\qedhere
\end{proof}

\section{Existence of $\epsilon$ geodesic}
\begin{prop}\label{prop eps sol}
\label{epsilon-existence}
    There exists a constant $\epsilon_0>0$. For any $\epsilon\in (0,\epsilon_0)$, there exists a smooth solution to (\ref{introepsilongeodesicequation}).
\end{prop}
\begin{proof}
    Consider a continuity path with $0 \le \delta \le 1$:
\begin{equation}\label{exs:family}
 (\frac{\partial_t^2 \varphi_{\delta}}{4})^{d-1}P_0(H_{\varphi_{\delta}})=\eps,\qquad H_\varphi\in\Dcal,
 \qquad \varphi(\cdot,j)=\delta \varphi_j \quad(j=0,1).
\end{equation}
with $\chi$ replaced by $\chi_{\delta}=(1-\delta)k \omega+ \delta \chi$, where $k$ is a big constant such that $kI\in \Upsilon_{g_2}$. Using Proposition \ref{in:full-hessian}, there exists a constant $C$ independent of $\epsilon$ such that:
\begin{equation}\label{H bound}
    ||H_{\varphi_{\delta}}||_{\infty}\le C.
\end{equation}
For fixed $\epsilon$, we can combine the above formula with the equation $(\frac{\partial_t^2 \varphi_{\delta}}{4})^{d-1}P_0(H_{\varphi_{\delta}})=\eps$ to get a positive constant $C_{\epsilon}$ such that $\partial^2_t \varphi_{\delta}> C_{\epsilon}$ and $P_0(H_{\varphi_{\delta}})> C_{\epsilon}$. This implies that there exists a compact set $K_{\epsilon}$ in $\Dcal$ such that $H_{\varphi_{\delta}} \in K_{\epsilon}$. This implies that (\ref{exs:family}) is uniformly elliptic.
Rewrite (\ref{exs:family}) as:
\begin{equation*}
    (d-1) \log \frac{\partial_t^2 \varphi_{\delta}}{4}+ \log P_0(H_{\varphi_{\delta}})= \log \epsilon.
\end{equation*}
The left-hand side of the above equation is concave, according to Lemma \ref{ex:schur}. Using Evans-Krylov estimate \cite{krylov1984boundedly} and (\ref{H bound}), we can get that there exists a constant $C$ independent of $\delta$ and $\alpha\in (0,1)$ such that:
\begin{equation*}
||\varphi_{\delta}||_{C^{2,\alpha}}\le C.
\end{equation*}
A uniform higher-order a priori estimate for $\varphi_{\delta}$ can be obtained by a standard elliptic bootstrapping argument. This shows the closedness of the continuity path (\ref{exs:family}). The openness of (\ref{exs:family}) is similar to the proof of Lemma \ref{lem:dependence}. At $\delta=0$, an exact solution to (\ref{exs:family}) is 
\[
 \varphi^{\eps,0}(t)=2\left(\frac{\eps}{g_2(kI)}\right)^{1/d}t(t-1).
\]
In conclusion, using the method of continuity path, we can solve for the $\epsilon$-geodesic.
\end{proof}

\section{Approximation of geodesic by $\epsilon$-geodesic}

We record the uniformity of the prior PDE input needed to vary endpoints.
All spatial differential operators below act on $X$.

\begin{prop}\label{prop:pde}
When the endpoints range over a smooth compact admissible family, there
are common constants $\epsilon_0,C>0$ such that
\ref{ex:family} has a unique smooth admissible solution $\varphi_{\epsilon}$ for
$0<\epsilon\leq\epsilon_0$, with
\begin{equation}\label{eq:pde-bounds}
 \begin{aligned}
 \|\varphi_{\epsilon}\|_\infty+\|\nabla^2\varphi_\epsilon\|_\infty
 &\leq C,\qquad 0<P_{0}(H_{\varphi_{\epsilon}})\leq C\epsilon^{\frac{1}{d}},\\
 P_{0}(H_{\varphi_{\epsilon}})
 &=\frac{\varphi_{\epsilon,tt}}{4}g(A_{\varphi_\epsilon})-\left(\frac{\partial_{X}\varphi_{\epsilon,t}}{2}\right)^{*}Dg_{2}(A_{\varphi_{\epsilon}})\left(\frac{\partial_{X}\varphi_{\epsilon,t}}{2}\right)
 \end{aligned}
\end{equation}
The constants are uniform on each fixed piecewise smooth compact family
as well.
\end{prop}
\begin{proof}
Smooth solvability is
from Proposition \ref{in:real-hessian}. $L^{\infty}$ estimate for $\varphi_{\epsilon,t}$ is proved by Lemma \ref{exs:barriers}. The full Hessian estimate is proved by Proposition \ref{in:full-hessian}. All these estimates are uniform when the endpoints range over a smooth compact admissible family.

 Since $Dg>0$, admissibility and the
bounded spatial matrices imply 
\[0<P_0(H_{\varphi_{\epsilon}})=ag_{2}(A)-b^{*}Dg_{2}(A_{\varphi_{\epsilon}})b\leq a g_{2}(A)\leq Ca,\] where $a=\frac{\varphi_{\epsilon,tt}}{4}$ and $b=\frac{\partial_{X}\varphi_{\epsilon,t}}{2}$. Hence
\[
 P_0(H_{\varphi_\epsilon})^d
   =\epsilon\left(\frac{P_0(H_{\varphi_\epsilon})}{a}\right)^{d-1}
   \leq C\epsilon,
\]
which proves the defect bound. For a piecewise smooth family, use the
minimum of the finitely many existence intervals and the maximum of the
corresponding constants.
\end{proof}

\begin{lem}
\label{c1alphaconvergence}
    Let $\varphi_{j}\in \mathcal H$ for $j=0,1$. There exists a function $\varphi\in C^{1,1}(\overline{\mathcal X})$  such that for any $\alpha \in (0,1)$, we have that
    \begin{equation*}
        \lim_{\epsilon \rightarrow 0} ||\varphi_{\epsilon}-\varphi||_{C^{1,\alpha}(\overline{\mathcal X})}=0.
    \end{equation*}
    Here $\varphi_{\epsilon}$ is the $\epsilon$-geodesic connecting $\varphi_0$ and $\varphi_1$.
\end{lem}
\begin{proof}
By \eqref{eq:pde-bounds}, the first derivatives of $\varphi_\epsilon$
are uniformly bounded and uniformly Lipschitz
up to the boundary. There exists a function $\varphi\in C^{1,1}(\overline{\mathcal{X}})$ and a subsequence $\epsilon_k \to 0$ such that for all $0<\alpha<1$, 
\[\lim_{k \rightarrow \infty} ||\varphi_{\epsilon_k}-\varphi||_{C^{1,\alpha}}=0.\]

Thus every sequence $\epsilon_j\to 0$ has
a subsequence converging in $C^1$, with a $C^{1,1}$ limit. The
uniform Lipschitz bound for the first derivatives and interpolation
give convergence of this subsequence in $C^{1,\alpha}$ for every
$0<\alpha<1$.

By comparison principle, we have that $\varphi_{\epsilon_1} \ge \varphi_{\epsilon_2}$ for any $\epsilon_1 \le \epsilon_2$. This shows that the limit of different subsequence of $\varphi_{\epsilon}$ are the same. As a result, we have
    \begin{equation*}
        \lim_{\epsilon \rightarrow 0} ||\varphi_{\epsilon}-\varphi||_{C^{1,\alpha}}=0.
    \end{equation*}
\end{proof}

\begin{prop}\label{lem:limit-geodesic-equation}
The function $\varphi$, obtained by Lemma \ref{c1alphaconvergence} as a $C^{1,1}$ limit of $\epsilon$-geodesic, is a geodesic connecting $\varphi_0$ and $\varphi_1$. In other words, $\varphi$ satisfies 
 $H_\varphi\in\overline{\mathcal D_0}$ almost everywhere and
\begin{equation}\label{eq:limit-geodesic-equation}
 \sum_{k=1}^{n+1}\binom{n+1}{k}\widetilde c_k
 \bigl(\pi_X^*\chi+\sqrt{-1}D\bar D\varphi\bigr)^k
 \wedge(\pi_X^*\omega)^{n+1-k}=0
\end{equation}
as currents on $\mathcal X$, and also almost everywhere.
Here $\widetilde c_k=c_{k-1}^{(2)}$, with
$c_r^{(2)}=0$ for $r>d$.
\end{prop}

\begin{proof}
The uniform Hessian bound by Proposition \ref{prop:pde} allows us to choose $C_1$,
independently of $\epsilon$, such that
\[
 C_1\Omega+\sqrt{-1}D\bar D\varphi_{\epsilon}
 \ge \Omega>0.
\]
On a sufficiently small coordinate ball, choose a smooth real
function $\rho$ satisfying $\sqrt{-1}D\bar D\rho=\Omega$.
The functions
\[
 u_{\epsilon}=\varphi_{\epsilon}+C_1\rho
 \longrightarrow u=\varphi+C_1\rho
\]
are plurisubharmonic and converge locally uniformly.
Applying B{\l}ocki's continuity theorem
\cite[Proposition~3.3]{blocki2005weak} in complex dimension
$n+1$, with $m=n+1$, gives
\[
 \bigl(C_1\Omega+\sqrt{-1}D\bar D\varphi_{\epsilon_j}\bigr)^\ell
 \rightharpoonup
 \bigl(C_1\Omega+\sqrt{-1}D\bar D\varphi\bigr)^\ell,
 \qquad 1\le\ell\le n+1,
\]
as currents.

For each $1\le k\le n+1$, expand
\[
\begin{aligned}
 (\pi_X^* \chi + \sqrt{-1}D \bar D \varphi_{\epsilon_j})^k 
 ={}&\sum_{\ell=0}^k\binom{k}{\ell}
 \bigl(C_1\Omega+\sqrt{-1}D\bar D
       \varphi_{\epsilon_j}\bigr)^\ell \wedge
 \bigl(\pi_X^*\chi-C_1\Omega\bigr)^{k-\ell}.
\end{aligned}
\]
All remaining factors are fixed smooth forms. Thus every term
converges as a current, also after wedging with
$(\pi_X^*\omega)^{n+1-k}$.
The local limits agree on overlaps because they are limits of
the same globally defined smooth forms.

Note that 
\begin{equation*}
    \sum_{k=1}^{n+1}\binom{n+1}{k}\widetilde c_k
 \bigl(\pi_X^*\chi+\sqrt{-1}D\bar D\varphi_{\epsilon}\bigr)^k
 \wedge(\pi_X^*\omega)^{n+1-k} = P_0(H_{\varphi_{\epsilon_j}})\Omega^{n+1}.
\end{equation*}
Proposition \ref{prop:pde} makes this
form converge weakly to zero. Passing to the limit proves
\eqref{eq:limit-geodesic-equation}.

Weak-star compactness preserves the complex-Hessian bound.
Local mollification of the fixed potential $u$ then identifies
the weak wedge products with their almost-everywhere
polynomial densities, so the equation also holds almost
everywhere. Mazur's lemma and the closed convexity of
$\overline{\mathcal D_0}$ give
$H_\varphi\in\overline{\mathcal D_0}$ almost everywhere.
Finally, uniform convergence preserves the endpoint values.
\end{proof}

\section{Uniqueness of geodesic}

\begin{prop}
\label{prop:weak-geodesic-uniqueness}
Let $u,v\in C^1(\overline{\mathcal X},\mathbb R)$ be
$S^1$-invariant functions with bounded distributional complex Hessians.
Suppose that
\[
 H_u,H_v\in\overline{\Dcal},
 \qquad P_0(H_u)=P_0(H_v)=0
 \quad\text{almost everywhere in }\mathcal X,
\]
and
\[
 u(\cdot,j)=v(\cdot,j)=\varphi_j,\qquad j=0,1.
\]
Then $u=v$. In particular, the weak geodesic constructed above is
unique in this class.
\end{prop}

\begin{proof}
Let $\underline{\varphi}$ be the smooth strictly admissible barrier
from Lemma~\ref{exs:barriers}, with boundary values
$\varphi_0,\varphi_1$. For $0<\theta<1$, set
\[
 u_\theta=(1-\theta)u+\theta\underline{\varphi}.
\]
Since $\Dcal$ is open and convex,
\[
 H_{u_\theta}
 =(1-\theta)H_u+\theta H_{\underline{\varphi}}
 \in\Dcal
 \quad\text{almost everywhere}.
\]
The boundedness of $H_u$ and the uniform strict admissibility of
$H_{\underline{\varphi}}$ imply that, for fixed $\theta$, these
matrices lie in a compact subset of $\Dcal$. Hence
\[
 P_0(H_{u_\theta})\ge c_\theta>0
 \quad\text{almost everywhere}.
\]

For $0\le s\le1$, write
\[
 H_s=(1-s)H_v+sH_{u_\theta},
 \qquad
 \mathcal L_\theta f
 =\int_0^1 DP_0(H_s)
       [\sqrt{-1}D\bar D f]\,ds.
\]
The coefficients are bounded and measurable.
Moreover, $DP_0\ge0$ on $\overline{\Dcal}$, while for
$s\in[1/2,1]$ the matrices $H_s$ lie in a common compact subset
of $\Dcal$: their coefficient of
$H_{\underline{\varphi}}$ is at least $\theta/2$.
Thus $\mathcal L_\theta$ is uniformly elliptic for fixed $\theta$.
The fundamental theorem of calculus gives
\[
 \mathcal L_\theta(u_\theta-v)
 =P_0(H_{u_\theta})-P_0(H_v)
 \ge c_\theta>0
 \quad\text{almost everywhere}.
\]

In each holomorphic coordinate chart, the bounded complex Hessians
give bounded Euclidean Laplacians. Interior estimates for the
Laplacian therefore imply
\[
 u_\theta-v\in W^{2,p}_{\mathrm{loc}}(\mathcal X)
 \qquad\text{for every }1<p<\infty.
\]
Taking $p>2n+2$, we may apply the strong maximum principle to
$\mathcal L_\theta$.
Since $u_\theta-v=0$ on $\partial\mathcal X$, a positive maximum
would occur in the interior. The strong maximum principle (see \cite[Theorem 9.6]{gilbarg1977elliptic}) would
then force local constancy, contradicting the strict differential
inequality. Consequently, $u_\theta\le v$.

Letting $\theta\downarrow0$ yields $u\le v$.
Interchanging $u$ and $v$ gives the reverse inequality, and hence
$u=v$.
\end{proof}

\section{Proof of Theorem \ref{thm existence geo}}
\begin{proof}
    [Proof of Theorem \ref{thm existence geo}] 
    

Proposition~\ref{prop eps sol} gives a smooth admissible
$\epsilon$-geodesic for every sufficiently small $\epsilon>0$.
Uniqueness for fixed $\epsilon$ follows by comparison.
Proposition~\ref{prop:pde} gives the uniform $C^2$ bound, and the
comparison argument in Lemma~\ref{c1alphaconvergence} gives the
monotonicity in $\epsilon$.

By Lemma~\ref{c1alphaconvergence}, the regularized geodesics converge
in $C^{1,\alpha}$, for every $0<\alpha<1$, to a function of class
$C^{1,1}$ with the prescribed endpoints.
Proposition~\ref{lem:limit-geodesic-equation} identifies the limit
as a weak geodesic, and Proposition~\ref{prop:weak-geodesic-uniqueness}
gives uniqueness in the stated class.
\end{proof}

 \section{The Functional and the Convexity Condition}
        
       For $i=1, 2$, we define the functional $I_i: C^{\infty}(X,\mathbb{R})\rightarrow \mathbb{R}$ by
        \begin{equation}
        \label{def-I}
            I_i(\varphi) = \int_{X}\varphi\sum_{k=0}^{d}\binom{n}{k}c_k^{(i)}\sum_{j=0}^{k}\frac{1}{j+1}\binom{k}{j}\chi^{k-j}\wedge (\sqrt{-1}\partial\bar{\partial}\varphi)^{j}\wedge\omega^{n-k},
        \end{equation}
        and define $J_i=-I_i$.

        Let $\varphi(t)$ be any smooth path. We compute the derivative and see that
         \begin{equation}
         \begin{split}
            \frac{d}{dt}I_i(\varphi) &= \int_{X}\frac{d\varphi}{dt}\sum_{k=0}^{d}\binom{n}{k}c_k^{(i)}\sum_{j=0}^{k}\frac{1}{j+1}\binom{k}{j}\chi^{k-j}\wedge (\sqrt{-1}\partial\bar{\partial}\varphi)^{j}\wedge\omega^{n-k}\\
            &+\int_{X}\varphi\sum_{k=1}^{d}\binom{n}{k}c_k^{(i)}\sum_{j=1}^{k}\frac{j}{j+1}\binom{k}{j}\chi^{k-j}\wedge (\sqrt{-1}\partial\bar{\partial}\varphi)^{j-1}\wedge\omega^{n-k}\wedge\sqrt{-1}\partial\bar{\partial}\frac{d\varphi}{dt}\\
            &=\int_{X}\frac{d\varphi}{dt}\sum_{k=0}^{d}\binom{n}{k}c_k^{(i)}\sum_{j=0}^{k}(\frac{1}{j+1}+\frac{j}{j+1})\binom{k}{j}\chi^{k-j}\wedge (\sqrt{-1}\partial\bar{\partial}\varphi)^{j}\wedge\omega^{n-k}\\
            &=\int_{X}\frac{d\varphi}{dt}\sum_{k=0}^{d}\binom{n}{k}c_k^{(i)}\chi_\varphi^{k}\wedge\omega^{n-k}\\
            &=\int_{X}\frac{d\varphi}{dt}g_i(\chi_\varphi).
        \end{split}
        \label{ddtI}
        \end{equation}
        Now suppose that $\varphi_0$ and $\varphi_1$ are smooth potentials in $\mathcal{H}$, and $\varphi_{\epsilon,t}$ is the smooth $\epsilon$-geodesic connecting them. Taking one more derivative of \eqref{ddtI}, we see that
        \[
            \frac{d^2I_i(\varphi_{\epsilon,t})}{dt^2}
            =\int_X\frac{d^2\varphi_{\epsilon,t}}{dt^2}g_i(\chi_{\varphi_{\epsilon,t}})
            +\int_X\frac{d\varphi_{\epsilon,t}}{dt}\sqrt{-1}\partial\bar{\partial}\frac{d\varphi_{\epsilon,t}}{dt}\wedge g_i'(\chi_{\varphi_{\epsilon,t}}),
        \]
        where
        \[
            g_i'(\chi_{\varphi_{\epsilon,t}})=\sum_{k=1}^d kc_k^{(i)}\binom{n}{k}\chi_{\varphi_{\epsilon,t}}^{k-1}\wedge\omega^{n-k}.
        \]
        Integrating by parts, we obtain
        \[
            \frac{d^2I_i(\varphi_{\epsilon,t})}{dt^2}
            =\int_X\frac{d^2\varphi_{\epsilon,t}}{dt^2}g_i(\chi_{\varphi_{\epsilon,t}})
            -\int_X\sqrt{-1}\partial\frac{d\varphi_{\epsilon,t}}{dt}\wedge\bar{\partial}\frac{d\varphi_{\epsilon,t}}{dt}\wedge g_i'(\chi_{\varphi_{\epsilon,t}}).
        \]

        From the Lemma \ref{newgeodesiclemma}, we have
        \[
            \frac{d^2I_2(\varphi_{\epsilon,t})}{dt^2}=\int_X 4P_{0}(H_{\varphi_{\epsilon,t}})\omega^{n}.
        \]
Using the $\epsilon$-geodesic equation \ref{ex:admissibility}, we obtain
\begin{equation}
\label{secondvariationI2}
\frac{d^{2}I_{2}(\varphi_{\epsilon,t})}{dt^{2}}=4\int_{X}P_{0}(H_{\varphi_{\epsilon}})\omega^{n}=\int_{X}4\epsilon a^{1-d}\omega^{n}\ge 0.
\end{equation}
We also have the following estimate ( see Proposition \ref{prop:pde} ) that 
\begin{equation}
\label{requiredestimate}
||P_{0}(H_{\varphi_{\epsilon}})||_{L^{\infty}}\le C\epsilon^{\frac{1}{d}}
\end{equation}
Using \ref{secondvariationI2}  and \ref{requiredestimate}, we have 
\[0\le \frac{d^{2}I_{2}(\varphi_{\epsilon,t})}{dt^{2}}\le C\epsilon^{\frac{1}{d}}.\]
Using basic real analysis, we see that 
\[|I_{2}(\varphi_{\epsilon,t})-(1-t)I_{2}(\varphi_{0})-tI_{2}(\varphi_{1})|\le \frac{C\epsilon^{\frac{1}{d}}}{2}t(1-t).\]
We can use the Bedford-Taylor products \cite{Bedford-Taylor} in \eqref{def-I} to define $I_i(\varphi_t)$ and $J_i(\varphi_t)=-I_i(\varphi_t)$.
        Theorem $2.6$ of \cite{Bedford-Taylor} implies
\[
I_i(\varphi_{\epsilon,t})\longrightarrow I_i(\varphi_t),
\qquad
J_i(\varphi_{\epsilon,t})\longrightarrow J_i(\varphi_t), \qquad \epsilon\to 0.
\]
Now taking $\epsilon\rightarrow 0$, we obtain 
\[I_{2}(\varphi_{t})=(1-t)I_{2}(\varphi_{0})+tI_{2}(\varphi_{1}).\]

For the $J_{1}$-functional, we first prove the following lemma.
\begin{lem}
\label{boundlemma}
   We have $\frac{g_{1}(\chi_{\varphi_{\epsilon,t}})}{g_{2}(\chi_{\varphi_{\epsilon,t}})}\le 1$ on $\Upsilon_{g_{2}}$. 
\end{lem}
  \begin{proof}
Let \(\lambda=(\lambda_1,\ldots,\lambda_n)\) denote the eigenvalues of \(\chi_{\varphi_{\epsilon,t}}\) with respect to the background metric. Now consider
$\lambda(s)=\lambda+s(1,\ldots,1),\ s\geq0$.
Assuming this ray stays in \(\Upsilon_{g_2}\), using Theorem \ref{PositiveDerivative},  we have $q(\lambda)\leq q(\lambda(s))$, where $q(\lambda)=\frac{g_{1}(\lambda)}{g_{2}(\lambda)}$.
Since \(g_1,g_2\) have the same degree and the leading coefficient is $1$, then
$\lim_{s\to\infty}
\frac{g_1(\lambda+s\mathbf 1)}
{g_2(\lambda+s\mathbf 1)}
=1$. Consequently,
\[
\frac{g_1(\chi_{\varphi_{\epsilon,t}})}
{g_2(\chi_{\varphi_{\epsilon,t}})}
\leq1.
\]
  \end{proof}
  Let us calculate
    \begin{equation}
            \begin{split}
                 -\frac{d^2J_1(\varphi_{\epsilon,t})}{dt^2}
            &=\frac{d^2I_1(\varphi_{\epsilon,t})}{dt^2}\\
            &=\int_{X}\frac{g_{1}(\chi_{\varphi_{\epsilon,t}})}{g_{2}(\chi_{\varphi_{\epsilon,t}})}\left(\frac{d^{2}\varphi_{\epsilon,t}}{dt^{2}}g_{2}(\chi_{\varphi_{\epsilon,t}})-\sqrt{-1}\partial\frac{d\varphi_{\epsilon,t}}{dt}\wedge\bar{\partial}(\frac{d\varphi_{\epsilon,t}}{dt})\wedge g_{2}'(\chi_{\varphi_{\epsilon,t}})\right)\\
    &+\int_{X}\sqrt{-1}\partial\frac{d\varphi_{\epsilon,t}}{dt}\wedge\bar{\partial}(\frac{d\varphi_{\epsilon,t}}{dt})\wedge (\frac{g_{1}(\chi_{\varphi_{\epsilon,t}})}{g_{2}(\chi_{\varphi_{\epsilon,t}})}g_{2}'(\chi_{\varphi_{\epsilon,t}})-g_{1}'(\chi_{\varphi_{\epsilon,t}})).
            \end{split}
        \end{equation}
       Using the proper position ( see Theorem \ref{PositiveDerivative} ), we see that the second term is non-positive. Hence we obtain 
       \begin{equation}
           \begin{split}
               &-\frac{d^2J_1(\varphi_{\epsilon,t})}{dt^2}\\
               &\le \int_{X}\frac{g_{1}(\chi_{\varphi_{\epsilon,t}})}{g_{2}(\chi_{\varphi_{\epsilon,t}})}\left(\frac{d^{2}\varphi_{\epsilon,t}}{dt^{2}}g_{2}(\chi_{\varphi_{\epsilon,t}})-\sqrt{-1}\partial\frac{d\varphi_{\epsilon,t}}{dt}\wedge\bar{\partial}(\frac{d\varphi_{\epsilon,t}}{dt})\wedge g_{2}'(\chi_{\varphi_{\epsilon,t}})\right).
           \end{split}
       \end{equation}
       From Lemma \ref{newgeodesiclemma}, we see that 
       \[\frac{d^{2}\varphi_{\epsilon,t}}{dt^{2}}g_{2}(\chi_{\varphi_{\epsilon,t}})-\sqrt{-1}\partial\frac{d\varphi_{\epsilon,t}}{dt}\wedge\bar{\partial}(\frac{d\varphi_{\epsilon,t}}{dt})\wedge g_{2}'(\chi_{\varphi_{\epsilon,t}})=4P_{0}(H_{\varphi_{\epsilon,t}})\omega^{n}>0.\]
       Using Lemma \ref{boundlemma} and \ref{prop:pde}, we have 
       \[-\frac{d^2J_1(\varphi_{\epsilon,t})}{dt^2}\le \int_{X}4P_{0}(H_{\varphi_{\epsilon,t}})\omega^{n}\le C\epsilon^{\frac{1}{d}}. \]
       Now let us take $0\le t_{1}<t_{2}<t_{3}\le 1$ and define $$L(t)=
\begin{cases}
\dfrac{t-t_1}{t_2-t_1}, & t_1\leq t\leq t_2, \\[6pt]
\dfrac{t_3-t}{t_3-t_2}, & t_2\leq t\leq t_3.
\end{cases}$$ We see that $L(t_{1})=L(t_{3})=0$ and $L(t_{2})=1$. First let us consider the function $$K:= \frac{J_{1}(\varphi_{\epsilon,t_{2}})-J_{1}(\varphi_{\epsilon,t_{1}})}{t_{2}-t_{1}}-\frac{J_{1}(\varphi_{\epsilon,t_{3}})-J_{1}(\varphi_{\epsilon,t_{2}})}{t_{3}-t_{2}}.$$ By the fundamental theorem of calculus, $$K=\int_{t_{1}}^{t_{3}}L'(t)\frac{dJ_{1}(\varphi_{\epsilon,t})}{dt}dt.$$ Now integrating by parts we obtain $$K=[L(t)\frac{dJ_{1}(\varphi_{\epsilon,t})}{dt}]|_{t_{1}}^{t_{3}}-\int_{t_{1}}^{t_{3}}L(t)\frac{d^{2}J_{1}(\varphi_{\epsilon,t})}{dt^{2}}dt=-\int_{t_{1}}^{t_{3}}L(t)\frac{d^{2}J_{1}(\varphi_{\epsilon,t})}{dt^{2}}dt$$ since $L(t_{1})=L(t_{3})=0$. Because $L(t)\ge 0$ on $[t_{1},t_{3}]$ and $\frac{d^{2}J_{1}(\varphi_{\epsilon,t})}{dt^{2}}\ge -C\epsilon^{\frac{1}{d}},$ we obtain $$K\le C\epsilon^{\frac{1}{d}}\int_{t_{1}}^{t_{3}}L(t)dt.$$ Now integrating and taking $\epsilon\rightarrow 0$, we obtain $$\frac{J_{1}(\varphi_{t_{2}})-J_{1}(\varphi_{t_{1}})}{t_{2}-t_{1}}-\frac{J_{1}(\varphi_{t_{3}})-J_{1}(\varphi_{t_{2}})}{t_{3}-t_{2}}\le 0.$$ Rearranging we obtain $$J_{1}(\varphi_{t_{2}})\le \frac{t_{3}-t_{2}}{t_{3}-t_{1}}J_{1}(\varphi_{t_{1}})+\frac{t_{2}-t_{1}}{t_{3}-t_{1}}J_{1}(\varphi_{t_{3}}).$$
In other words, the $I_2$-functional is affine and the $J_1$-functional is convex on weak geodesics.
Now  using \ref{ddtI}, $|\frac{d\varphi_{\epsilon,t}}{dt}|\le C$ by Lemma \ref{exs:barriers}, \ref{in:full-hessian}, we see that 
\[|\frac{dI_{i}(\varphi_{\epsilon,t})}{dt}|\le C.\]
Taking $\epsilon\rightarrow 0$, we see that $J_{1}$ and $I_{2}$ are Lipschitz.

\section{$\mathcal{E}^{p}$\  Theory}
We use the structural facts $\Upsilon_{g_{2}}$ is convex, $g_{2}>0$ and $Dg_{2}>0$ on $\Upsilon_{g_{2}}$,
and the block domain is convex with a strictly elliptic logarithmic
operator.\\
For $\varphi\in\mathcal H$, set
\begin{equation}\label{eq:measure}
 g_{2}(\chi_{\varphi})=g_{2}(A_\varphi)\omega^n
 =\sum_{j=0}^d c_j^{(2)}\binom nj
       (\chi+\ddc\varphi)^j\wedge\omega^{n-j},
 \qquad V=\int_X g_{2}(\chi_{\varphi})>0.
\end{equation}
Closedness of $\chi$ and $\omega$ and Stokes' theorem show that $V$ is
independent of $\varphi$. Let $\varphi(t),\ t\in[0,1]$, be a piecewise smooth path in $\mathcal{H}$. For $1\leq p<\infty$, define
\begin{equation}\label{eq:path-distance}
 \begin{gathered}
 \|\xi\|_{p,\varphi}
   =\left(\int_X|\xi|^p g_{2}(\chi_{\varphi})\right)^{1/p}, \qquad E_{p}(\varphi,t)=\int_{X}\big|\frac{d\varphi}{dt}\big|^{p}g_{2}(\chi_{\varphi})\\
 \ell_p(\varphi)=\int_0^1E_{p}(\varphi,t)^{\frac{1}{p}}dt, \qquad
 d_p(\varphi_0,\varphi_1)=\inf_\psi\ell_p(\psi),
 \end{gathered}
\end{equation}
where the infimum is over piecewise smooth $\mathcal H$-valued paths with
the prescribed endpoints.
\subsection{Approximation and variation}

\begin{lem}\label{lem:dependence}
Let $\psi(s)$ be a smooth admissible path and let
$\widehat\varphi\in\mathcal H$ be fixed. For each fixed positive
$\epsilon$ in the common existence interval, the solutions
$\varphi_\epsilon(s,t)$ from $\psi(s)$ to $\widehat\varphi$ depend smoothly on
$s$ and satisfy
\begin{equation}\label{eq:parameter-bound}
 \|\varphi_{\epsilon,s}(s,\cdot)\|_\infty
   \leq\|\psi_s(s)\|_\infty.
\end{equation}
This includes parameters for which $\psi(s)=\widehat\varphi$.
\end{lem}
\begin{proof}
The operator
$\widetilde{G}(H)=\log P_0(H)+(d-1)\log a$ is smooth and strictly elliptic
on $\mathcal D_0$ by Lemma \ref{ex:schur}.     For fixed \(\epsilon>0\), the equation is
\[\widetilde{G}(H_{\varphi_\epsilon})= \log\epsilon.\]
The matrix \(H_{\varphi_\epsilon}\) is in the admissible domain \(\mathcal D_0\). Since \(X\times\bar{\mathcal{A}}\) is compact and the solution is smooth, its matrix range is a compact subset of \(\mathcal D_0\).

By strict ellipticity of \(\widetilde{G}\), its linearization at the solution is a strictly elliptic operator. More precisely, the linearization has the form
\begin{align*}
    &\mathcal L_\epsilon:
C^{2,\alpha}_0(\overline{\mathcal{X}})
\longrightarrow C^{\alpha}(\overline{\mathcal{X}}),\\
&\mathcal L_\epsilon v
=\tilde{G}^{\alpha\bar\beta}(H_{\varphi_\epsilon})
\,\partial_\alpha\partial_{\bar\beta}v,
\end{align*}
with no zeroth-order term and the subscript $0$ denotes zero boundary values.
Using maximum principle, we see that the Dirichlet operator is injective and Schauder theory and the Fredholm alternative give surjectivity. Hence it is invertible between the corresponding H\"older spaces.
 After subtracting the endpoint interpolation, the implicit function theorem applies in
zero-boundary H\"older spaces. This proves the local smooth dependence on $s$. Uniqueness of the solution patches the local parameter-dependent solutions into a smooth family along the entire parameter interval.

Differentiate the nonlinear equation with respect to \(s\) and applying the chain rule, we get
\[
D\widetilde{G}\big|_{H_{\varphi_\epsilon}}
\left(H_{\varphi_{\epsilon,s}}\right)=0.
\]
Since \(H_\varphi\) depends linearly on the second derivatives of \(\varphi\), this becomes
\[
\mathcal L_\epsilon \varphi_{\epsilon,s}=0.
\]
Because of the boundary values, we have
\[
\varphi_{\epsilon,s}(s,0)=\psi_s(s),
\qquad
\varphi_{\epsilon,s}(s,1)=0.
\]
The maximum principle gives
\[\|\varphi_{\epsilon,s}(s,\cdot)\|_{\infty}
\le \|\psi_s(s)\|_{\infty}.\]
\end{proof}

For a two-parameter family, write $\varphi=\varphi_\epsilon(s,t)$ while keeping
$\epsilon$ fixed and for $h\in C^2(\mathbb R)$, set
\[
 E_h(s,t)=\int_Xh(\varphi_t)g_{2}(\chi_{\varphi}),\qquad
 J_h(s,t)=\int_Xh'(\varphi_t)\varphi_s\ g_{2}(\chi_{\varphi}).
\]
\begin{lem}\label{lem:variation}
The two energy variations satisfy
\begin{equation}\label{eq:variation}
 \begin{aligned}
 \partial_tE_h&=\int_Xh'(\varphi_t)4P_{0}(H_{\varphi_{\epsilon}})\omega^n,\\
 \partial_sE_h&=\partial_tJ_h
       -\int_Xh''(\varphi_t)\varphi_s4P_{0}(H_{\varphi_{\epsilon}})\omega^n.
 \end{aligned}
\end{equation}
\end{lem}
\begin{proof}
We calculate
\[
 \begin{aligned}
 \partial_tE_h
 &=\int_Xh'(\varphi_t)\varphi_{tt}g_{2}(\chi_{\varphi})
       +\int_Xh(\varphi_t)\ddc \varphi_t\wedge g_{2}'(\chi_{\varphi})\\
 &=\int_Xh'(\varphi_t)
       \bigl(\varphi_{tt}g_{2}(\chi_{\varphi})-\sqrt{-1}\partial \varphi_t\wedge\bar\partial \varphi_t\wedge g_{2}'(\chi_{\varphi})\bigr),
 \end{aligned}
\]
which proves the first identity by the Lemma \ref{newgeodesiclemma}.
Let us calculate
\[
 \partial_sE_h
   =\int_Xh'(\varphi_t)\varphi_{st}g_{2}(\chi_{\varphi})
      +\int_X\varphi_s\ddc h(\varphi_t)\wedge g_{2}'(\chi_{\varphi})
\]
and 
\[\partial_{t}J_{h}=\int_{X}h''(\varphi_{t})\varphi_{tt}\varphi_{s}g_{2}(\chi_{\varphi})+\int_{X}h'(\varphi_{t})\varphi_{ts}g_{2}(\chi_{\varphi})+\int_{X}h'(\varphi_{t})\varphi_{s}\sqrt{-1}\partial\bar{\partial}\varphi_{t}\wedge g_{2}'(\chi_{\varphi}).\]
Subtracting, we get
\[
 \partial_tJ_h-\partial_sE_h
  =\int_Xh''(\varphi_t)\varphi_s
       \bigl(\varphi_{tt}g_{2}(\chi_{\varphi})-i\partial \varphi_t\wedge\bar\partial \varphi_t\wedge g_{2}'(\chi_{\varphi})\bigr).
\]
Lemma \ref{newgeodesiclemma} proves the second identity.
\end{proof}

\subsection{Regularized length and endpoint comparison}
Fix $1\leq p<\infty$. For $0<\delta\leq1$, use the smooth function
\begin{equation}\label{eq:regularization}
 \begin{gathered}
 h_\delta(r)=(r^{2}+\delta^{2})^{\frac{p}{2}} ,\qquad
 E_{p,\delta}(\varphi,t)=\int_Xh_\delta(\frac{d\varphi}{dt})g_{2}(\chi_{\varphi}),\\
 L_\delta(\varphi)=\int_0^1E_{p,\delta}(\varphi,t)^{\frac{1}{p}}\ dt.
 \end{gathered}
\end{equation}
Let $J_\delta=J_{h_\delta}$ in Lemma~\ref{lem:variation}.
\begin{lem}
\label{lem:regularization}
Uniformly on a fixed compact endpoint family,
\begin{equation}\label{eq:length-regularization}
 0\leq L_\delta(\varphi_\epsilon)-\ell_p(\varphi_\epsilon)
       \leq V^{\frac{1}{p}}\delta,
\end{equation}
and
\begin{equation}\label{eq:energy-oscillation}
 \operatorname*{osc}_{t\in[0,1]}E_p(\varphi_\epsilon,t)
       \leq C\epsilon^{\frac{1}{d}}.
\end{equation}
\end{lem}
\begin{proof}
The pointwise inequalities
$|r|\leq\sqrt{r^{2}+\delta^2}\leq|r|+\delta$ 
and Minkowski's inequality give \eqref{eq:length-regularization}. We will also use
\begin{equation}\label{eq:regularization-floor}
 E_{p,\delta}\geq V\delta^{p},\qquad
 \frac{|h_\delta'(r)|}{p}\leq h_\delta(r)^{1-\frac{1}{p}}.
\end{equation}
On bounded $r$ intervals, $|h_\delta'(r)|$ is bounded independently of
$0<\delta\leq1$, including when $p=1$. Proposition~\ref{prop:pde} and
the first identity in \eqref{eq:variation} therefore imply
$|\partial_tE_{p,\delta}|\leq C\epsilon^{1/d}$, with $C$ independent of
$\delta$. As $\delta\to0$, $h_\delta(r)\to|r|^p$ uniformly on these
bounded intervals. The fixed mass $V$ makes $E_{p,\delta}\to E_p$ uniform
in $t$, so the same oscillation bound holds for $E_p$.
\end{proof}
\begin{prop}
\label{prop:comparison}
Let $\psi:[0,1]\to\mathcal H$ be piecewise smooth and
$\widehat\varphi\in\mathcal H$ be fixed. Let $\varphi_\epsilon(s,t)$ solve
\ref{ex:family} from $\psi(s)$ to $\widehat\varphi$. For every fixed
$0<\delta\leq1$, there is $C_\delta$, independent of $\epsilon$, such that
\begin{equation}\label{eq:comparison}
 L_\delta(\varphi_\epsilon(0,\cdot))
  \leq\ell_p(\psi)+L_\delta(\varphi_\epsilon(1,\cdot))
           +C_\delta\epsilon^{1/d}.
\end{equation}
Note that $\widehat{\varphi}$ may be on the path $\psi$.
\end{prop}
\begin{proof}
We first work on a smooth piece and suppress only the index $\epsilon$. Because of Lemma \ref{lem:dependence}, we can differentiate in $s$. Using 
Lemma~\ref{lem:variation} and integration by parts in $t$, we obtain
\begin{equation}\label{eq:length-variation}
 \begin{aligned}
 \partial_s L_\delta(\varphi(s,\cdot))
 &=\frac1p\bigl[E_{p,\delta}^{1/p-1}J_\delta\bigr]_{t=0}^{t=1}\\
 &\quad-\frac1p\left(\frac1p-1\right)
    \int_0^1 E_{p,\delta}^{\frac{1}{p}-2}(\partial_tE_{p,\delta})J_\delta\,dt\\
 &\quad-\frac1p\int_0^1 E_{p,\delta}^{1/p-1}
       \int_Xh_\delta''(\varphi_t)\varphi_s4P_{0}(H_{\varphi_{\epsilon}})\omega^n\,dt.
 \end{aligned}
\end{equation}
In the above, we have used the notation $J_{\delta}=J_{h_{\delta}}$.  For fixed $\delta>0$, the first inequality in \eqref{eq:regularization-floor}, the
bounds for $\varphi_s$ and $\varphi_t$, boundedness of $h_\delta''$ on bounded
intervals, and $|\partial_tE_{p,\delta}|\leq C\epsilon^{1/d}$ ( see Lemma \ref{lem:regularization} ) bound both integral terms by $C_\delta\epsilon^{1/d}$ in
absolute value. Observe that $J_\delta$ is bounded using $V$ and
\eqref{eq:parameter-bound}. If $p=1$, the first integral term vanishes.

At $t=1$, $\varphi_s=0$. At $t=0$, we have $\varphi_s=\psi_s$ and
$g_{2}(\chi_\varphi)=g_{2}(\chi_\psi)$. H\"older's inequality and
\eqref{eq:regularization-floor} give, for $p>1$,
\begin{align*}
    &\frac1p E_{p,\delta}(s,0)^{\frac{1}{p}-1}|J_\delta(s,0)|\\
    &\leq\frac{1}{p} E_{p,\delta}(s,0)^{\frac{1}{p}-1}\times 
\left(\int_X|\psi_s|^p g_2(\chi_{\psi(s)})\right)^{1/p}\times
\left(
\int_X |h_{\delta}'(\varphi_{\epsilon,t}(s,0))|^{p/(p-1)}
g_2(\chi_{\psi(s)})
\right)^{(p-1)/p}\\
&\leq \frac{1}{p} E_{p,\delta}(s,0)^{\frac{1}{p}-1}\times 
\left(\int_X|\psi_s|^p g_2(\chi_{\psi(s)})\right)^{1/p}\times
p(E_{p,\delta}(s,0))^{\frac{p-1}{p}}\\
      & \leq\|\psi_s\|_{p,\psi(s)}.
\end{align*}
For $p=1$, the same bound follows directly from $|h_\delta'|\leq1$.
Consequently
\[
 \partial_sL_\delta(\varphi(s,\cdot))
       \geq-\|\psi_s\|_{p,\psi(s)}-C_\delta\epsilon^{1/d}.
\]
Integrate on each smooth piece and add. The regularized lengths agree at
the break points by uniqueness. This proves \eqref{eq:comparison}.
\end{proof}

\begin{lem}\label{lem:diagonal}
Suppose both endpoints equal $\varphi\in\mathcal H$, and let
$m=\min_Xg_{2}(A_\varphi)>0$. Then
\begin{equation}\label{eq:diagonal}
 \|\varphi_{\epsilon,t}\|_\infty\leq2(\epsilon/m)^{1/d},\qquad
 L_\delta(\varphi_\epsilon)
       \leq V^{1/p}\bigl(2(\epsilon/m)^{1/d}+\delta\bigr).
\end{equation}
\end{lem}
\begin{proof}
The function
\[
 \underline \varphi(x,t)=\varphi(x)+2(\epsilon/m)^{1/d}t(t-1)
\]
has blocks $a=(\epsilon/m)^{1/d}$, $b=0$ and $A=A_\varphi$. It is
admissible and satisfies
\[
 a^{d-1}P_0(H_{\underline \varphi})
       =(\epsilon/m)g_{2}(A_\varphi)\geq\epsilon.
\]
Comparison gives $\underline \varphi\leq \varphi_\epsilon$. Indeed, the
linearization of $\log P_0+(d-1)\log a$ in the convex domain
$\mathcal D_0$ for each fixed positive $\epsilon$ is
strictly elliptic. Time convexity also gives
$\varphi_\epsilon\leq\varphi$. The endpoint derivatives of these inequalities,
followed by monotonicity of $\varphi_{\epsilon,t}$, give
\[
 -2(\epsilon/m)^{1/d}
 \leq \varphi_{\epsilon,t}(x,0)
 \leq \varphi_{\epsilon,t}(x,t)
 \leq \varphi_{\epsilon,t}(x,1)
 \leq2(\epsilon/m)^{1/d}.
\]
This proves the first bound in \eqref{eq:diagonal}. The second follows
from $\sqrt{|\varphi_{\epsilon,t}|^{2}+\delta^2}\leq|\varphi_{\epsilon,t}|+\delta$ and
the fixed mass $V$.
\end{proof}

\subsection{Length minimization and the metric}
\begin{proof}[Proof of Theorem~\ref{thm:metric}]
We know that every $\varphi_\epsilon$ is an admissible smooth path, so
$d_p(\varphi_0,\varphi_1)\leq\ell_p(\varphi_\epsilon)$. Fix any piecewise
smooth path $\psi$ from $\varphi_0$ to $\varphi_1$, and apply
Proposition~\ref{prop:comparison} with
$\widehat\varphi=\varphi_1$. The terminal family member has identical
endpoints. Lemmas~\ref{lem:regularization} and \ref{lem:diagonal} give,
with $\delta>0$ fixed,
\[
 \limsup_{\epsilon\to0}\ell_p(\varphi_\epsilon)
       \leq\ell_p(\psi)+V^{\frac{1}{p}}\delta.
\]
Let $\delta\to0$ and take the infimum over $\psi$. This proves
\eqref{eq:length-limit}.

Since $\ell_p(\varphi_\epsilon)^p$ lies between the minimum and maximum of
$E_p(\varphi_\epsilon,t)$, \eqref{eq:energy-oscillation} implies
\[
 \sup_{t\in[0,1]}
   |E_p(\varphi_\epsilon,t)-\ell_p(\varphi_\epsilon)^p|
       \leq C\epsilon^{1/d}.
\]
Together with \eqref{eq:length-limit}, this proves
\eqref{eq:energy-limit}.

Pointwise time convexity, the prescribed endpoints and the fundamental theorem of calculus imply
\[
 \varphi_{\epsilon,t}(x,0)\leq\varphi_1(x)-\varphi_0(x)
                   \leq \varphi_{\epsilon,t}(x,1).
\]
On $\{\varphi_0>\varphi_1\}$, the first inequality gives
$|\varphi_{\epsilon,t}(x,0)|\ge |\varphi_1(x)-\varphi_0(x)|$. On $\{\varphi_1>\varphi_0\}$, the
last inequality gives $|\varphi_{\epsilon,t}(x,1)|\ge |\varphi_1(x)-\varphi_0(x)|$. Now using 
\eqref{eq:energy-limit}, we obtain 
\[d_{p}(\varphi_{0},\varphi_{1})^{p}=\lim_{\epsilon\rightarrow 0} \sup_{t\in[0,1]} E_{p}(\varphi_{\epsilon},t)\ge \max\left\{\begin{aligned}&
  \int_{\{\varphi_0>\varphi_1\}}
       |\varphi_0-\varphi_1|^p g_{2}(\chi_{\varphi_0}),\\&
  \int_{\{\varphi_1>\varphi_0\}}
       |\varphi_1-\varphi_0|^p g_{2}(\chi_{\varphi_1})
       \end{aligned}\right\}. \]
       This proves \eqref{eq:endpoint-bound}.
\end{proof}

\begin{proof}[Proof of Corollary~\ref{cor:metric}]
 The triangle inequality, symmetry and positivity are standard. \\Since $\mathcal{H}$ is path connected, the reversed path proves symmetry and the glued path proves the triangle inequality.
The only non-trivial thing is that $d_{p}(\varphi_{0},\varphi_{1})=0$ implies that $\varphi_{0}=\varphi_{1}$. This is true by \eqref{eq:endpoint-bound}.
\end{proof}

\bibliographystyle{plain}
\bibliography{reference1}

\begin{center}
\begin{minipage}{0.88\textwidth}
\small
\textbf{Kartick Ghosh}\\[0.35em]
University of Science and Technology of China\\
No.\ 96 Jinzhai Road, Baohe District, Hefei, Anhui 230000, P.\,R.\ China\\[0.35em]
\href{mailto:kghosh@ustc.edu.cn}{\texttt{kghosh@ustc.edu.cn}}\\
\textbf{XiShen Jin}\\[0.35em]
Renmin University of China\\
 No.59 Zhongguancun Street, Haidian District, Beijing 100872, China\\[0.35em]
\href{mailto:jinxishen@ruc.edu.cn}{\texttt{jinxishen@ruc.edu.cn}}\\
\textbf{Ziyue Wang}\\[0.35em]
University of Science and Technology of China\\
No.\ 96 Jinzhai Road, Baohe District, Hefei, Anhui 230000, P.\,R.\ China\\[0.35em]
\href{mailto:wangziyue@mail.ustc.edu.cn}{\texttt{wangziyue@mail.ustc.edu.cn}}\\
\textbf{Yulun Xu}\\[0.35em]
University of Toronto\\
 40 St. George Street, Toronto, ON, Canada\\[0.35em]
\href{mailto:yulun.1.xu@gmail.com}{\texttt{yulun.1.xu@gmail.com}}
\end{minipage}
\end{center}

\end{document}